\documentclass[utf8]{article}

\usepackage[left=2.5cm,right=2.5cm,top=2.5cm,bottom=2.5cm]{geometry}
\usepackage{float} 
\usepackage{threeparttable}
\usepackage{enumerate} 
    \usepackage{mathtools,amssymb,amsthm}
    \usepackage{natbib} 

    \allowdisplaybreaks 

    \newcommand{\oto}{\leftrightarrow}
    
    \newcommand{\tto}{\Rightarrow}
    
    \newcommand{\otto}{\Leftrightarrow}

    \newtheorem{definition}{Definition}[subsection] 
    \newtheorem{lemma}{Lemma}[subsection]
    \newtheorem{theorem}{Theorem}[subsection]
    \newtheorem{metatheorem}{Meta Theorem}[subsection]
    \newtheorem{corollary}{Corollary}[subsection]
    \newtheorem{remark}{Remark}[subsection]
    \newtheorem{example}{Example}[subsection]
    \newtheorem{derivedrule}{Rule}[subsection]
    \newtheorem{metarule}{Meta Rule}[subsection]

    \title{Defining implication relation for classical logic}
    \author{
        Li Fu \\
        \small
        School of Software, Chongqing University, Chongqing, China \\
        \small
        fuli@cqu.edu.cn
    }
    \date{} 

\begin{document}
        \maketitle
        \begin{abstract}
            In classical logic, ``P implies Q'' is equivalent to ``not-P or Q''. It is well known that the equivalence is problematic. Actually, from ``P implies Q'', ``not-P or Q'' can be inferred (``Implication-to-Disjunction'' is valid), whereas from ``not-P or Q'', ``P implies Q'' cannot be inferred in general (``Disjunction-to-Implication'' is not generally valid), so the equivalence between them is invalid in general. This work aims to remove the incorrect Disjunction-to-Implication from classical logic (CL). The logical system (the logic IRL) this paper proposes has the expected properties: (a) CL is obtained by adding Disjunction-to-Implication to IRL, and (b) Disjunction-to-Implication is not derivable in IRL; while (c) fundamental laws in classical logic, including law of excluded middle (LEM) and principle of double negation, law of non-contradiction (LNC) and ex contradictione quodlibet (ECQ), conjunction elimination and disjunction introduction, and hypothetical syllogism and disjunctive syllogism, are all retained in IRL.

        \end{abstract}


        \section{Introduction}
        \label{sec:Introduction}

        In classical logic, ``$P$ implies $Q$'' ($P\to Q$) is taken to be logically equivalent to ``it cannot be $P$ and not $Q$'' ($\neg(P\wedge\neg Q)$) or ``not-$P$ or $Q$'' ($\neg P\vee Q$) by duality. This equivalence can be viewed as the definition of implication by disjunction and negation ($P\to Q=\neg P\vee Q$, called \emph{material implication}), or used to define disjunction by implication and negation, so $P\to Q$ and $\neg P\vee Q$ can be replaced with one another in a derivation.

        It is well known that this definition of implication is problematic in that it leads to counter-intuitive results or ``paradoxes''. There are many paradoxes of material implication \citep[see e.g.][]{bronstein}. For instance, the theorem
        $$(P\to Q)\vee(Q\to P)$$
        of classical logic means that any two propositions must be related by implication (total order). This does not agree with a large number of situations in the real world. For example, let $P=$ ``Alice is in Athens'' and $Q=$ ``Alice is in London'' (given an instant in the future), then both $P\to Q$ and $Q\to P$ cannot be true, so $(P\to Q)\vee(Q\to P)$ cannot be true.

        Many efforts have been made to resolve this problem and important non-classical logics are developed, such as modal logics which introduce the ``strict implication'', relevance logics which require that the antecedent and consequent of implications to be relevantly related (so they are paraconsistent logics which reject \emph{ex contradictione quodlibet (ECQ)}).

        However, classical logic exists on its own reason \citep[see e.g.][]{fulda}. And implication is the kernel concept in logic: ``implies'' means the same as that ``of ordinary inference and proof'' \citep{lewis}. So it is not satisfactory that the uniquely fundamental, simple and useful classical logic has a wrong definition of implication -- the material implication.

        The paper explains why ``Implication-to-Disjunction'' ($(P\to Q)\to \neg P\vee Q$) is valid whereas ``Disjunction-to-Implication'' ($\neg P\vee Q\to(P\to Q)$) is not, by describing firstly the original meaning of ``implication''. So the goal of this work is a logic in which Disjunction-to-Implication is removed while fundamental laws in classical logic such as law of excluded middle (LEM) and principle of double negation, law of non-contradiction (LNC) and ex contradictione quodlibet (ECQ), conjunction elimination and disjunction introduction, and hypothetical syllogism (transitivity of implication) and disjunctive syllogism, are all retained.


        \section{Original meaning of implication}

        \subsection{Four situations}
        \label{subsec:Four situations}

        Let $P$ and $Q$ be propositions. Consider how $P$ affects $Q$. There are four possible situations about $Q$ when $P$ is true:
        \begin{enumerate}
            \item when $P$ is true, $Q$ is \emph{necessarily true};
            \item when $P$ is true, $Q$ is \emph{coincidentally true};
            \item when $P$ is true, $Q$ is \emph{necessarily false};
            \item when $P$ is true, $Q$ is \emph{coincidentally false}.
        \end{enumerate}

        Where:
        \begin{itemize}
            \item ``necessarily'' means ``with a certain mechanism'';
            \item  ``coincidentally'' means ``without a certain mechanism''.
        \end{itemize}

        And  ``certain mechanism'' means one of the two:
        \begin{enumerate}
            \item It is certain literally or physically (in some scope, cf. Subsection \ref{subsec:Truth values of implication}), for example, ``in the open air'' is a physical mechanism of ``if it rains, then the ground gets wet'' (cf. Example \ref{exam:OPQ});
            \item It is certain formally or logically, for example, ``from $(P\vee Q)\wedge\neg P$ to $Q$'' is a logical mechanism.
        \end{enumerate}

        \subsection{Original meaning of implication}
        \label{subsec:OriginalMeaningOfImplication}

        The original meaning of implication $P\to Q$ is described as follows.
        \begin{itemize}
            \item $P\to Q$ is true, if and only if whenever $P$ is true, $Q$ is \emph{necessarily true}, no matter whether $Q$ is true or false when $P$ is false. Equivalently,
            \item $P\to Q$ is false, if and only if whenever $P$ is true, $Q$ is \emph{not necessarily true} (i.e. $Q$ is \emph{coincidentally true}, \emph{necessarily false}, or \emph{coincidentally false}), no matter whether $Q$ is true or false when $P$ is false.
        \end{itemize}

        According to this description, implication from propositions $P$ to $Q$ is (can be defined as) a relation determined by \emph{a certain mechanism} (physical or logical) from $P$ to $Q$, which guarantees that $Q$ is true whenever $P$ is true. Without such a mechanism, even if $Q$ is true when $P$ is true, $P\to Q$ is still false since $Q$ is just coincidentally true when $P$ is true.

        The original meaning of bi-implication $P\oto Q$ is defined similarly as that ``$P\oto Q$ is true'' means ``there is a certain mechanism between $P$ and $Q$ in some scope, which guarantees bi-directional truth-preserving''.

        \subsection{Non-truth-functionality of implication}

        Syntactically, all the logical connectives $\{\top,\bot,\neg,\wedge,\vee,\to,\oto\}$ are functions mapping zero, one or two formulas to a formula. But semantically, implication $\to$ and bi-implication $\oto$ are essentially different from the others. The following is a comparison of implication $\to$ and conjunction $\wedge$ as an example.

        Let $L$ be the set of all the formulas representing a propositional language. Let $U$ be the set of all valuation functions from $L$ to the set of truth values $\{\mathrm{T,F}\}$. Then:

        According to the meaning of conjunction, there exists a function $f:\{\mathrm{T,F}\}^2\to\{\mathrm{T,F}\}$ such that for \emph{any} $v\in U$ and \emph{any} $\phi,\psi\in L$, $v(\phi\wedge\psi)=f(v(\phi),v(\psi))$, That is, $\wedge$ is truth-functional.

        In contrast, according to the original meaning of implication, there exists no function $f:\{\mathrm{T,F}\}^2\to\{\mathrm{T,F}\}$ such that for \emph{any} $v\in U$ and \emph{any} $\phi,\psi\in L$, $v(\phi\to\psi)=f(v(\phi),v(\psi))$. That is, $\to$ is not truth-functional.

        For instance, in terms of the meaning of conjunction, if we know that $P$ is true and $Q$ is true, then we know for sure that $P\wedge Q$ is true; whereas, in terms of the original meaning of implication, even if we know that $P$ is true and $Q$ is true, we are not sure about whether $P\to Q$ is true or false, because this depends on whether there is ``a certain mechanism'' from $P$ to $Q$ to guarantee that whenever $P$ is true $Q$ must be true.

        In fact, we know surely the truth value of $P\wedge Q$ for all possible four combinations of truth values $\{\mathrm{TT,TF,FT,FF}\}$ of $P$ and $Q$, whereas we do not know the truth value of $P\to Q$ from that of $P$ and $Q$ except for the only one combination $\mathrm{TF}$: only when $P$ is true and $Q$ is false, we know for sure that $P\to Q$ is false, because (only) $P\wedge\neg Q\to \neg(P\to Q)$ is valid according to the original meaning of implication (cf. Subsection \ref{subsec:Incorrectness of material implication}).

        From the semantic viewpoint, implication $\to$ is a binary \emph{relation} rather than a binary operation like conjunction $\wedge$ and disjunction $\vee$. Specifically, the implication $\to$, as a binary relation on $L$, the truth value of $\phi\to\psi$ is determined by whether $(\phi,\psi)$ belongs to the set of the relation $\to$ that is a subset of $L\times L$ (semantically, it is just the T-set of the implication, cf. Subsection \ref{subsec:Semantics of IRL}), rather than by the truth value of $\phi$ and the truth value of $\psi$. We call $\to$ \emph{implication relation} when emphasizing this fact.

        So is for bi-implication $\oto$.

        \subsection{Incorrectness of ``material implication'' as implication}
        \label{subsec:Incorrectness of material implication}

        Consider the relation between $P\to Q$ and $\neg P\vee Q$ in terms of the semantics of implication defined above.

        ``$\neg P\vee Q$ is true'' consists of three cases: $P$ is true and $Q$ is true, $P$ is false and $Q$ is true, $P$ is false and $Q$ is false. The first case ``$P$ is true and $Q$ is true'' includes not only the sole situation where $P\to Q$ is true (i.e. $Q$ is necessarily true when $P$ is true), but also the the situation ``$P$ is true and $Q$ is coincidentally true'' which is not a context where $P\to Q$ is true. In other words, ``$\neg P\vee Q$ is true'' contains not only the case where ``$P\to Q$ is true'' but also more cases. In symbols, $(P\to Q)\to\neg P\vee Q$ holds whereas $\neg P\vee Q \to (P\to Q)$ does not.

        ``$\neg P\vee Q$ is false'' contains only one case: $P$ is true and $Q$ is false that includes two sub-cases, i.e. ``$P$ is true and $Q$ is necessarily false'' and ``$P$ is true and $Q$ is coincidentally false''. The two sub-cases is a subset of situations where $P\to Q$ is false (cf. Subsection \ref{subsec:OriginalMeaningOfImplication}). In symbols, $\neg(\neg P\vee Q)\to\neg(P\to Q)$ (i.e. $P\wedge\neg Q\to\neg(P\to Q)$) holds, whereas $\neg(P\to Q)\to \neg(\neg P\vee Q)$ does not. This result is the same as above in terms of contraposition.

        In summary, $(P\to Q)\oto\neg P\vee Q$ is generally not a tautology since
        \begin{equation*}
            \begin{aligned}
                (P\to Q)\to\neg P\vee Q &\text{ (Implication-to-Disjunction) is valid while} \\
                \neg P\vee Q \to (P\to Q) &\text{ (Disjunction-to-Implication) is generally invalid.}
            \end{aligned}
        \end{equation*}

        Some researchers have addressed this issue with different motivations or explanations, such as \citet{maccoll}, \citet{bronstein}, \citet{woods}, and \citet{dale}.\\

        \subsection{Difference from ``strict implication''}

        The meaning of ``strict implication'' of modal logics is NOT the same as the original meaning of implication described above. In a classical modal logic:
        \begin{itemize}
            \item $\Box (P\to Q)$ is true means ``necessarily $P$ implies $Q$'' is true. Equivalently,
            \item $\Box (P\to Q)$ is false, or $\neg\Box (P\to Q)$ is true, or $\Diamond\neg (P\to Q)$ is true, means ``not necessarily $P$ implies $Q$'' is true, or ``possibly not $P$ implies $Q$'' is true.
        \end{itemize}
        The negation of strict implication, ``it is \emph{possible not} that $P$ implies $Q$'', is not excluding the situation ``it is \emph{possible} that $P$ implies $Q$''. This is different from the original meaning of implication defined in previous Subsection \ref{subsec:OriginalMeaningOfImplication}.

        As a matter of fact, ``strict implication'' does not define the meaning of implication itself: let $S=P\to Q$, the strict implication $\Box (P\to Q)=\Box S$ says just ``necessarily $S$'', nothing about the implication $P\to Q$ itself.

        \subsection{Construction of the right logic}
        \label{subsec:Construction of the right logic}

        The goals of this work in construction of the logic are:
        \begin{enumerate}
            \item Classical scope. The logical system to be constructed is ``in the classical scope''. For example, it is a two-valued logic and a sub-system of classical logic, so that any operators new to classical logic like modal operators are not considered.
            \item Not conflicting with the original meaning of implication. The system should not contain invalid theorems in terms of the meaning of operators, especially the original meaning of implication $\to$. For example, Disjunction-to-Implication $\neg\phi\vee\psi\to(\phi\to\psi)$ must be rejected.
            \item Strength. Under the above two conditions, the system should be as strong as possible. For example, both LEM and ECQ are not rejected. LEM is irrelevant to the original meaning of implication. ECQ, which is necessary for retaining both disjunction introduction and disjunctive syllogism (cf. Section \ref{sec:Discussion}), does not conflict with the original meaning of implication (cf. Subsection \ref{subsec:Semantics of IRL}). And Implication-to-Disjunction $(\phi\to\psi)\to\neg\phi\vee\psi$ is kept since it follows directly from the original meaning of implication as explained previously.
        \end{enumerate}

        There are many choices of ways for constructing logics. Two main approaches are:
        \begin{enumerate}
            \item Choose bi-implication $\oto$ as the primitive one, then define implication $\to$ by bi-implication with the rule $\phi\to\psi\otto\phi\wedge\psi\oto\phi$ or $\phi\to\psi\otto\phi\vee\psi\oto\psi$.
            \item Choose implication $\to$ as the primitive one, then define bi-implication $\oto$ by implication with the rule $\phi\oto\psi\otto(\phi\to\psi)\wedge(\psi\to\phi)$.
        \end{enumerate}

        Because of its symmetry, bi-implication-based approach has an advantage in structure. The logic IRL presented in Section \ref{sec:Propositional logic IRL} and its algebraic system IRLA in Section \ref{sec:Semantics and the algebra IRLA}, are constructed in the bi-implication-based approach. The natural deduction for IRL is constructed in the implication-based approach that is more suitable to it.


        \section{Propositional logic IRL}
        \label{sec:Propositional logic IRL}

        \textbf{IRL} (Implication-Relation Logic) is used to refer the logic presented in this section.

        \subsection{Definitions and notations}
        \label{subsec:Definitions and notations 1}

        The following definitions and notations are adopted for the systems proposed in this paper. The description of concepts is general, e.g. ``formula'' in this subsection may mean not only the ``(well-formed) formula'' of a propositional language or first-order language, but also a statement like a rule expression or a sentence of some meta-language being used.

        \subsubsection*{Definitions}

        \begin{definition}[Rule of inference]
            \label{def:Rule of inference}
            Let $X$ be a (normally non-empty) finite set of formulas (premises) and $y$ be a formula (conclusion).

            A (single-conclusion) \emph{rule of inference} is an expression of the form
            $$X\tto y$$
            that is interpreted as that if every element of $ X $ holds, then $y$ holds.

            Let $Y=\{y_1,...,y_n\}$ be a non-empty finite set of formulas (conclusions).

            A \emph{multiple-conclusion rule}
            $$X\tto Y$$ is defined as the set of single-conclusion rules $\{X\tto y_1,...,X\tto y_n\}$.

            A bi-inference rule
            $$X\otto Y$$
            is defined as the union of $X\tto Y$ and $Y\tto X$ (treat a single-conclusion rule as a singleton for the union operation).
        \end{definition}

        \begin{remark}\;

            \emph{Empty-premise rule}: In Definition \ref{def:Rule of inference}, a rule with an empty set of premises represents actually axiom(s) or theorem(s).

            \emph{Domain-specific rule}: A domain-specific rule is sometimes called an axiom or a theorem in that domain. For example, in a Boolean lattice extended from a Boolean algebra, $x\le y\otto x\wedge y=x$ is usually listed as the defining axiom for its partial order ($\le$), although it is essentially a (bi-inference) domain-specific rule.
        \end{remark}

        \begin{definition}[Proof and theorem]
            \label{def:Proof and theorem}
            Let $X$ be a (possibly empty) finite set of formulas (premises) and $y$ be a formula (conclusion).

            A \emph{proof} of $y$ from $X$, denoted
            $$X\vdash y,$$
            is a non-empty finite sequence of formulas $(y_1,...,y_n)$ with $y_n=y$, where $y_k\;(1\le k\le n)$ is an axiom, a theorem, or a result from a subset of $X\cup\{y_1,...,y_{k-1}\}$ by applying a rule of inference.
            \begin{itemize}
                \item  An \emph{unconditional proof} is a proof with an empty set of premises. Formula $y$ is called a \emph{theorem} if there exists an unconditional proof of $y$, denoted
                $$\vdash y.$$
                \item A \emph{Conditional proof} is a proof with a non-empty set of premises, i.e. $X\vdash y$ where $X$ is non-empty. Note that a member of the set $X$ of premises is not a member of its \emph{proof sequence} $(y_1,...,y_n)$.
            \end{itemize}
            The length or steps of a proof is the length of its proof sequence.
        \end{definition}

        \begin{remark}
            \label{rem:Proof sequences}
            According to Definition \ref{def:Proof and theorem}:
            \begin{itemize}
                \item All axioms are theorems. Therefore, any member of a proof sequence is \emph{just either a theorem or a result} from a subset of premises and previous members of the proof sequence \emph{by applying a rule}.
                \item In a proof sequence $(y_1,...,y_n)$ of an unconditional proof, each subsequence of the form $(y_1,...,y_k)\,(1\le k\le n)$ is a unconditional proof, thus every member $y_k\,(1\le k\le n)$ is a theorem.
            \end{itemize}
        \end{remark}

        \begin{definition}[Derivable rule]
            \label{def:Derivable rule}
            Let $X$ and $Y$ be finite sets of formulas where $Y$ is non-empty.

            A rule $X\tto Y$ is \emph{derivable} if and only if there exists a proof $X\vdash Y$ where rules used are only existing ones. An existing rule is a primitive rule or a derived rule.
        \end{definition}

        In terms of the definitions of proof and derivable rule, any primitive rule is derivable.

        \begin{remark}[Meta rule and meta proof]
            If a premise or a conclusion of a rule or a proof contains a symbol not belonging to the object language, then the rule or the proof is a \emph{meta rule} or \emph{meta proof}. For example, the rule $(\phi\vdash\psi)\tto(\vdash\phi\to\psi)$ of conditional proof is a meta rule, a proof $(\phi\vdash\psi)\vdash(\vdash\phi\to\psi)$ of the deduction theorem is a meta proof.
        \end{remark}

        \begin{definition}[Uniform substitution]
            \label{def:Uniform substitution}
            Let $x_1,...,x_k$ be variables, $t_1,...,t_k$ be terms, and $E$ be an expression containing variables.

            A \emph{uniform substitution} is a set $\{t_1/x_1,...,t_k/x_k\}$ and $E\{t_1/x_1,...,t_k/x_k\}$ means to (simultaneously) substitute $t_i$ for every occurrence of each $x_i$ in $E$. A member $x_i/x_i$ in a substitution $\{...,x_i/x_i,...\}$ can be omitted, e.g. $\{t_1/x_1,x_2/x_2,t_3/x_3\}$ is the same as $\{t_1/x_1,t_3/x_3\}$.

            The \emph{rule of uniform substitution} is the form $E\tto E\{t_1/x_1,...,t_k/x_k\}$, where $E$ is an axiom, a theorem, a rule, or a meta rule of a formal system.
        \end{definition}

        The rule of uniform substitution is common to logical systems for instantiating (schemas of) axioms, theorems, rules and meta rules, and it is often used tacitly.

        \begin{definition}[Subformula replacement]
            \label{def:Subformula replacement}
            Let $x,y,z$ be formulas.

            The expression $z[x\mapsto y]$ denotes \emph{any} member of the set of all possible results where each result is obtained via replacing one or more \emph{arbitrarily selected} occurrence(s) of $x$ by $y$ in $z$.

            For instance, $(x\wedge y\to x\wedge y)[x\wedge y\mapsto y]$ represents any result in the set $\{y\to x\wedge y,\, x\wedge y\to y,\, y\to y\}$.
        \end{definition}

        \subsubsection*{Presentation of proofs}

        The presentation format of a proof is as follows.

        An unconditional proof $\vdash\psi$ of length $n$ is presented in the format:

        \begin{tabular}{lll}
            $1$&$\psi_1$&reason\\
            $\vdots$&$\vdots$&$\vdots$\\
            $n-1$&$\psi_{n-1}$&reason\\
            $n$&$\psi$&reason
        \end{tabular}\\

        A conditional proof $\{\phi_1,...,\phi_m\}\vdash\psi$ of length $n$ is presented in the format:

        \begin{tabular}{lll}
            $1$&$\phi_1$&Premise (Assumption)\\
            $\vdots$&$\vdots$&$\vdots$\\
            $m$&$\phi_m$&Premise (Assumption)\\
            $m+1$&$\psi_1$&reason\\
            $\vdots$&$\vdots$&$\vdots$\\
            $m+n-1$&$\psi_{n-1}$&reason\\
            $m+n$&$\psi$&reason
        \end{tabular}\\

        A proof (unconditional or conditional) can has other proofs  (unconditional or conditional) embedded within it. In the presentations, this nested structure of proofs must be clearly indicated in a consistent manner such as indentation used in this paper.

        An example that an unconditional proof has an inner conditional proof:

        \begin{tabular}{lll}
            1 &$\psi_{11}$& theorem-name\\
            2 &\qquad$\phi_{21}$&Premise\\
            3 &\qquad$\psi_{21}$&theorem-name\\
            4 &\qquad$\psi_{22}$&2,3: rule-name (From $\{2,3\}$ by rule-name)\\
            5 &$\phi_{21}\to\psi_{22}$& 2--4: Rule of conditional proof\\
            6 &$\psi_{12}$&1,5: rule-name
        \end{tabular}\\

        Note that \emph{the result of a conditional proof (by the rule of conditional proof) is not a member of its own proof sequence}, so it must be \emph{unindented} to just one level higher (one layer outer) as shown in the example (line 5).

        \vspace{\baselineskip}
        In some expressions, ``$\alpha_0\backsim\alpha_1\backsim\cdots\backsim\alpha_n$'' is a shorthand of ``$\alpha_0\backsim\alpha_1$, $\alpha_1\backsim\alpha_2$, ..., $\alpha_{n-1}\backsim\alpha_n$'', where $\backsim$ stands for some binary relations or operators such as $\oto$, $\to$, $\otto$, $\tto$, $=$, $\le$.

        Hence, in the presentation of a proof, the format

        \begin{tabular}{lll}
            &$\alpha_0$&\\
            $1$&\qquad$\backsim\alpha_1$&\\
            $\vdots$&\qquad$\vdots$&\\
            $n$&\qquad$\backsim\alpha_n$&
        \end{tabular}

        is a shorthand of the format

        \begin{tabular}{lll}
            $1$&$\alpha_0\backsim\alpha_1$&\\
            $\vdots$&$\vdots$&\\
            $n$&$\alpha_{n-1}\backsim\alpha_n$&
        \end{tabular}

        \subsubsection*{A note on nested proofs}

        By \emph{definition} (Definition \ref{def:Proof and theorem}, and cf. Remark \ref{rem:Proof sequences}), any member of a proof sequence must be either a theorem or a result from a subset of premises and previous members by a rule. Specifically, let $X\vdash y_n$ has a proof sequence $(y_1,...,y_n)$, then any $y_k$ must be in one of the following two cases:
        \begin{itemize}
            \item $y_k$ is a theorem, or it is assumed be a theorem in a premise of an \emph{outer proof}, or it is a theorem derived in an \emph{outer proof};
            \item $y_k$ is obtained from a subset of $X\cup\{y_1,...,y_{k-1}\}$ by applying a rule, or by an assumed rule in a premise of an \emph{outer proof}, or by a rule derived in an \emph{outer proof}.
        \end{itemize}

        \emph{No other cases are allowed.}

        \vspace{\baselineskip}
        This is further illustrated in Table \ref{tab:Only two cases allowed for members of proof sequence}.

        Attention must be paid to this, especially for a \emph{nested conditional proof}.

        \begin{table}[H] 
            \centering
            \small
            \setlength{\tabcolsep}{6pt} 
            \renewcommand{\arraystretch}{1.5} 
            \begin{threeparttable}
                \caption{Only two cases allowed for members of proof sequence}
                \label{tab:Only two cases allowed for members of proof sequence}
                \begin{tabular}{p{0.125\linewidth}|p{0.155\linewidth}|p{0.65\linewidth}}
                    \hline
                    An outer proof & The inner proof & Comment\\
                    \hline
                    $A_1$ &  &  \\
                    $\vdots$ &  &  \\
                    $A_m$    &  & $\{A_1,...,A_m\}$ is the set of premises/assumptions of the outer proof \tnote{a} \\
                    \hline
                    $B_1$ &  &  \\
                    $\vdots$ &  &  \\
                    $B_i$    &  & For $B_i$, a member of the outer proof sequence, \emph{only two cases are allowed}: it is a theorem, or it is a result from a subset of $\{A_1,...,A_m\}\cup \{B_1,...,B_{i-1}\}$ by applying a rule \tnote{b} \\
                    $\vdots$ &  &  \\
                    $B_j$ &  &  \\
                    \hline
                    & $P_1$ &  \\
                    & $\vdots$ &  \\
                    & $P_n$ & $\{P_1,...,P_n\}$ is the set of premises/assumptions of the inner proof \tnote{a} \\
                    \hline
                    & $Q_1$ &  \\
                    & $\vdots$ &  \\
                    & $Q_k$ & For $Q_k$, a member of the inner proof sequence, \emph{only two cases are allowed}: (a) it is a theorem, or it is assumed be a theorem in $\{A_1,...,A_m\}$, or it is a theorem in $\{B_1,...,B_j\}$; or (b) it is a result from a subset of $\{P_1,...,P_n\}\cup \{Q_1,...,Q_{k-1}\}$ by a rule, or by a rule assumed in $\{A_1,...,A_m\}$, or by a rule in $\{B_1,...,B_j\}$ \\
                    & $\vdots$ &  \\
                    \hline
                \end{tabular}
                \begin{tablenotes}
                    \item[a] If the set of premises is empty, it is an unconditional proof; otherwise, it is a conditional proof.
                    \item[b] If the outer proof itself is embedded in another proof, treat it the same way as the inner proof in the table.
                \end{tablenotes}
            \end{threeparttable}
        \end{table}

        \begin{example}[An \emph{illegal} proof]
            Proof of $\phi\to(\psi\to\phi)$:\\
            \begin{tabular}{lll}
                1&\qquad$\phi$&Premise (outer)\\
                2&\qquad\qquad$\psi$&Premise (inner)\\
                3&\qquad\qquad$\phi$&1: Reiteration -- \emph{Not allowed}\\
                4&\qquad$\psi\to\phi$&2--3: Rule of conditional proof\\
                5&$\phi\to(\psi\to\phi)$&1--4: Rule of conditional proof
            \end{tabular}\\

            Line 3 is \emph{not allowed} because $\phi$ in this line does not belong to any of the allowed cases: (a) it is a theorem or a theorem indicated in an outer proof ($\phi$ is not assumed be a theorem in the outer premise i.e. line 1); (b) it is a result from a subset of premises and previous members of its own layer (that must be a subset of $\{\psi\}$ from line 2) by applying a rule or by an rule in an outer proof.
        \end{example}

        \begin{example}[A legal proof]
            Proof of $(\top\to\phi)\to(\psi\to\phi)$:\\
            \begin{tabular}{lll}
                1&\qquad$\top\to\phi$&Premise (outer)\\
                2&\qquad\qquad$\psi$&Premise (inner)\\
                3&\qquad\qquad$\top$&A theorem ($\top$ is indeed a theorem)\\
                4&\qquad$\psi\to\top$&2--3: Rule of conditional proof\\
                5&\qquad$\psi\to\phi$&4,1: Rule of transitivity\\
                6&$(\top\to\phi)\to(\psi\to\phi)$&1--5: Rule of conditional proof
            \end{tabular}\\
        \end{example}

        Sometimes, a meta proof (for a meta language) contains an ordinary proof (for an object language) such as the following example, which details a meta proof of a meta rule.

        \begin{example}[A legal meta proof]
            Proof of $(\vdash\psi)\tto(\vdash\phi\to\psi)$:\\
            \begin{tabular}{lll}
                1&\qquad$\vdash\psi$&Premise (meta proof)\\
                2&\qquad\qquad $\phi$&Premise (proof)\\
                3&\qquad\qquad $\psi$&1 ($\psi$ is assumed in line 1 be a theorem)\\
                4&\qquad$\phi\vdash\psi$&2--3: Definition of conditional proof\\
                5&\qquad$\vdash\phi\to\psi$&4: Rule of conditional proof\\
                6&$(\vdash\psi)\vdash(\vdash\phi\to\psi)$&1--5: Definition of conditional proof\\
                7&$(\vdash\psi)\tto(\vdash\phi\to\psi)$&6: Definition of derivable rule\\
            \end{tabular}\\
        \end{example}

        More examples can be found in Subsections \ref{subsec:Elementary derived results} and \ref{subsec:More derivations in IRL}.

        \subsection{Logic IRL}
        \label{subsec:Logic IRL}

        Let $\{\top,\bot,\neg,\wedge,\vee,\oto,\to\}$ be the set of logical connectives (operations, where $\top$ and $\bot$ are nullary operators).

        Suppose there is a countable (finite or countably infinite) set $K$ of atomic or primitive propositions. Let $A,B,C,...$ be \emph{primitive propositional variables} ranging over $K$.

        A (well-formed) formula is an expression built as usual from a finite number of primitive propositional variables $(A_1,...A_k)$ using operations $\{\top,\bot,\neg,\wedge,\vee,\oto,\to\}$ finite times. Let $L$ be the set of all formulas and $\phi,\psi,\chi,...$ be \emph{formula variables} ranging over $L$.

        A \emph{formula schema} is a function expression $F(\phi_1,...,\phi_n)$ (where $F:L^n \to L$) built from a finite number of formula variables $(\phi_1,...,\phi_n)$ using operations $\{\top,\bot,\neg,\wedge,\vee,\oto,\to\}$ finite times.

        Thus: the set $K$ of primitive propositions is a generating set of the set $L$ of all formulas; the set $L$ is countably infinite and is closed under the set of operations $\{\top,\bot,\neg,\wedge,\vee,\oto,\to\}$, so the set of all schemas is also countably infinite and closed under all the operations.

        \vspace{\baselineskip}
        The propositional logic IRL is defined by the following axioms and primitive rules of inference, where $\phi,\psi,\chi,...$ are arbitrary formulas for which one can substitute any schemas.

        \subsubsection*{Axioms}

        \begin{enumerate}
            \item $\phi\oto\phi$ (Reflexivity of bi-implication). Symmetry and transitivity can be derived.
            \item $\phi\wedge\psi \oto \psi\wedge\phi$, $\phi\vee\psi \oto \psi\vee\phi$ (Commutativity).
            \item $(\phi\wedge\psi)\wedge\chi\oto \phi\wedge(\psi\wedge\chi)$, $(\phi\vee\psi)\vee\chi\oto \phi\vee(\psi\vee\chi)$ (Associativity). Redundant.
            \item $\phi\wedge(\psi\vee\chi) \oto (\phi\wedge\psi)\vee(\phi\wedge\chi)$, $\phi\vee(\psi\wedge\chi) \oto (\phi\vee\psi)\wedge(\phi\vee\chi)$ (Distributivity).
            \item $\phi\wedge\top\oto \phi$, $\phi\vee\bot\oto\phi$ (Identity elements, top and bottom).
            \item $\phi\wedge\neg\phi \oto \bot$, $\phi\vee\neg\phi \oto \top$ (Complement/Negation, LNC and LEM).
            \item $(\phi \rightarrow \psi) \oto (\phi\wedge\psi \oto \phi)$ (Definition of implication).
        \end{enumerate}

        \subsubsection*{Primitive rules of inference}

        \begin{enumerate}
            \item $\{\phi,\psi\}\tto\phi\wedge\psi$ (Conjunction introduction).
            \item $\phi\wedge\psi\tto\phi$ (Conjunction elimination).
            \item $\{\phi\oto\psi,\chi\}\tto\chi[\phi\mapsto\psi]$ (Replacement property of bi-implication).
            \item $(\phi\vdash\psi)\tto(\vdash\phi\to\phi)$ (Rule of conditional proof). Meta Rule.
        \end{enumerate}

        \begin{remark}\;
            \begin{itemize}
                \item Symmetry and transitivity of bi-implication are proved in the next subsection.
                \item Associativity can be proved following a procedure similar to that by \citet{huntington}.
                \item In a proof, by definition, any member of its proof sequence must be either a theorem (or a theorem indicated in an outer layer), or a result from a subset of the premises and previous members in its own layer by a rule (or by a rule indicated in an outer layer), no other cases are allowed. It is important to pay attention to this, especially for a nested conditional proof (cf. Subsection \ref{subsec:Definitions and notations 1}).
                \item Rule of conditional proof of the general form $(\Gamma\cup\{\phi\}\vdash\psi)\tto(\Gamma\vdash\phi\to\psi)$, where $\Gamma$ is a set of formulas, does not hold for IRL: if it holds, then from $\{\phi,\psi\}\vdash\phi$ it follows $\phi\vdash\psi\to\phi$, hence $\vdash\phi\to(\psi\to\phi)$, contrary to that $\nvdash_{IRL}\phi\to(\psi\to\phi)$ (cf. Subsection \ref{subsec:D2I is underivable in IRL}).
                \item The logic \emph{IRL is consistent} as a sub-system of classical (propositional) logic (cf. Remark \ref{rem:CL is a super-system of IRL} and Subsection \ref{subsec:IRL-CL relationship}).
            \end{itemize}
        \end{remark}

        \subsection{Elementary derived results}
        \label{subsec:Elementary derived results}

        \subsubsection*{Symmetry and transitivity of bi-implication}

        \begin{derivedrule}[Symmetry of bi-implication]
            \label{derule:Symmetry of bi-implication}
            $$\phi\oto \psi\tto \psi\oto \phi.$$
        \end{derivedrule}

        \begin{proof}\;\\
            \begin{tabular}{lll}
                1 &\qquad $\phi\oto \psi$ & Premise\\
                2 &\qquad $\phi\oto \phi$ & Reflexivity\\
                3 &\qquad $(\phi\oto \phi)[\phi\mapsto \psi]$ & 1,2: Replacement\\
                4 &\qquad $\phi[\phi\mapsto \psi]\oto\phi$ & 3: \emph{Selected replacement}\\
                5 &\qquad $\psi\oto \phi$ & 4: Result of the replacement\\
                6 & $\phi\oto \psi\vdash \psi\oto \phi$ & 1-5: Definition of conditional proof\\
                7 & $\phi\oto \psi\tto \psi\oto \phi$ & 6: Definition of derivable rule\\
            \end{tabular}\\
        \end{proof}

        Note: The definitions of conditional proof and derivable rule can be used tacitly in a proof as usual. For example, the last two lines in the above proof can be omitted.

        \begin{derivedrule}[Transitivity of bi-implication]
            \label{derule:Transitivity of bi-implication}
            $$\{\phi\oto \psi,\psi\oto \chi\}\tto \phi\oto \chi.$$
        \end{derivedrule}

        \begin{proof}\;\\
            \begin{tabular}{lll}
                1 & $\phi\oto \psi$ & Premise\\
                2 & $\psi\oto \chi$ & Premise\\
                3 & $(\phi\oto \psi)[\psi\mapsto \chi]$ & 2,1: Replacement\\
                4 & $\phi\oto \chi$ & 3: Result of the replacement\\
            \end{tabular}\\
        \end{proof}

        \subsubsection*{Equivalence of bi-implication theorem and bi-inference rule}

        \begin{metarule}[Bi-implication theorem to bi-inference rule]
            \label{metarule:Bi-implication theorem to bi-inference rule}
            $$(\vdash\phi\oto\psi)\tto(\phi\otto\psi).$$
        \end{metarule}

        \begin{proof}\;\\
            \begin{tabular}{lll}
                1 & $\vdash\phi\oto\psi$ & Premise (begin meta proof)\\
                2 & \qquad $\phi$ & Premise (begin proof 1)\\
                3 & \qquad $\phi\oto\psi$ & 1 (outside assumed theorem)\\
                4 & \qquad $\phi[\phi\mapsto\psi]$ & 3,2: Replacement\\
                5 & \qquad $\psi$ & 4: Result of the replacement (end proof 1)\\
                6 & $\phi\vdash\psi$ & 2--5: Definition of conditional proof\\
                7 & $\phi\tto\psi$ & 6: Definition of derivable rule\\
                8 & \qquad $\psi$ & Premise (begin proof 2)\\
                9 & \qquad $\phi\oto\psi$ & 1 (outside assumed theorem)\\
                10 & \qquad $\psi\oto\phi$ & 9: Symmetry of bi-implication\\
                11 & \qquad $\psi[\psi\mapsto\phi]$ & 10,8: Replacement\\
                12 & \qquad $\phi$ & 11: Result of the replacement (end proof 2)\\
                13 & $\psi\vdash\phi$ & 8--12: Definition of conditional proof\\
                14 & $\psi\tto\phi$ & 13: Definition of derivable rule\\
                15 & $\psi\otto\phi$ & 7,14: Definition of bi-inference rule (end meta proof)
            \end{tabular}\\
        \end{proof}

        By this meta rule, from bi-implication axioms we obtain bi-inference rules directly. For example, from $\phi\wedge\psi\oto\psi\wedge\phi$, $(\phi\wedge\psi)\wedge\chi\oto \phi\wedge(\psi\wedge\chi)$ and $(\phi\to\psi)\oto(\phi\wedge\psi\oto\phi)$, we obtain rules $\phi\wedge\psi\otto\psi\wedge\phi$, $(\phi\wedge\psi)\wedge\chi\otto \phi\wedge(\psi\wedge\chi)$ and $\phi\to\psi\otto\phi\wedge\psi\oto\phi$, respectively.

        \begin{derivedrule}[Bi-implication introduction]
            \label{derivedrule:Bi-implication introduction}
            $$\{\phi\to \psi,\psi\to \phi\}\tto \phi\oto \psi.$$
        \end{derivedrule}

        \begin{proof}\;\\
            \begin{tabular}{lll}
                1 & $\phi\to \psi$ & Premise\\
                2 & $\psi\to \phi$ & Premise\\
                3 & $\phi\wedge\psi\oto \phi$ & 1: Definition of $\to$\\
                4 & $\phi\oto \phi\wedge\psi$ & 3: Symmetry of $\oto$\\
                5 & $\psi\wedge\phi\oto \psi$ & 2: Definition of $\to$\\
                6 & $\phi\wedge\psi\oto \psi\wedge\phi$ & Commutativity\\
                7 & $\phi\wedge\psi\oto \psi$ & 6,5: Transitivity of $\oto$\\
                8 & $\phi\oto \psi$ & 4,7: Transitivity of $\oto$
            \end{tabular}\\
        \end{proof}

        \begin{metarule}[Bi-inference rule to bi-implication theorem]
            \label{metarule:Bi-inference rule to bi-implication theorem}
            $$(\phi\otto \psi)\tto (\vdash \phi\oto \psi).$$
        \end{metarule}

        \begin{proof}\;\\
            \begin{tabular}{lll}
                1 & $\phi\otto \psi$ & Premise\\
                2 & \qquad\qquad $\phi$ & Premise\\
                3 & \qquad\qquad $\psi$ & 2: The rule assumed in 1\\
                4 & \qquad$\phi\to \psi$ & 2--3: Conditional proof\\
                5 & \qquad$\psi\to \phi$ & 1--4: Similarly\\
                6 & \qquad $\phi\oto \psi$ &4,5: Bi-implication introduction\\
                7 & $\vdash \phi\oto \psi$ &4--6: Unconditional proof
            \end{tabular}\\
        \end{proof}

        \begin{corollary}[Equivalence of bi-implication theorem and bi-inference rule]
            \label{corollary:Equivalence of bi-implication theorem and bi-inference rule}
            $$(\vdash \phi\oto \psi)\otto(\phi\otto \psi).$$
        \end{corollary}

        \begin{proof}
            It follows from Meta Rules \ref{metarule:Bi-implication theorem to bi-inference rule} and \ref{metarule:Bi-inference rule to bi-implication theorem} by the definition of bi-inference rule.
        \end{proof}

        Thus, we do not need to distinguish bi-implication theorem $\vdash\phi\oto\psi$ and bi-inference rule $\phi\otto\psi$: if we have one of the two, we have both.

        \subsubsection*{Tacit use of symmetry properties}

        \begin{metarule} Let $F(\phi\circ\phi)$ be a schema containing $\phi\circ\psi$, where $\circ$ stands for $\oto$, $\wedge$, or $\vee$. Then
            $$F(\phi\circ\psi)\otto F(\psi\circ\phi).$$
        \end{metarule}

        \begin{proof} We already have $(\phi\circ\psi)\oto(\psi\circ\phi)$ for $\oto$, $\wedge$ and $\vee$, so the proof can be done below:\\
            \begin{tabular}{lll}
                1 & \qquad $F(\phi\circ\psi)$ & Premise\\
                2 & \qquad $(\phi\circ\psi)\oto(\psi\circ\phi)$ & Known theorem\\
                3 & \qquad $F(\psi\circ\phi)$ & 2,1: Replacement\\
                4 & $F(\phi\circ\psi)\tto F(\psi\circ\phi)$ & 1--3: Definition of derivable rule\\
                5 & $F(\psi\circ\phi)\tto F(\phi\circ\psi)$ & 1--4: Similarly\\
                6 & $F(\psi\circ\phi)\otto F(\phi\circ\psi)$ & 4,5: Definition of bi-inference rule.
            \end{tabular}\\
        \end{proof}

        Thus, the symmetry properties of $\oto$, $\wedge$ and $\vee$ can be used tacitly in a proof, viewing $F(\phi\circ\psi)$ and $F(\psi\circ\phi)$ as the same thing.

        \subsubsection*{Modus ponens}

        \begin{derivedrule}[Modus ponens]
            \label{derule:Modus ponens}
            $$\{\phi\to\psi,\phi\}\tto\psi.$$
        \end{derivedrule}

        \begin{proof}\;\\
            \begin{tabular}{lll}
                1&$\phi\to\psi$&Premise\\
                2&$\phi$&Premise\\
                3&$\phi\oto\phi\wedge\psi$&1: Definition of implication\\
                4&$\phi\wedge\psi$&3,2: Replacement\\
                5&$\psi$&4: Conjunction elimination\\
            \end{tabular}\\
        \end{proof}

        \subsubsection*{Equivalence of the three forms}

        With the rules of commutativity and associativity for conjunction, $\phi_1\wedge\cdots\wedge\phi_n$ is well defined for any $n\ge 1$. And from the bidirectional reflexivity rule $\phi\otto\phi$ we get its unidirectional form $\phi\tto\phi$ by the definition of bi-inference rule. Therefore, the rules of conjunction introduction and elimination can be generalized to $\{\phi_1,...,\phi_n\}\tto\phi_1\wedge\cdots\wedge\phi_n$ and $\phi_1\wedge\cdots\wedge\phi_n\tto\phi_i$ respectively, for all finite $n\ge 1$.

        \begin{metarule}[Rule of conditional proof with multiple premises]
            \label{metarule:Rule of conditional proof with multiple premises}
            $$(\{\phi_1,...,\phi_n\}\vdash\psi)\tto(\vdash\phi_1\wedge\cdots\wedge\phi_n\to\psi).$$
        \end{metarule}
        \begin{proof}
            Suppose $\{\phi_1,...,\phi_n\}\vdash\psi$ has a proof sequence $(\psi_1,...,\psi)$, then $(\phi_1,...,\phi_n,\psi_1,...,\psi)$ is a proof sequence of $\phi_1\wedge\cdots\wedge\phi_n\vdash\psi$ as shown below:\\
            \begin{tabular}{lll}
                0 & $\phi_1\wedge\cdots\wedge\phi_n$ & Premise\\
                1 & $\phi_1$ & 0: Conjunction elimination\\
                \vdots  & $\vdots$ & \vdots\\
                n & $\phi_n$ & 0: Conjunction elimination\\
                & $\psi_1$ & 1...n: First proof-sequence member of $\{\phi_1,...,\phi_n\}\vdash\psi$\\
                \vdots & $\vdots$ & \vdots\\
                & $\psi$ & 1...n: Last proof-sequence member of $\{\phi_1,...,\phi_n\}\vdash\psi$\\
            \end{tabular}\\
            Thus, from $\phi_1\wedge\cdots\wedge\phi_n\vdash\psi$ we obtain $\vdash\phi_1\wedge\cdots\wedge\phi_n\to\psi$ by the rule of conditional proof with single premise.
        \end{proof}

        \begin{metarule}(Converse of CP rule)
            \label{metarule:Converse of CP rule}
            $$(\vdash\phi_1\wedge\cdots\wedge\phi_n\to\psi)\tto(\{\phi_1,...,\phi_n\}\vdash\psi).$$
        \end{metarule}

        \begin{proof}\;\\
            \begin{tabular}{lll}
                1 & $\vdash\phi_1\wedge\cdots\wedge\phi_n\to\psi$ & Premise (meta proof)\\
                2 & \qquad $\phi_1$ & Premise (proof)\\
                \vdots & \qquad $\vdots$ & \vdots\\
                3 & \qquad $\phi_n$ & Premise (proof)\\
                4 &\qquad $\phi_1\wedge\cdots\wedge\phi_n$ & 2...3: Conjunction introduction\\
                5 & \qquad $\phi_1\wedge\cdots\wedge\phi_n\to\psi$ & 1 (assumed theorem)\\
                6 & \qquad $\psi$ & 5,4: Modus ponens\\
                7 & $\{\phi_1,...,\phi_n\}\vdash\psi$ & 2--6: Definition of conditional proof
            \end{tabular}\\
        \end{proof}

        \begin{metarule}[Equivalence of the three forms]
            \label{metarule:Equivalence of the tree forms}
            $$(\vdash\phi_1\wedge\cdots\wedge\phi_n\to\psi)\otto(\{\phi_1,...,\phi_n\}\vdash\psi)\otto(\{\phi_1,...,\phi_n\}\tto\psi).$$
        \end{metarule}

        \begin{proof}
            It follows from the definition of derivable rule, Meta Rule \ref{metarule:Rule of conditional proof with multiple premises} and Meta Rule \ref{metarule:Converse of CP rule}, and the definition of bi-inference rule.
        \end{proof}

        Therefore, if we have any one of the three forms, we have the other two. So we just need to find either an unconditional proof  $\vdash\phi_1\wedge\cdots\wedge\phi_n\to\psi$ or a conditional proof $\{\phi_1,...,\phi_n\}\vdash\psi$ to obtain both the theorem $\phi_1\wedge\cdots\wedge\phi_n\to\psi$ and the rule $\{\phi_1,...,\phi_n\}\tto\psi$.

        \subsubsection*{Other elementary results}

        \begin{derivedrule}[Leibniz's replacement]
            \label{derule:Leibniz's replacement}
            $$\phi\oto \psi\tto \chi\oto \chi[\phi\mapsto \psi].$$
        \end{derivedrule}

        \begin{proof}\;\\
            \begin{tabular}{lll}
                1 & $\phi\oto \psi$ & Premise\\
                2 & $\chi\oto \chi$ & Reflexivity\\
                3 & $(\chi\oto \chi)[\phi\mapsto \psi]$ & 1,2: Replacement\\
                4 & $\chi\oto \chi[\phi\mapsto \psi]$ & 3: \emph{Selected replacement}
            \end{tabular}\\
        \end{proof}

        \begin{theorem}[Top is a theorem]
            \label{thrm:Top is a theorem}
            $$\top.$$
        \end{theorem}

        \begin{proof}\;\\
            \begin{tabular}{lll}
                1&$\phi\oto\phi$&Reflexivity\\
                2&$(\phi\oto\phi)\wedge\top$&1: Identity element\\
                3&$\top$&2: Conjunction elimination
            \end{tabular}\\
        \end{proof}

        \begin{metarule}[Theorem like top]
            \label{derule:Theorem like top}
            $$(\vdash\psi)\tto(\vdash\phi\to\psi).$$
        \end{metarule}

        \begin{proof}\;\\
            \begin{tabular}{lll}
                1&$\vdash\psi$&Premise (meta proof)\\
                2&\qquad $\phi$&Premise (proof)\\
                3&\qquad $\psi$&1 (assumed theorem)\\
                4&$\vdash\phi\to\psi$&2--3: Conditional proof\\
            \end{tabular}\\
        \end{proof}

        A theorem functions like the top $\top$.

        \begin{metarule}[Theorem equivalence]
            \label{metarule:Theorem equivalence}
            $$\{\vdash\phi,\vdash\psi\}\tto(\vdash\phi\oto\psi).$$
        \end{metarule}

        \begin{proof}\;\\
            \begin{tabular}{lll}
                1&$\vdash\phi$& Premise\\
                2&$\vdash\psi$& Premise\\
                3&$\vdash\psi\to\phi$& 1: Theorem like top\\
                4&$\vdash\phi\to\psi$& 2: Theorem like top\\
                5&\qquad $\psi\to\phi$& 3 (theorem from an outer layer)\\
                6&\qquad $\phi\to\psi$& 4 (theorem from an outer layer)\\
                7&\qquad $\phi\oto\psi$& 5,6: Bi-implication introduction\\
                8&$\vdash\phi\oto\psi$&5--7: Unconditional proof
            \end{tabular}\\
        \end{proof}

        \begin{metarule}[Theorem representation]
            \label{derule:Theorem representation}
            $$(\vdash\phi\oto\top)\otto(\vdash\phi).$$
        \end{metarule}

        \begin{proof}\;\\
            \begin{tabular}{lll}
                1&$\vdash\phi\oto\top$&Premise\\
                2&\qquad $\phi\oto\top$&1 (assumed theorem)\\
                3&\qquad $\top$&Top is a theorem\\
                4&\qquad $\phi$&2,3: Replacement $[\top\mapsto\phi]$\\
                5&$\vdash\phi$&2--4: Unconditional proof\\
                &&\\
                1&$\vdash\phi$&Premise\\
                2&\qquad $\phi$&1 (assumed theorem)\\
                3&\qquad $\top$&Top is a theorem\\
                4&\qquad $\phi\oto\top$&2,3: Theorem equivalence\\
                5&$\vdash\phi\oto\top$&2--4: Unconditional proof\\
            \end{tabular}\\
        \end{proof}

        For example, $\phi\vee\neg\phi\oto\top$ and $\phi\vee\neg\phi$ are equivalent representations of LEM as a theorem.

        \subsection{More derivations in IRL}
        \label{subsec:More derivations in IRL}

        \subsubsection*{Some basic theorems or rules}

        Note that no matter what is proved is a theorem or a rule, we have both a theorem and a rule according to the equivalence of the three forms (Meta Rule \ref{metarule:Equivalence of the tree forms}) and the equivalence of bi-implication theorem and bi-inference rule (Corollary \ref{corollary:Equivalence of bi-implication theorem and bi-inference rule}).

        \begin{theorem}[Complementarity of top and bottom]
            \label{thrm:Complementarity of top and bottom}
            $$\neg\top\oto\bot,\;\neg\bot\oto\top.$$
        \end{theorem}

        \begin{proof}\;\\
            \begin{tabular}{lll}
                1&$\neg\top\oto\neg\top\wedge\top$&Identity element\\
                2&$\neg\top\wedge\top\oto\bot$&LNC\\
                3&$\neg\top\oto\bot$& 1,2: Transitivity of bi-implication\\
            \end{tabular}\\
            The second is proved dually.
        \end{proof}

        \begin{theorem}[Double negation]
            \label{thrm:Double Negation}
            $$\neg\neg\phi\oto\phi.$$
        \end{theorem}

        \begin{proof}\;\\
            \begin{tabular}{lll}
                &$\neg\neg\phi$&\\
                1&\qquad$\oto\neg\neg\phi\wedge\top$&Identity element\\
                2&\qquad$\oto\neg\neg\phi\wedge(\phi\vee\neg\phi)$&LEM: Leibniz's replacement\\
                3&\qquad$\oto(\neg\neg\phi\wedge\phi)\vee(\neg\neg\phi\wedge\neg\phi)$&Distributivity\\
                4&\qquad$\oto(\neg\neg\phi\wedge\phi)\vee\bot$&LNC: Leibniz's replacement\\
                5&\qquad$\oto(\neg\neg\phi\wedge\phi)\vee(\neg\phi\wedge\phi)$&LNC: Leibniz's replacement\\
                6&\qquad$\oto(\neg\neg\phi\vee\neg\phi)\wedge\phi$&Distributivity\\
                7&\qquad$\oto\top\wedge\phi$&LEM: Leibniz's replacement\\
                8&\qquad$\oto\phi$&Identity element\\
                9&$\neg\neg\phi\oto\phi$&1--8: Transitivity of bi-implication\\
            \end{tabular}\\
        \end{proof}

        Note: Leibniz's replacement may be used tacitly, e.g. ``LEM: Leibniz's replacement'' can be written simply as ``LEM'' in a reason of a proof line.

        \begin{theorem}[Idempotence]
            \label{thrm:Idempotence}
            $$\phi\wedge\phi\oto\phi,\;\phi\vee\phi\oto\phi.$$
        \end{theorem}

        \begin{proof}\;\\
            \begin{tabular}{lll}
                &$\phi\wedge\phi$&\\
                1&\qquad$\oto (\phi\wedge\phi) \vee \bot$&Identity\\
                2&\qquad$\oto (\phi\wedge\phi) \vee (\phi\wedge\neg\phi)$&LNC\\
                3&\qquad$\oto \phi\wedge(\phi\vee\neg\phi)$&Distributivity\\
                4&\qquad$\oto \phi\wedge\top$&LEM\\
                5&\qquad$\oto \phi$&Identity\\
                6&$\phi\wedge\phi\oto\phi$&1--5: Transitivity\\
            \end{tabular}\\
            The second is proved similarly.
        \end{proof}

        \begin{theorem}[Annihilator]
            \label{thrm:Annihilator}
            $$\bot\wedge\phi\oto\bot,\;\top\vee\phi\oto\top.$$
        \end{theorem}

        \begin{proof}\;\\
            \begin{tabular}{lll}
                &$\bot\wedge\phi$&\\
                1&\qquad$\oto(\bot\wedge\phi)\vee\bot$&Identity\\
                2&\qquad$\oto(\bot\wedge\phi)\vee(\neg\phi\wedge\phi)$&LNC\\
                3&\qquad$\oto(\bot\vee\neg\phi)\wedge\phi$&Distributivity\\
                4&\qquad$\oto\neg\phi\wedge\phi$&Identity\\
                5&\qquad$\oto\bot$&LNC\\
                6&$\bot\wedge\phi\oto\bot$&1--5: Transitivity\\
            \end{tabular}\\
            The second is proved similarly.
        \end{proof}

        \begin{theorem}[Absorption laws]
            \label{thrm:Absorption laws}
            $$\phi\wedge(\phi\vee\psi)\oto\phi,\;\phi\vee(\phi\wedge\psi)\oto\phi.$$
        \end{theorem}

        \begin{proof}\;\\
            \begin{tabular}{lll}
                &$\phi\wedge(\phi\vee\psi)$&\\
                1&\qquad$\oto(\phi\vee\bot)\wedge(\phi\vee\psi)$&Identity\\
                2&\qquad$\oto\phi\vee(\bot\wedge\psi)$&Distributivity\\
                3&\qquad$\oto\phi\vee\bot$&Annihilator\\
                4&\qquad$\oto\phi$&Identity\\
                5&$\phi\wedge(\phi\vee\psi)\oto\phi$&1--4: Transitivity
            \end{tabular}\\
            The second is proved similarly.
        \end{proof}

        \subsubsection*{Basic properties of implication}

        \begin{theorem}[Reflexivity of implication]
            \label{thrm:ReflexityOfImplication}
            $$\phi\to\phi.$$
        \end{theorem}

        \begin{proof}\;\\
            \begin{tabular}{lll}
                1&$\phi\wedge\phi\oto\phi$&Idempotence\\
                2&$\phi\to\phi$&1: Definition of implication
            \end{tabular}\\
        \end{proof}

        Antisymmetry of implication is the same as Bi-implication introduction (Rule \ref{derivedrule:Bi-implication introduction}).

        \begin{derivedrule}[Transitivity of implication]
            \label{derivedrule:Transitivity of implication}
            $$\{\phi\to\psi,\psi\to\chi\} \tto \phi\to\chi.$$
        \end{derivedrule}

        \begin{proof}\;\\
            \begin{tabular}{lll}
                1&$\phi\to\psi$&Premise\\
                2&$\psi\to\chi$&Premise\\
                3&$\phi\wedge\psi\oto\phi$&1: Definition of implication\\
                4&$\psi\wedge\chi\oto\psi$&2: Definition of implication\\
                &$\phi\wedge\chi$&\\
                5&\qquad$\oto\phi\wedge\psi\wedge\chi$&3: Replacement $[\phi\mapsto\phi\wedge\psi]$\\
                6&\qquad$\oto\phi\wedge\psi$&4: Replacement $[\psi\wedge\chi\mapsto\psi]$\\
                7&\qquad$\oto\phi$&3: Replacement\\
                8&$\phi\wedge\chi\oto\phi$&5--7: Transitivity of bi-implication\\
                9&$\phi\to\chi$&8: Definition of implication\\
            \end{tabular}\\
        \end{proof}

        \subsubsection*{Equivalent definitions of implication}

        \begin{lemma}
            \label{lemma:Absorption with negation}
            $$\phi\wedge(\neg\phi\vee\psi)\oto\phi\wedge\psi,\;\phi\vee(\neg\phi\wedge\psi)\oto\phi\vee\psi.$$
        \end{lemma}

        \begin{proof}\;\\
            \begin{tabular}{lll}
                &$\phi\wedge(\neg\phi\vee\psi)$&\\
                1&\qquad$\oto(\phi\wedge\neg\phi)\vee(\phi\wedge\psi)$&Distributivity\\
                2&\qquad$\oto\bot\vee(\phi\wedge\psi)$&LNC\\
                3&\qquad$\oto\phi\wedge\psi$&Identity\\
                4&$\phi\wedge(\neg\phi\vee\psi)\oto\phi\wedge\psi$&1--3: Transitivity\\
            \end{tabular}\\
            The second is proved dually.
        \end{proof}

        \begin{derivedrule}[Equivalent definitions of implication]
            \label{derivedrule:Equivalent definitions of implication}
            \begin{equation*}
                \begin{aligned}
                    &\phi\to\psi\\
                    \otto\;&\phi\wedge\psi\oto\phi\\
                    \otto\;&\phi\vee\psi\oto\psi\\
                    \otto\;&\phi\wedge\neg\psi\oto\bot\\
                    \otto\;&\neg\phi\vee\psi\oto\top.
                \end{aligned}
            \end{equation*}
        \end{derivedrule}

        \begin{proof}\;\\
            \begin{tabular}{lll}
                1&$\phi\to\psi$& Premise\\
                2&$\phi\wedge\psi\oto\phi$&1: Definition of implication\\
                & & \\
            \end{tabular}\\

            \begin{tabular}{lll}
                1&$\phi\wedge\psi\oto\phi$&Premise\\
                &$\phi\vee\psi$&\\
                2&\qquad$\oto(\phi\wedge\psi)\vee\psi$&1: Replacement $[\phi\mapsto\phi\wedge\psi]$\\
                3&\qquad$\oto\psi$&Absorption\\
                4&$\phi\vee\psi\oto\psi$&2--3: Transitivity\\
                & & \\
            \end{tabular}\\

            \begin{tabular}{lll}
                1&$\phi\vee\psi\oto\psi$&Premise\\
                &$\phi\wedge\neg\psi$&\\
                2&\qquad$\oto(\phi\wedge\neg\psi)\vee\bot$&Identity\\
                3&\qquad$\oto(\phi\wedge\neg\psi)\vee(\psi\wedge\neg\psi)$&LNC\\
                4&\qquad$\oto(\phi\vee\psi)\wedge\neg\psi$&Distributivity\\
                5&\qquad$\oto\psi\wedge\neg\psi$&1: Replacement $[\phi\vee\psi\mapsto\psi]$\\
                6&\qquad$\oto\bot$&LNC\\
                7&$\phi\wedge\neg\psi\oto\bot$&2--6: Transitivity\\
                & & \\
            \end{tabular}\\

            \begin{tabular}{lll}
                1&$\phi\wedge\neg\psi\oto\bot$&Premise\\
                &$\neg\phi\vee\psi$&\\
                2&\qquad$\oto\neg\phi\vee\bot\vee\psi$&Identity\\
                3&\qquad$\oto\neg\phi\vee(\phi\wedge\neg\psi)\vee\psi$&1: Replacement $[\bot\mapsto\phi\wedge\neg\psi]$\\
                4&\qquad$\oto\neg\phi\vee\phi\vee\psi$&Lemma \ref{lemma:Absorption with negation}\\
                5&\qquad$\oto\top\vee\psi$&LEM\\
                6&\qquad$\oto\top$&Annihilator\\
                7&$\neg\phi\vee\psi\oto\top$&2--6: Transitivity\\
                & & \\
            \end{tabular}\\

            \begin{tabular}{lll}
                1&$\neg\phi\vee\psi\oto\top$&Premise\\
                &$\phi\wedge\psi$&\\
                2&\qquad$\oto\bot\vee(\phi\wedge\psi)$&Identity\\
                3&\qquad$\oto(\phi\wedge\neg\phi)\vee(\phi\wedge\psi)$&LNC\\
                4&\qquad$\oto\phi\wedge(\neg\phi\vee\psi)$&Distributivity\\
                5&\qquad$\oto\phi\wedge\top$&1: Replacement\\
                6&\qquad$\oto\phi$&Identity\\
                7&$\phi\wedge\psi\oto\phi$&2--6: Transitivity\\
                & & \\
            \end{tabular}\\

            \begin{tabular}{lll}
                1&$\phi\wedge\psi\oto\phi$&Premise\\
                2&$\phi\to\psi$&1: Definition of implication\\
            \end{tabular}\\
        \end{proof}

        Hereinafter, any one of the four is simply called the definition of implication (by bi-implication).

        \subsubsection*{Implication-to-Disjunction}

        \begin{theorem}[Implication-to-Disjunction / Upper bound of implication]
            \label{thrm:Implication-to-Disjunction}
            $$(\phi\to\psi)\to\neg\phi\vee\psi.$$
        \end{theorem}

        \begin{proof}\;\\
            \begin{tabular}{lll}
                1&$\phi\to\psi$&Premise\\
                2&$\neg\phi\vee\psi\oto\top$& 1: Definition of implication\\
                3&$\top$& Top is a theorem\\
                4&$\neg\phi\vee\psi$& 2,3: Replacement $[\top\mapsto\neg\phi\vee\psi]$\\
            \end{tabular}\\
        \end{proof}

        \subsubsection*{De Morgan's laws}

        \begin{theorem}[De Morgan's laws]
            \label{thrm:De Morgan's laws}
            $$\neg(\phi\wedge\psi)\oto(\neg\phi\vee\neg\psi),\;\neg(\phi\vee\psi)\oto(\neg\phi\wedge\neg\psi).$$
        \end{theorem}

        \begin{proof}\;\\
            \begin{tabular}{lll}
                &$\neg\neg(\phi\wedge\psi)\vee(\neg\phi\vee\neg\psi)$&\\
                1&\qquad $\oto(\phi\wedge\psi)\vee(\neg\phi\vee\neg\psi)$&Double negation\\
                2&\qquad $\oto(\phi\vee\neg\phi\vee\neg\psi)\wedge(\psi\vee\neg\phi\vee\neg\psi)$&Distributivity\\
                3&\qquad $\oto(\top\vee\neg\psi)\wedge(\top\vee\neg\phi)$&LEM\\
                4&\qquad $\oto\top\wedge\top$&Annihilator\\
                5&\qquad $\oto\top$&Identity element\\
                6&$\neg\neg(\phi\wedge\psi)\vee(\neg\phi\vee\neg\psi)\oto\top$&1--5: Transitivity\\
                7&$\neg(\phi\wedge\psi)\to(\neg\phi\vee\neg\psi)$&6: Definition of implication\\
                &$(\neg\phi\vee\neg\psi)\wedge\neg\neg(\phi\wedge\psi)$&\\
                8&\qquad $\oto(\neg\phi\vee\neg\psi)\wedge(\phi\wedge\psi)$&Double negation\\
                9&\qquad $\oto(\neg\phi\wedge\phi\wedge\psi)\vee(\neg\psi\wedge\phi\wedge\psi)$&Distributivity\\
                10&\qquad $\oto(\bot\wedge\psi)\vee(\bot\wedge\phi)$&LNC\\
                11&\qquad $\oto\bot\vee\bot$&Annihilator\\
                12&\qquad $\oto\bot$&Identity element\\
                13&$(\neg\phi\vee\neg\psi)\wedge\neg\neg(\phi\wedge\psi)\oto\bot$&8--12: Transitivity\\
                14&$\neg\phi\vee\neg\psi\to\neg(\phi\wedge\psi)$&13: Definition of implication\\
                15&$\neg(\phi\wedge\psi)\oto(\neg\phi\vee\neg\psi)$&7,14: Bi-implication introduction\\
            \end{tabular}\\
            The second is proved dually.
        \end{proof}

        \subsubsection*{Lattice-like properties}

        \begin{theorem}[Conjunction elimination and disjunction introduction]
            \label{thrm:Conjunction elimination and disjunction introduction}
            $$\phi\wedge\psi\to\phi,\;\phi\to\phi\vee\psi.$$
        \end{theorem}

        \begin{proof}
            In terms of equivalent definitions of implication, conjunction elimination and disjunction introduction are equivalent to absorption laws $\phi\vee(\phi\wedge\psi)\oto\phi$ and $\phi\wedge(\phi\vee\psi)\oto\phi$ respectively.
        \end{proof}

        \begin{derivedrule}[One-to-two and Two-to-one]
            \label{derivedrule:One-to-two and Two-to-one}
            $$\{\chi\to\phi,\chi\to\psi\}\otto \chi\to\phi\wedge\psi,\;\{\phi\to\chi,\psi\to\chi\}\otto \phi\vee\psi\to\chi.$$
        \end{derivedrule}

        \begin{proof}\;\\
            \begin{tabular}{lll}
                1&$\chi\to\phi$&Premise\\
                2&$\chi\to\psi$&Premise\\
                3&$\chi\wedge\phi\oto\chi$&1: Definition of implication\\
                4&$\chi\wedge\psi\oto\chi$&2: Definition of implication\\
                5&$\chi\wedge\phi\wedge\psi\oto\chi$&3,4: Replacement $[\chi\mapsto\chi\wedge\phi]$\\
                6&$\chi\to\phi\wedge\psi$&5: Definition of implication\\
                & & \\
                1&$\chi\to\phi\wedge\psi$&Premise\\
                2&$\phi\wedge\psi\to\phi$&Conjunction elimination\\
                3&$\chi\to\phi$&1--2: Transitivity\\
                4&$\chi\to\psi$&1--3: Similarly\\
            \end{tabular}\\
            The second is proved dually.
        \end{proof}

        \begin{derivedrule}[Strengthening]
            \label{derivedrule:Strengthening}
            $$\phi\to\psi\tto \phi\to\phi\wedge\psi,\;\phi\to\psi\tto \phi\vee\psi\to\psi.$$
        \end{derivedrule}

        \begin{proof}\;\\
            \begin{tabular}{lll}
                1&$\phi\to\psi$&Premise\\
                2&$\phi\to\phi$&Reflexivity\\
                3&$\phi\to\phi\wedge\psi$&1,2: One-to-two\\
            \end{tabular}\\
            The second is proved dually.
        \end{proof}

        \begin{derivedrule}[Weakening]
            \label{derivedrule:Weakening}
            $$\phi\to\psi\tto \phi\wedge\chi\to\psi,\;\phi\to\psi\tto \phi\to\psi\vee\chi.$$
        \end{derivedrule}

        \begin{proof}\;\\
            \begin{tabular}{lll}
                1&$\phi\to\psi$&Premise\\
                2&$\phi\wedge\chi\to\phi$&Conjunction elimination\\
                3&$\phi\wedge\chi\to\psi$&2,1: Transitivity\\
            \end{tabular}\\
            The second is proved dually.
        \end{proof}

        \begin{derivedrule}[Compatibility of implication with conjunction and disjunction]
            $$\phi\to\psi\tto \phi\wedge\chi\to\psi\wedge\chi,\;\phi\to\psi\tto \phi\vee\chi\to\psi\vee\chi.$$
        \end{derivedrule}

        \begin{proof}\;\\
            \begin{tabular}{lll}
                1&$\phi\to\psi$&Premise\\
                2&$\phi\wedge\chi\to\psi$&1: Weakening\\
                3&$\phi\wedge\chi\to\chi$&Conjunction elimination\\
                4&$\phi\wedge\chi\to\psi\wedge\chi$&2,3: One-to-two\\
            \end{tabular}\\
            The second is proved dually.
        \end{proof}

        \begin{derivedrule}[Monotonicity of conjunction and disjunction w.r.t. implication]
            \begin{equation*}
                \begin{aligned}
                    &\{\phi_1\to\psi_1,\phi_2\to\psi_2\}\tto \phi_1\wedge\phi_2\to\psi_1\wedge\psi_2,\\
                    &\{\phi_1\to\psi_1,\phi_2\to\psi_2\}\tto \phi_1\vee\phi_2\to\psi_1\vee\psi_2.
                \end{aligned}
            \end{equation*}
        \end{derivedrule}

        \begin{proof}\;\\
            \begin{tabular}{lll}
                1&$\phi_1\to\psi_1$&Premise\\
                2&$\phi_2\to\psi_2$&Premise\\
                3&$\phi_1\wedge\phi_2\to\psi_1$&1: Weakening\\
                4&$\phi_1\wedge\phi_2\to\psi_2$&2: Weakening\\
                5&$\phi_1\wedge\phi_2\to\psi_1\wedge\psi_2$&3,4: One-to-two\\
            \end{tabular}\\
            The second is proved dually.
        \end{proof}

        \begin{theorem}[Contraposition]
            \label{thrm:Contraposition}
            $$(\phi\to\psi)\oto(\neg\psi\to\neg\phi).$$
        \end{theorem}

        \begin{proof}\;\\
            \begin{tabular}{lll}
                &$(\phi\to\psi)$&\\
                1&\qquad $\oto(\phi\wedge\neg\psi\oto\bot)$& Definition of implication\\
                2&\qquad $\oto(\neg\psi\wedge\neg\neg\phi\oto\bot)$& Commutativity, Double negation\\
                3&\qquad $\oto(\neg\psi\to\neg\phi)$& Definition of implication\\
                4&$(\phi\to\psi)\oto(\neg\psi\to\neg\phi)$& 1--3: Transitivity
            \end{tabular}\\
        \end{proof}

        \begin{derivedrule}[Modus tollens]
            \label{derule:Modus tollens}
            $$\{\phi\to\psi,\neg\psi\}\tto\neg\phi.$$
        \end{derivedrule}

        \begin{proof}\;\\
            \begin{tabular}{lll}
                1&$\phi\to\psi$&Premise\\
                2&$\neg\psi$&Premise\\
                3&$\neg\psi\to\neg\phi$&1: Contraposition\\
                4&$\neg\phi$&2,3: Modus ponens
            \end{tabular}\\
        \end{proof}

        \begin{theorem}[ECQ and its dual]
            \label{thrm:ECQ and its dual}
            $$\bot\to\phi,\;\phi\to\top.$$
        \end{theorem}

        \begin{proof}
            In terms of equivalent definitions of implication, ECQ and its dual are equivalent to Annihilators $\bot\wedge\phi\oto\bot$ and $\top\vee\phi\oto\top$ respectively.
        \end{proof}

        \subsubsection*{Some properties of bi-implication}

        \begin{derivedrule}[Bi-implication elimination]
            $$\phi\oto\psi\tto \phi\to\psi.$$
        \end{derivedrule}

        \begin{proof}\;\\
            \begin{tabular}{lll}
                1&$\phi\oto\psi$&Premise\\
                2&$\phi\to\phi$&Reflexivity of implications\\
                3&$\phi\to\psi$&1,2: Replacement\\
            \end{tabular}\\
        \end{proof}

        \begin{derivedrule}[Bi-implication by implications]
            $$\phi\oto\psi\otto(\phi\to\psi)\wedge(\psi\to\phi).$$
        \end{derivedrule}

        \begin{proof}\;\\
            \begin{tabular}{lll}
                1&$\phi\oto\psi$&Premise\\
                2&$\phi\to\psi$&1: Bi-implication elimination\\
                3&$\psi\to\phi$&1: Bi-implication elimination\\
                4&$(\phi\to\psi)\wedge(\psi\to\phi)$&2,3: Conjunction introduction\\
                &&\\
                1&$(\phi\to\psi)\wedge(\psi\to\phi)$&Premise\\
                2&$\phi\oto\psi$&1: Bi-implication introduction\\
            \end{tabular}\\
        \end{proof}

        \begin{derivedrule}[Bi-implication of negation]
            $$\phi\oto\psi\otto \neg\phi\oto\neg\psi.$$
        \end{derivedrule}

        \begin{proof}\;\\
            \begin{tabular}{lll}
                &$(\phi\oto\psi)$&\\
                1&\qquad$\oto(\phi\to\psi)\wedge(\psi\to\phi)$&Bi-implication by implications\\
                2&\qquad$\oto(\neg\psi\to\neg\phi)\wedge(\neg\phi\to\neg\psi)$&Contraposition\\
                3&\qquad$\oto(\neg\phi\oto\neg\psi)$&Bi-implication by implications
            \end{tabular}\\
        \end{proof}

        \subsubsection*{More properties of implication}

        \begin{theorem}[Generalized ECQ of implication]
            \label{thrm:Generalized ECQ of implication}
            $$(\phi\to\psi)\to(\chi\to(\phi\to\psi)).$$
        \end{theorem}

        \begin{proof}\;\\
            \begin{tabular}{lll}
                1&$\phi\to\psi$&Premise\\
                2&$\phi\wedge\psi\oto\phi$& 1: Definition of implication\\
                3&$\phi\wedge\psi\to\psi$&Conjunction elimination\\
                4&$\chi\to(\phi\wedge\psi\to\psi)$&3: Theorem like top\\
                5&$\chi\to(\phi\to\psi)$&2,4: Replacement $[\phi\wedge\psi\mapsto\phi]$\\
            \end{tabular}\\
            It can also be derived in a meta proof:\\
            \begin{tabular}{lll}
                1 & $\vdash\phi\wedge\psi\to\psi$ & Conjunction elimination\\
                2 & $\vdash\chi\to(\phi\wedge\psi\to\psi)$ & 1: Theorem like top\\
                3 & \qquad $\phi\to\psi$ & Premise\\
                4 & \qquad $\phi\wedge\psi\oto\phi$& 3: Definition of implication\\
                5 & \qquad $\chi\to(\phi\wedge\psi\to\psi)$ & 2 (theorem from an outer layer)\\
                6 & \qquad $\chi\to(\phi\to\psi)$ & 4,5: Replacement $[\phi\wedge\psi\mapsto\phi]$\\
                7 & $\vdash(\phi\to\psi)\to(\chi\to(\phi\to\psi))$ & 3--6: Conditional proof\\
            \end{tabular}\\
        \end{proof}

        Note that Generalized ECQ of the general form $\psi\to(\phi\to\psi)$ does not hold for IRL (cf. Subsection \ref{subsec:D2I is underivable in IRL}).

        \begin{theorem}[Safe generalized ECQ]
            \label{thrm:Safe generalized ECQ}
            $$(\top\to\phi)\to(\psi\to\phi),\;(\phi\to\bot)\to(\phi\to\psi).$$
        \end{theorem}
        \begin{proof}\;\\
            \begin{tabular}{lll}
                1&$\top\to\phi$ & Premise\\
                2&$\psi\to\top$ & Dual of ECQ\\
                3&$\psi\to\phi$ & 2,1: Transitivity
            \end{tabular}\\
            The second is proved dually.\\
        \end{proof}

        \begin{theorem}[Importation]
            \label{thrm:Importation}
            $$(\phi\to(\psi\to\chi))\to(\phi\wedge\psi\to\chi).$$
        \end{theorem}

        \begin{proof}\;\\
            \begin{tabular}{lll}
                1&$\phi\to(\psi\to\chi)$&Premise\\
                &$\phi\wedge\psi$&\\
                2&\qquad$\to(\psi\to\chi)\wedge\psi$&1: Compatibility of $\to$ with $\wedge$\\
                3&\qquad$\to\chi$&Modus ponens\\
                4&$\phi\wedge\psi\to\chi$&2--3: Transitivity
            \end{tabular}\\
        \end{proof}

        \begin{corollary}[Contraction]
            \label{corollary:Contraction}
            $(\phi\to(\phi\to\psi))\to(\phi\to\psi)$.
        \end{corollary}
        \begin{proof}
            It follows from the above theorem and Idempotence.
        \end{proof}

        Exportation (converse of importation) $ (\phi\wedge\psi\to\chi)\to(\phi\to(\psi\to\chi))$ does not hold for IRL except for some special cases as follows.

        \begin{theorem}[Exported antisymmetry of implication / bi-implication introduction]
            \label{thrm:Exported antisymmetry of implication}
            $$(\phi\to\psi)\to((\psi\to\phi)\to(\phi\oto\psi)).$$
        \end{theorem}

        \begin{proof}\;\\
            \begin{tabular}{lll}
                1&$\phi\to\psi$&Premise\\
                2&$(\psi\to\phi)\to(\phi\to\psi)$&1: Generalized ECQ of implication\\
                3&$(\psi\to\phi)\to(\phi\to\psi)\wedge(\psi\to\phi)$&2: Strengthening\\
                4&$(\phi\to\psi)\wedge(\psi\to\phi)\to(\phi\oto\psi)$&Bi-implication introduction\\
                5&$(\psi\to\phi)\to(\phi\oto\psi)$&3,4: Transitivity\\
            \end{tabular}\\
        \end{proof}

        \begin{theorem}[Exported transitivity of implication / Suffixing and Prefixing]
            \label{thrm:Exported transitivity of implication}
            $$(\phi\to\psi)\to((\psi\to\chi)\to(\phi\to\chi)),\;(\phi\to\psi)\to((\chi\to\phi)\to(\chi\to\psi)).$$
        \end{theorem}

        \begin{proof}\;\\
            \begin{tabular}{lll}
                1&$\phi\to\psi$&Premise\\
                2&$(\psi\to\chi)\to(\phi\to\psi)$&1: Generalized ECQ of implication\\
                3&$(\psi\to\chi)\to(\phi\to\psi)\wedge(\psi\to\chi)$&2: Strengthening\\
                4&$(\phi\to\psi)\wedge(\psi\to\chi)\to(\phi\to\chi)$&Transitivity of implication\\
                5&$(\psi\to\chi)\to(\phi\to\chi)$&3--4: Transitivity of implication\\
            \end{tabular}\\
            The second is obtained via uniform substitution and contraposition.\\
        \end{proof}

        \begin{theorem}[Exported one-to-two and two-to-one]
            \label{thrm:Exported one-to-two and two-to-one}
            $$(\chi\to\phi)\to((\chi\to\psi)\to(\chi\to\phi\wedge\psi)),\;(\phi\to\chi)\to((\psi\to\chi)\to(\phi\vee\psi\to\chi)).$$
        \end{theorem}

        \begin{proof}\;\\
            \begin{tabular}{lll}
                1&$\chi\to\phi$&Premise\\
                2&$(\chi\to\psi)\to(\chi\to\phi)$&1: Generalized ECQ of implication\\
                3&$(\chi\to\psi)\to(\chi\to\phi)\wedge(\chi\to\psi)$&2: Strengthening\\
                4&$(\chi\to\phi)\wedge(\chi\to\psi)\to(\chi\to\phi\wedge\psi)$&One-to-two\\
                5&$(\chi\to\psi)\to(\chi\to\phi\wedge\psi))$&3--4: Transitivity\\
            \end{tabular}\\
            The second is obtained via contraposition.\\
        \end{proof}

        \begin{theorem}[Self-distributivity of implication]
            \label{thrm:Self-distributivity of implication}
            $$(\chi\to(\phi\to\psi))\to((\chi\to\phi)\to(\chi\to\psi)).$$
        \end{theorem}

        \begin{proof}\;\\
            \begin{tabular}{lll}
                1&$\chi\to(\phi\to\psi)$& Premise\\
                2&$\chi\wedge\phi\to\psi$& 1: Importation\\
                3&$(\chi\to\phi)\to(\chi\wedge\phi\to\psi)$& 2: Generalized ECQ of implication\\
                4&$(\chi\to\phi)\to(\chi\to\chi\wedge\phi)$& Strengthening\\
                5&$(\chi\to\phi)\to(\chi\to\chi\wedge\phi)\wedge(\chi\wedge\phi\to\psi)$& 4,3: One-to-two\\
                6&$(\chi\to\chi\wedge\phi)\wedge(\chi\wedge\phi\to\psi)\to(\chi\to\psi)$& Transitivity\\
                7&$(\chi\to\phi)\to(\chi\to\psi)$& 5,6: Transitivity\\
            \end{tabular}\\
        \end{proof}

        \subsubsection*{Implication with one or no variable}

        \begin{theorem}[Univariate implication expressions]
            \label{thrm:Univariate implication expressions}
            \begin{equation*}
                \begin{aligned}
                    &(\phi\to\neg\phi) \oto (\phi\oto\bot) \oto (\phi\to\bot);\\
                    &(\neg\phi\to\phi) \oto (\phi\oto\top) \oto (\top\to\phi).
                \end{aligned}
            \end{equation*}
        \end{theorem}

        \begin{proof}\;\\
            \begin{tabular}{lll}
                &$(\phi\to\neg\phi)$&\\
                1&\qquad$\oto(\phi\wedge\neg\neg\phi\oto\bot)$&Definition of implication\\
                2&\qquad$\oto(\phi\wedge\phi\oto\bot)$&Double negation\\
                3&\qquad$\oto(\phi\oto\bot)$&Idempotence\\
                4&\qquad$\oto(\phi\to\bot)\wedge(\bot\to\phi)$&Bi-implication by implication\\
                5&\qquad$\oto(\phi\to\bot)$&Theorem (ECQ) like top\\
                6&$(\phi\to\neg\phi) \oto (\phi\oto\bot) \oto (\phi\to\bot)$&1--3--5: Transitivity\\
            \end{tabular}\\
            The second is proved dually.
        \end{proof}

        \begin{theorem}[Upper bound of univariate implication]
            \label{thrm:Upper bound of univariate implication}
            $$(\top\to\phi)\to\phi,\;
            (\phi\to\bot)\to\neg\phi.$$
        \end{theorem}

        \begin{proof}
            They follow from the upper bound of implication $(\phi\to\psi)\to\neg\phi\vee\psi$ (Theorem \ref{thrm:Implication-to-Disjunction}) by uniform substitution, and properties of top and bottom.
        \end{proof}

        \begin{corollary}[Contradictory implication]
            \label{corollary:Contradictory implication}
            $$(\top\to\bot)\to\bot.$$
        \end{corollary}

        \begin{proof}
            It follows from Theorem \ref{thrm:Upper bound of univariate implication} by uniform substitution.
        \end{proof}

        \subsubsection*{Proof methods related}

        \begin{derivedrule}[Constructive dilemma]
            \label{derivedrule:Constructive dilemma}
            $$\{\phi_1\vee\phi_2,\phi_1\to\psi_1,\phi_2\to\psi_2\}\tto\psi_1\vee\psi_2.$$
        \end{derivedrule}

        \begin{proof}\;\\
            \begin{tabular}{lll}
                1&$\phi_1\vee\phi_2$&Premise\\
                2&$\phi_1\to\psi_1$&Premise\\
                3&$\phi_2\to\psi_2$&Premise\\
                4&$\phi_1\vee\phi_2\to\psi_1\vee\psi_2$&2,3: Monotonicity of $\vee$ w.r.t. $\to$\\
                5&$\psi_1\vee\psi_2$&1,4: Modus penens\\
            \end{tabular}\\
        \end{proof}

        \begin{derivedrule}[Proof by exhaustion]
            \label{derivedrule:Proof by exhaustion}
            $$\{\phi\vee\psi,\phi\to\chi,\psi\to\chi\}\tto\chi.$$
        \end{derivedrule}

        \begin{proof}\;\\
            \begin{tabular}{lll}
                1&$\phi\vee\psi$& Premise\\
                2&$\phi\to\chi$& Premise\\
                3&$\psi\to\chi$& Premise\\
                4&$\phi\vee\psi\to\chi$& 2,3: Two-to-one\\
                5&$\chi$& 1,4: Modus ponens
            \end{tabular}\\
        \end{proof}

        \begin{derivedrule}[LEM exhaustion]
            \label{derivedrule:LEM exhaustion}
            $$\{\phi\to\psi,\neg\phi\to\psi\}\tto\psi.$$
        \end{derivedrule}

        \begin{proof}\;\\
            \begin{tabular}{lll}
                1&$\phi\to\psi$& Premise\\
                2&$\neg\phi\to\psi$& Premise\\
                3&$\psi\vee\neg\phi$& LEM\\
                4&$\psi$& 3,1,2: Proof by exhaustion
            \end{tabular}\\
        \end{proof}

        \begin{derivedrule}[Reductio ad absurdum (RAA)]
            \label{derivedrule:RAA}
            $$\{\phi\to\psi,\phi\to\neg\psi\}\to\neg\phi.$$
        \end{derivedrule}

        \begin{proof}\;\\
            \begin{tabular}{lll}
                1&$\phi\to\psi$&Premise\\
                2&$\phi\to\neg\psi$&Premise\\
                3&$\neg\psi\to\neg\phi$&1: Contraposition\\
                4&$\psi\to\neg\phi$&2: Contraposition, Double negation\\
                5&$\neg\phi$& 3,4: LEM exhaustion
            \end{tabular}\\
        \end{proof}

        The rules $\{\phi\to\psi,\phi\to\neg\psi\}\tto\neg\phi$, $\phi\to\neg\phi\tto\neg\phi$ and $\phi\to\bot\tto\neg\phi$ are all forms of RAA, while the rules $\{\neg\phi\to\psi,\neg\phi\to\neg\psi\}\tto\phi$, $\neg\phi\to\phi\tto\phi$ (consequentia mirabilis) and $\neg\phi\to\bot\tto\phi$ are forms of proof by contradiction.

        \begin{derivedrule}[Disjunctive syllogism]
            \label{derivedrule:Disjunctive syllogism}
            $$\{\phi\vee\psi,\neg\phi\}\tto\psi.$$
        \end{derivedrule}

        \begin{proof}\;\\
            \begin{tabular}{lll}
                1&$\phi\vee\psi$& Premise\\
                2&$\neg\phi$& Premise\\
                3&$(\phi\vee\psi)\wedge\neg\phi$& 1,2: Conjunction introduction\\
                4&$(\phi\wedge\neg\phi)\vee(\psi\wedge\neg\phi)$& 3: Distributivity\\
                5&$\phi\wedge\neg\phi\to\psi$& ECQ\\
                6&$\psi\wedge\neg\phi\to\psi$& Conjunction elimination\\
                7&$\psi$& 4,5,6: Proof by exhaustion
            \end{tabular}\\
        \end{proof}

        \begin{derivedrule}[Implication and inconsistency]
            \label{derivedrule:Implication and inconsistency}
            $$\phi\to\psi\otto \phi\wedge\neg\psi\to\bot.$$
        \end{derivedrule}

        \begin{proof}\;\\
            \begin{tabular}{lll}
                &$(\phi\to\psi)$& \\
                1&\qquad$\oto(\phi\wedge\neg\psi\oto\bot)$& Definition of implication\\
                2&\qquad$\oto(\phi\wedge\neg\psi\to\bot)$& Univariate implication expressions\\
                3&$(\phi\to\psi)\oto (\phi\wedge\neg\psi\to\bot)$& 1,2: Transitivity\\
            \end{tabular}\\
        \end{proof}

        \begin{derivedrule}[LEM exhaustion for implication]
            \label{derivedrule:LEM exhaustion for implication}
            $$\{\chi\wedge\phi\to\psi,\neg\chi\wedge\phi\to\psi\}\tto\phi\to\psi.$$
        \end{derivedrule}

        \begin{proof}\;\\
            \begin{tabular}{lll}
                1&$\chi\wedge\phi\to\psi$& Premise\\
                2&$\neg\chi\wedge\phi\to\psi$& Premise\\
                3&$(\chi\wedge\phi)\vee(\neg\chi\wedge\phi)\to\psi$&  1,2: Two-to-one\\
                &$(\chi\wedge\phi)\vee(\neg\chi\wedge\phi)$& \\
                4&\qquad$\oto(\chi\vee\neg\chi)\wedge\phi$& Distributivity\\
                5&\qquad$\oto\phi$& Theorem (LEM) like top\\
                6&$(\chi\wedge\phi)\vee(\neg\chi\wedge\phi)\oto\phi$& 4--5: Transitivity\\
                7&$\phi\to\psi$&6,3: Replacement
            \end{tabular}\\
        \end{proof}

        \begin{derivedrule}[Resolution]
            \label{derivedrule:Resolution}
            $$\{\phi\vee\chi,\psi\vee\neg\chi\}\tto\phi\vee\psi.$$
        \end{derivedrule}

        \begin{proof}\;\\
            \begin{tabular}{lll}
                1&$\neg\chi\wedge(\phi\vee\chi)\to\phi$& Disjunctive syllogism\\
                2&$\chi\wedge(\psi\vee\neg\chi)\to\psi$& Disjunctive syllogism\\
                3&$\neg\chi\wedge(\phi\vee\chi)\to\phi\vee\psi$& 1: Strengthening\\
                4&$\chi\wedge(\psi\vee\neg\chi)\to\phi\vee\psi$& 2: Strengthening\\
                5&$\neg\chi\wedge(\phi\vee\chi)\wedge(\psi\vee\neg\chi)\to\phi\vee\psi$& 3: Weakening\\
                6&$\chi\wedge(\phi\vee\chi)\wedge(\psi\vee\neg\chi)\to\phi\vee\psi$& 4: Weakening\\
                7&$(\phi\vee\chi)\wedge(\psi\vee\neg\chi)\to\phi\vee\psi$& 5,6: LEM exhaustion for implication
            \end{tabular}\\
        \end{proof}

        \begin{derivedrule}[Cut]
            \label{thrm:Cut}
            $$\{\phi_1\to\chi,\phi_2\wedge\chi\to\psi\}\tto \phi_1\wedge\phi_2\to\psi.$$
        \end{derivedrule}

        \begin{proof}\;\\
            \begin{tabular}{lll}
                1&$\phi_1\to\chi$& Premise\\
                2&$\phi_2\wedge\chi\to\psi$& Premise\\
                3&$\phi_1\wedge\phi_2\to\chi$& 1: Weakening\\
                4&$\phi_1\wedge\phi_2\to\phi_1\wedge\phi_2\wedge\chi$& 3: Strengthening\\
                5&$\phi_1\wedge\phi_2\wedge\chi\to\psi$& 2: Weakening\\
                6&$\phi_1\wedge\phi_2\to\psi$& 4,5: Transitivity\\
            \end{tabular}\\
        \end{proof}

        \begin{metarule}[Equivalence of implication theorem and disjunction theorem]
            \label{metarule:Equivalence of implication theorem and disjunction theorem}
            $$(\vdash\phi\to\psi)\otto(\vdash\neg\phi\vee\psi).$$
        \end{metarule}

        \begin{proof}\;\\
            \begin{tabular}{lll}
                1&$\vdash\phi\to\psi$&Premise\\
                2&\qquad $(\phi\to\psi)\to\neg\phi\vee\psi$&Implication-to-Disjunction\\
                3&\qquad $\phi\to\psi$&1 (assumed theorem)\\
                4&\qquad $\neg\phi\vee\psi$&2,3: Modus ponens\\
                5&$\vdash\neg\phi\vee\psi$&2--4: Unconditional proof\\
                &&\\
                1&$\vdash\neg\phi\vee\psi$&Premise\\
                2&\qquad $\phi$&Premise\\
                3&\qquad $\neg\phi\vee\psi$&1 (assumed theorem)\\
                4&\qquad $\psi$&2,3: Disjunctive syllogism\\
                5&$\vdash\phi\to\psi$&3--5: Conditional proof
            \end{tabular}\\
        \end{proof}


        \section{Semantics and the algebra IRLA}
        \label{sec:Semantics and the algebra IRLA}

        \subsection{Semantics of IRL}
        \label{subsec:Semantics of IRL}

        An interpretation of a formal language assigns a meaning (semantic value) to each expression within that language. However, logic is fundamentally concerned with truth values rather than specific contents of statements of a language. So interpretations that assign the same truth value (even if they assign different contents) to a statement can be viewed as the same in a logic. Hence, an interpretation of a propositional language can be represented by just a truth assignment or valuation function.

        \vspace{\baselineskip}
        Let $\{\mathrm{T,F}\}$ be the set of truth values for a two-valued logic.

        Let $L$ be the set of all formulas in the language of IRL.

        Let $U = \{\mathrm{T,F}\}^L$, the set of all valuation functions for $L$.

        The set $L$ is countably infinite (even in a case when the countable set of primitive propositions is finite), thus $|U|=|\{\mathrm{T,F}\}|^{|L|}=2^{|L|} =2^{\aleph_0}=\mathfrak{c}$, i.e. $U$ is generally uncountable with the cardinality of $\mathbb{R}$. (In classical logic, as all the operations -- including particularly the implication -- are truth-functional, the number of all interpretations for a formula with $n$ propositional variables is only $2^n$.)

        \begin{definition}[T-set]
            \label{def:T-set}
            The T-set of a formula $\phi\in L$ is defined as
            $$V_{\phi} = \{v\in U \mid v(\phi)=\mathrm{T}\}.$$
        \end{definition}

        \subsubsection*{Semantics of IRL}
        \label{subsubsec:Semantics of IRL}

        The semantic system for IRL is defined via T-set as follows.

        \begin{enumerate}
            \item The structure $(2^U,U,\emptyset,{}^\complement,\cap,\cup,\subseteq)$ is the power set Boolean lattice of $U$. (For any formula $\phi$ its T-set $V_{\phi}$ is a member of $2^U$ by definition.)
            \item For any primitive proposition $A$: $V_A$ is given (may be unknown, so underivable).
            \item For the truth-functional operations $\{\top,\bot,\neg,\wedge,\vee\}$:
            \subitem $V_{\top}=U$, $V_{\bot}=\emptyset$.
            \subitem $V_{\neg\phi} = V_{\phi}^{\complement}$.
            \subitem $V_{\phi\wedge\psi} = V_{\phi}\cap V_{\psi}$, $V_{\phi\vee\psi} = V_{\phi}\cup V_{\psi}$.
            \item For the non-truth-functional bi-implication $\oto$:
            \subitem $V_{\phi\oto\phi} = U$, $V_{\phi\oto\psi} = V_{\psi\oto\phi}$, $V_{\phi\oto\psi}\cap V_{\psi\oto\chi} \subseteq V_{\phi\oto\chi}$;
            \subitem $V_{\phi\oto\psi}\cap V_{\chi}\subseteq V_{\chi[\phi\mapsto \psi]}$;
            \subitem $V_{\phi\oto\psi} = U$ if and only if $V_{\phi} = V_{\psi}$.
            \item For the non-truth-functional implication $\to$: $V_{\phi\to\psi} = V_{\phi\wedge\psi\oto\phi}$.
        \end{enumerate}

        \begin{remark}\;
            \begin{itemize}
                \item Semantics of bi-implication, as it is not truth-functional, is given by a collection of characterizing conditions. Implication is defined by bi-implication.
                \item A derived result for implication is $V_{\phi\to\psi} = U \text{ if and only if } V_{\phi}\subseteq V_{\psi}$. In terms of the original meaning of implication, this indicates that $V_{\phi}\subseteq V_{\psi}$ is just the ``certain mechanism'' of the tautological implication from $\phi$ to $\psi$ (the scope of this mechanism is $U$).
                \item By definition of T-set, for any formula $\phi$, $\emptyset\subseteq V_{\phi}\subseteq U$, thus $V_{\bot\to\psi}=U$ and $V_{\phi\to\top}=U$, indicating that ECQ and its dual have ``certain mechanisms'' of tautological implications, so they do not conflict with the original meaning of implication.
            \end{itemize}
        \end{remark}

        \begin{definition}[Validity and semantic consequence]
            \label{def:Validity and semantic consequence}
            With respect to the semantic system, validity/tautology, semantic consequence and equivalence, are defined below:
            \begin{itemize}
                \item A formula $\phi$ is valid or is a \emph{tautology}, denoted $\vDash\phi$, if and only if $V_{\phi}=U$.
                \item A formula $\psi$ is a \emph{semantic consequence} of a formula $\phi$, denoted $\phi\vDash\psi$, if and only if $V_{\phi}\subseteq V_{\psi}$.
                \item Formulas $\phi$ and $\psi$ are semantically equivalent if  $\phi\vDash\psi$ and $\psi\vDash\phi$, or equivalently $V_{\phi}=V_{\psi}$.
            \end{itemize}
        \end{definition}

        The semantic system of IRL defined in this subsection can be viewed as an extension of the power set Boolean lattice of $U$, $(2^U,U,\emptyset,{}^\complement,\cap,\cup,\subseteq)$, by adding semantic characteristics of propositional formulas of IRL in the T-set notation. The algebraic system IRLA developed in the next subsection, an extension of the Boolean lattice $(X,1,0,\neg,\wedge,\vee,\le)$, is simply a ``concise version'' of the semantic system via T-set in this subsection. Symbols of IRLA language are defined in Table \ref{tab:Definition of IRLA symbols}.

        \begin{table}[H] 
            \centering
            \small
            \setlength{\tabcolsep}{8pt} 
            \renewcommand{\arraystretch}{1.28} 
            \begin{threeparttable}
                \caption{Definition of IRLA symbols}
                \label{tab:Definition of IRLA symbols}
                \begin{tabular}{c|cc}
                    \hline
                    \emph{IRL} & \emph{T-set} & \emph{IRLA} \\
                    \hline
                    Primitive propositions: & & Constants:\\
                    $A,B,...$ & $V_A,V_B,...$ & $a,b,...$\\
                    \hline
                    Formulas: & & Variables:\\
                    $\phi,\psi,...$ & $V_{\phi},V_{\psi},...$ & $x,y,...$\\
                    \hline
                    Formula schema: & & Term:\\
                    $F(\phi,\psi,...)$ & $V_{F(\phi,\psi,...)}$ & $t(x,y,...)$\\
                    \hline
                    Operations other than $\oto,\to$: & & \\
                    $\top,\bot,\neg,\wedge,\vee$ & $U,\emptyset,\cdot^\complement,\cap,\cup$ & $1,0,\neg,\wedge,\vee$\\
                    \hline
                    Bi-implication: & & \\
                    $\phi\oto\psi$ & $V_{\phi\oto\psi}$ \tnote{a} & $x\oto y$ \tnote{b}\\
                    \hline
                    Implication: & & \\
                    $\phi\to\psi$ & $V_{\phi\to\psi}$ \tnote{a} & $x\to y$ \tnote{b}\\
                    \hline
                    Tautology: & & \\
                    $\vDash\phi$ & $V_{\phi}=U$ & $x=1$\\
                    \hline
                    Contradiction: & & \\
                    $\vDash\neg\phi$ & $V_{\phi}=\emptyset$ & $x=0$\\
                    \hline
                    Contingency: & & \\
                    $\nvDash\phi$ and $\nvDash\neg\phi$ & $V_{\phi}\ne U$ and $V_{\phi}\ne \emptyset$ & $x\ne 1$ and $x\ne 0$\\
                    \hline
                    Tautological bi-implication: & & \\
                    $\vDash\phi\oto\psi$ & $V_{\phi\oto\psi} =U$ i.e. $V_{\phi} = V_{\psi}$ & $x\oto y=1$ i.e. $x=y$\\
                    \hline
                    Tautological implication: & & \\
                    $\vDash\phi\to\psi$ & $V_{\phi\to\psi}=U$ i.e. $V_{\phi}\subseteq V_{\psi}$ & $x\to y=1$ i.e. $x\le y$\\
                    \hline
                \end{tabular}
                \begin{tablenotes}
                    \item[a] $V_{\phi\oto\psi}$ and $V_{\phi\to\psi}$ can not be generally represented by $V_{\phi}$ and $V_{\psi}$.
                    \item[b] The IRLA representation of the T-sets $x=V_{\phi}$, $y=V_{\psi}$ and $x\oto y=V_{\phi\oto\psi}$ implies the definition $V_{\phi}\oto V_{\psi}=V_{\phi\oto\psi}$, and similarly $V_{\phi}\to V_{\psi}=V_{\phi\to\psi}$.
                \end{tablenotes}
            \end{threeparttable}
        \end{table}

        \subsection{Algebra IRLA}
        \label{subsec:Algebra IRLA}

        The algebra presented in this subsection, called \textbf{IRLA}, is an extension of Boolean lattice that is an extension of Boolean algebra.

        Let $(X,1,0,\neg,\wedge,\vee)$ be an \emph{algebraic system}, where $X$ is the underlying set, $(1,0,\neg,\wedge,\vee)$ are operations of arities $(0,0,1,2,2)$ on $X$. Note that ``an operation on a set'' implies that the set is \emph{closed} under the operation. So $1,0\in X$, i.e. $X$ contains at least two constants $1$ and $0$. And $X$ may contains other countably many constants $a,b,...$.

        Let $x,y,...$ be variables and $s,t,...$ terms for the algebraic language.

        \subsubsection*{Boolean algebra}
        \label{subsubsec:Boolean algebra}

        A Boolean algebra $(X,1,0,\neg,\wedge,\vee)$ is defined by the following axioms and rules of inference.

        \begin{itemize}
            \item Axioms: commutativity, associativity and distributivity for $\wedge$ and $\vee$; definitions of identity elements $1$ and $0$; definition of complement or negation $\neg$. Associativity is redundant (\cite{huntington}).

            \item Rules of inference:
            \subitem Primitive:
            \begin{enumerate}
                \item $s=s$. (Reflexivity)
                \item $\{s=t,u=v\}\tto (u=v)[s\mapsto t]$. (Replacement property)
            \end{enumerate}
            \subitem Derived:
            \begin{enumerate}
                \setcounter{enumi}{3}
                \item $s=t\otto t=x$. (Symmetry)
                \item $\{s=t,t=u\}\tto s=u$. (Transitivity)
                \item $s=t\tto u=u[s\mapsto t]$. (Leibniz's replacement)
            \end{enumerate}
        \end{itemize}
        \emph{}

        Although the above rules are actually default for (a \emph{variety} of) algebras such as Boolean algebra -- so they are usually not listed explicitly, it is best to specify the set of logical rules especially for an algebraic system that represents a logical system (e.g. Boolean algebra and Boolean lattice), even if these rules may be used tacitly.

        \subsubsection*{Boolean lattice}
        \label{subsubsec:Boolean lattice}

        A Boolean lattice $(X,1,0,\neg,\wedge,\vee,\le)$ is an extension of Boolean algebra $(X,1,0,\neg,\wedge,\vee)$ by adding the defining axiom/rule for the partial order ($\le$): $x\le y \otto x\wedge y=x$. The following can be derived:

        \begin{enumerate}
            \item $s\le s$.
            \item $\{s\le t,t\le s\}\tto s=t$.
            \item $\{s\le t,t\le u\}\tto s\le u$, $\{s=t,t\le u\}\tto s\le u$, $\{s\le t,t= u\}\tto s\le u$.
            \item $\{s=t,u\le v\}\tto (u\le v)[s\mapsto t]$.
            \item $s=t\tto u[s\mapsto t]\le u$, $s=t\tto u\le u[s\mapsto t]$.
        \end{enumerate}

        \subsubsection*{Definition of IRLA}

        \begin{definition}[Algebra IRLA]
            \label{def:IRLA}
            The algebraic system $IRLA=(X,1,0,\neg,\wedge,\vee,\oto,\to,\le)$ is the extension of Boolean lattice $BL=(X,1,0,\neg,\wedge,\vee,\le)$ by defining the new operations ($\oto$ and $\to$) with the axioms below:
            \begin{enumerate}
                \item $(x\oto y)\wedge z \le z[x\mapsto y]$. (Replacement property)
                \item $x=y\tto x\oto y=1$. (Tautological bi-implication)
                \item $x\to y = x\wedge y \oto x$. (Definition of implication)
            \end{enumerate}
        \end{definition}

        \subsection{Some theorems of IRLA}
        \label{subsec:Some theorems of IRLA}

        \begin{theorem}[Reflexivity of bi-implication in IRLA]
            \label{thrm:Reflexivity of bi-implication in IRLA}
            $x\oto x = 1$.
        \end{theorem}
        \begin{proof}\;\\
            \begin{tabular}{lll}
                1 & $x = x$ & Reflexivity of $=$\\
                2 & $x\oto x = 1$ & 1: Tautological bi-implication (Axiom 2)
            \end{tabular}\\
        \end{proof}

        \begin{theorem}[Symmetry of bi-implication in IRLA]
            \label{thrm:Symmetry of bi-implication in IRLA}
            $x\oto y = y\oto x$.
        \end{theorem}

        \begin{proof}\;\\
            \begin{tabular}{lll}
                1 & $x\oto x = 1$ & Reflexivity of $\oto$\\
                & $x\oto y$ & \\
                2 & \qquad$= (x\oto y)\wedge 1$ & Identity element\\
                3 & \qquad$= (x\oto y)\wedge (x\oto x)$ & 1: Leibniz's replacement of $=$\\
                4 & \qquad$\le (x\oto x)[x\mapsto y]$ & Replacement property of $\oto$ (Axiom 1)\\
                5 & \qquad$= x[x\mapsto y]\oto x$ & \emph{Selected replacement}\\
                6 & \qquad$= y\oto x$ & Result of the replacement\\
                7 & $x\oto y \le y\oto x$ & 2--6: Transitivity of $=$ and $\le$\\
                8 & $y\oto x \le x\oto y$ & 1--7: Similarly\\
                9 & $x\oto y = y\oto x$ & 7,8: Antisymmetry of $\le$
            \end{tabular}\\
        \end{proof}

        \begin{theorem}[Transitivity of bi-implication in IRLA]
            \label{thrm:Transitivity of bi-implication in IRLA}
            $$(x\oto y)\wedge(y\oto z) \le x\oto z.$$
        \end{theorem}

        \begin{proof}\;\\
            \begin{tabular}{lll}
                & $(x\oto y)\wedge(y\oto z)$ & \\
                1 & \qquad$= (y\oto z)\wedge(x\oto y) $ & Commutativity of $\wedge$\\
                2 & \qquad$\le (x\oto y)[y\mapsto z]$ & Replacement property of $\oto$ (Axiom 1)\\
                3 & \qquad$=(x\oto z)$ & Result of the replacement\\
                4 & $(x\oto y)\wedge(y\oto z) \le x\oto z$ & 1--3: Transitivity of $=$ and $\le$\\
            \end{tabular}\\
        \end{proof}

        \begin{theorem}[Leibniz's replacement in IRLA]
            \label{thrm:Leibniz's replacement in IRLA}
            $$x\oto y \le z\oto z[x\mapsto y].$$
        \end{theorem}

        \begin{proof}\;\\
            \begin{tabular}{lll}
                1 & $z\oto z = 1$ & Reflexivity\\
                & $x\oto y$ & \\
                2 & \qquad$= (x\oto y)\wedge 1$ & Identity element\\
                3 & \qquad$= (x\oto y)\wedge (z\oto z)$ & 1: Replacement property of $=$\\
                4 & \qquad$\le (z\oto z)[x\mapsto y]$ & Replacement property of $\oto$ (Axiom 1)\\
                5 & \qquad$= z\oto z[x\mapsto y]$ & \emph{Selected replacement}\\
                6 & $x\oto y \le z\oto z[x\mapsto y]$ & 2--5: Transitivity of $=$ and $\le$\\
            \end{tabular}\\
        \end{proof}

        \begin{theorem}[Tautological bi-implication and equality]
            \label{thrm:Tautological bi-implication and equality}
            $$x\oto y=1 \otto x=y.$$
        \end{theorem}

        \begin{proof}\;\\
            \begin{tabular}{lll}
                1&$x\oto y=1$&Premise\\
                &$x$&\\
                2&\qquad$=1\wedge x$&Identity element\\
                3&\qquad$=(x\oto y)\wedge x$&1: Replacement property of $=$\\
                4&\qquad$\le y$&Replacement property of $\oto$ (Axiom 1)\\
                5&$x\le y$&2--4: Transitivity of $=$ and $\le$\\
                6&$y\oto x=x\oto y$&Symmetry of $\oto$\\
                7&$y\oto x=1$&6,1:Transitivity of $=$\\
                8&$y\le x$ &7,1--5: Similarly\\
                9&$x=y$ & 5,8: Antisymmetry of $\le$
            \end{tabular}\\
            Thus we have $x\oto y=1\tto x=y$, and combining this with $x=y \tto x\oto y = 1$ (Axiom 2) we get $x\oto y = 1 \otto x=y$
        \end{proof}

        \begin{corollary}[Tautological implication and inequality]
            \label{corollary:Tautological implication and inequality}
            $$x\to y=1 \otto x\le y.$$
        \end{corollary}

        \begin{proof}
            By definition of $\to$ (Axiom 3), $x\to y=1\otto x\wedge y \oto x =1$, and $\le$ is defined by $x \le y \otto x\wedge y = x$. In terms of the above theorem, $x\wedge y \oto x =1 \otto x\wedge y = x$, so $x\to y=1 \otto x\le y$.
        \end{proof}

        In terms of Theorem \ref{thrm:Tautological bi-implication and equality} and Corollary \ref{corollary:Tautological implication and inequality}, an IRLA theorem of the form $s\oto t=1$ or $s\to t=1$ follows immediately from a corresponding theorem of Boolean lattice (BL) of the form $s = t$ or $s \le t$, so the latter is also a theorem in IRLA. Some examples of such theorems of IRLA and BL are listed in Table \ref{tab:Some theorems of IRLA and BL}.\\

        \begin{table}[H] 
            \centering
            \small
            \setlength{\tabcolsep}{8pt} 
            \renewcommand{\arraystretch}{1.28} 
            \begin{threeparttable}
                \caption{Some theorems of IRLA and BL}
                \label{tab:Some theorems of IRLA and BL}
                \begin{tabular}{r|rl}
                    \hline
                    Name &Form only in IRLA&Form in BL and IRLA\\
                    \hline
                    Idempotence of $\wedge$&$x\wedge x \oto x=1$&$x\wedge x = x$\\
                    Idempotence of  $\vee$&$x\vee x \oto x=1$&$x\vee x = x$\\
                    Commutativity of $\wedge$&$x\wedge y \oto y\wedge x=1$&$x\wedge y = y\wedge x$\\
                    Commutativity of $\vee$&$x\vee y \oto y\vee x=1$&$x\vee y = y\vee x$\\
                    De Morgan's law 1&$\neg(x\wedge y)\oto\neg x\vee \neg y=1$&$\neg(x\wedge y)=\neg x\vee \neg y$\\
                    De Morgan's law 2&$\neg(x\vee y)\oto\neg x\wedge \neg y=1$&$\neg(x\vee y)=\neg x\wedge \neg y$\\
                    Reflexivity of implication&$x\to x=1$&$x\le x$\\
                    Conjunction elimination&$x\wedge y\to x=1$&$x\wedge y\le x$\\
                    Disjunction introduction&$x\to x\vee y=1$&$x\le x\vee y$\\
                    Boundedness of $\to$ &$0\to x=1,\;x\to 1=1$&$0\le x,\;x\le 1$\\
                    \hline
                \end{tabular}
            \end{threeparttable}
        \end{table}

        \begin{theorem}[Bi-implication introduction / Antisymmetry of implication in IRLA]
            \label{thrm:Bi-implication introduction in IRLA}
            $$(x\to y)\wedge (y\to x)\le x\oto y.$$
        \end{theorem}

        \begin{proof}\;\\
            \begin{tabular}{lll}
                &$(x\to y)\wedge (y\to x)$&\\
                1&\qquad$= (x\oto x\wedge y)\wedge (x\wedge y\oto y)$&Definition of $\to$\\
                2&\qquad$\le x\oto y$&Transitivity of $\oto$\\
                3&$(x\to y)\wedge (y\to x)\le x\oto y$&1--2: Transitivity of $=$ and $\le$
            \end{tabular}\\
        \end{proof}

        \begin{theorem}[Bi-implication elimination in IRLA]
            \label{thrm:Bi-implication elimination in IRLA}
            $$x\oto y\le x\to y,\:x\oto y\le y\to x.$$
        \end{theorem}

        \begin{proof}\;\\
            \begin{tabular}{lll}
                &$x\oto y$&\\
                1&\qquad$= (x\oto y)\wedge 1$&Identity element\\
                2&\qquad$= (x\oto y)\wedge(x\to x)$&Reflexivity of $\to$ ($x\oto x=1$)\\
                3&\qquad$\le (x\to x)[x\mapsto y]$& Replacement property of $\oto$ (Axiom 1)\\
                4&\qquad$= x\to x[x\mapsto y]$& Selected replacement\\
                5&\qquad$= x\to y$& Result of the replacement\\
                6&$x\oto y \le x\to y$&1--5: Transitivity of $=$ and $\le$
            \end{tabular}\\
            The second is obtained by symmetry of $\oto$.
        \end{proof}

        \begin{theorem}[Bi-implication by implication in IRLA]
            \label{thrm:Bi-implication by implication in IRLA}
            $$x\oto y=(x\to y)\wedge(y\to x).$$
        \end{theorem}

        \begin{proof}\;\\
            \begin{tabular}{lll}
                1&$x\oto y \le x\to y$&Bi-implication elimination\\
                2&$x\oto y \le y\to x$&Bi-implication elimination\\
                3&$x\oto y \le (x\to y)\wedge(y\to x)$&1,2: One-to-two of $\le$\\
                4&$(x\to y)\wedge(y\to x) \le x\oto y$&Bi-implication introduction\\
                5&$x\oto y=(x\to y)\wedge(y\to x)$&3,4: Antisymmetry of $\le$
            \end{tabular}\\
        \end{proof}

        \begin{theorem}[Transitivity of implication in IRLA]
            \label{thrm:Transitivity of implication in IRLA}
            $$(x\to y)\wedge (y\to z)\le x\to z.$$
        \end{theorem}

        \begin{proof}\;\\
            \begin{tabular}{lll}
                &$(x\to y)\wedge (y\to z)$&\\
                1&\qquad$= (x\wedge y\oto x)\wedge (y\wedge z\oto y)$&Definition of $\to$\\
                2&\qquad$= (x\wedge y\oto x)\wedge(x\wedge y\oto x)\wedge(y\wedge z\oto y)$&Idempotence\\
                3&\qquad$\le (x\wedge y\oto x)\wedge(x\wedge (y\wedge z)\oto x)$&Replacement $[y\mapsto y\wedge z]$\\
                4&\qquad$= (x\wedge y\oto x)\wedge((x\wedge y)\wedge z\oto x)$&Associativity\\
                5&\qquad$\le x\wedge z\oto x$&Replacement $[x\wedge y\mapsto x]$\\
                6&\qquad$=  x\to z$&Definition of $\to$\\
                7&$(x\to y)\wedge (y\to z)\le x\to z$&1--6: Transitivity of $=$ and $\le$
            \end{tabular}\\
        \end{proof}

        \begin{theorem}[One-to-two and Two-to-one in IRLA]
            \label{thrm:One-to-two and Two-to-one in IRLA}
            \begin{equation*}
                \begin{aligned}
                    (z\to x)\wedge(z\to y) &= z\to x\wedge y.\\
                    (x\to z)\wedge(y\to z) &= x\vee y\to z.
                \end{aligned}
            \end{equation*}
        \end{theorem}

        \begin{proof}\;\\
            \begin{tabular}{lll}
                &$(z\to x)\wedge(z\to y)$&\\
                1&\qquad$=(x\wedge z\oto z)\wedge(y\wedge z\to z)$&Definition of $\to$\\
                2&\qquad$\le x\wedge y\wedge z\oto z$&Replacement $[z\mapsto y\wedge z]$\\
                3&\qquad$= z\to x\wedge y$&Definition of $\to$\\
                4&$(z\to x)\wedge(z\to y)\le z\to x\wedge y$&1--3: Transitivity of $=$ and $\le$\\
                &$z\to x\wedge y$&\\
                5&\qquad$=(z\to x\wedge y)\wedge(x\wedge y\to x)$&Conjunction elimination ($x\wedge y\to x=1$)\\
                6&\qquad$\le z\to x$&Transitivity of $\to$\\
                7&$z\to x\wedge y\le z\to x$&5--6: Transitivity of $=$ and $\le$\\
                8&$z\to x\wedge y\le z\to y$&5--7: Similarly\\
                9&$z\to x\wedge y\le (z\to x)\wedge(z\to y)$&7,8: One-to-two of $\le$\\
                10&$(z\to x)\wedge(z\to y) = z\to x\wedge y$&4,9: Antisymmetry of $\le$
            \end{tabular}\\
            The second is proved dually.
        \end{proof}

        \begin{theorem}[Modus ponens in IRLA]
            \label{thrm:Modus ponens in IRLA}
            $$(x\to y)\wedge x\le y.$$
        \end{theorem}

        \begin{proof}\;\\
            \begin{tabular}{lll}
                &$(x\to y)\wedge x$&\\
                1&\qquad$=(x\oto x\wedge y)\wedge x$& Definition of $\to$\\
                2&\qquad$\le x\wedge y$& Replacement property of $\oto$ ($[x\mapsto x\wedge y]$)\\
                3&\qquad$\le y$& Conjunction elimination\\
                4&$(x\to y)\wedge x\le y$& 1--3: Transitivity of $=$ and $\le$
            \end{tabular}\\
        \end{proof}

        \begin{corollary}[Upper bound of implication]
            \label{corollary:Upper bound of implication}
            $$x\to y\le \neg x\vee y.$$
        \end{corollary}

        \begin{proof}
            In BL, $x\wedge y\le z\otto x\le \neg y\vee z$, so modus ponens $(x\to y)\wedge x\le y$ is equivalent to the upper bound of implication $x\to y\le \neg x\vee y$ in IRLA, hence the latter follows from the above theorem.
        \end{proof}

        \begin{theorem}[Tautological implication and disjunction]
            \label{thrm:Tautological implication and disjunction}
            $$x\to y=1\otto \neg x\vee y=1.$$
        \end{theorem}

        \begin{proof}
            It follows from $x\to y=1 \otto x\le y$ of IRLA and $x\le y\otto \neg x\vee y=1$ of BL.
        \end{proof}

        \begin{theorem}[Contradictory implication]
            \label{thrm:Contradictory implication}
            $$1\to 0=0.$$
        \end{theorem}

        \begin{proof}
            It follows from boundedness $0\le 1\to 0$ and $1\to 0\le \neg 1 \vee 0 = 0\vee 0 =0 $ by Corollary \ref {corollary:Upper bound of implication}.
        \end{proof}

        \subsection{Soundness and completeness}
        \label{subsec:Soundness and completeness}

        By definition (cf. Subsection \ref{subsec:Semantics of IRL}), the T-set of IRL formula schema $F(\phi,\psi,...)$ is the IRLA term $t(x,y,...)$; and $\vDash\phi$ is the same as $x=1$,  $\vDash\phi\oto\psi$ is the same as $x=y$ or $x\oto y=1$, and $\vDash\phi\to\psi$ is the same as $x\le y$ or $x\to y=1$.

        \subsubsection*{Soundness of IRL}
        \label{subsubsec:Soundness of IRL}

        \begin{metatheorem}[Axioms are tautologies]
            \label{metathrm:Axioms are tautologies}
            All axioms of IRL are tautologies.
        \end{metatheorem}
        \begin{proof}[Sketch of proof]
            For each axiom schema $F(\phi,\psi,...)$ of IRL, verify $t(x,y,...)=1$ in IRLA. For example, $\phi\oto\phi$ is a tautology since $x\oto x=1$ in IRLA.
        \end{proof}

        \begin{metatheorem}[Primitive rules preserve tautologousness]
            \label{metathrm:Primitive rules preserve tautologousness}
            All the three primitive rules (except for the meta rule of conditional proof) of IRL, i.e. Conjunction introduction $\{\phi,\psi\}\tto\phi\wedge\psi$, Conjunction elimination $\phi\wedge\psi\tto\phi$ and Replacement property $\{\phi\oto\psi,\chi\}\tto\chi[\phi\mapsto\psi]$, are tautologousness-preserving, that is, if their premises are tautologies, so are their conclusions.
        \end{metatheorem}

        \begin{proof}\;

            Conjunction introduction $\{\phi,\psi\}\tto\phi\wedge\psi$:

            \begin{tabular}{lll}
                1 & $\vDash\phi$ & Known\\
                2 & $\vDash\psi$ & Known\\
                3 & $x=1$ & 1: Definition\\
                4 & $y=1$ & 2: Definition\\
                & $x\wedge y$ &\\
                5 & \qquad $= 1\wedge 1$ & 3,4: Replacement\\
                6 & \qquad $= 1$ & Identity element\\
                7 & $x\wedge y = 1$ & 5,6: Transitivity\\
                8 & $\vDash\phi\wedge\psi$ & 7: Definition\\
            \end{tabular}

            Conjunction elimination $\phi\wedge\psi\tto\phi$:

            \begin{tabular}{lll}
                1 & $\vDash\phi\wedge\psi$ & Known\\
                2 & $x\wedge y = 1$ & 1: Definition\\
                & $1$ & \\
                3 & \qquad $=x\wedge y$ & 2: Symmetry\\
                4 & \qquad $\le x$ & Conjunction elimination\\
                5 & $1\le x$ & 3,4: Transitivity of $=$ and $\le$\\
                6 & $x \le 1$ & Boundedness\\
                7 & $x = 1$ & 5,6: Antisymmetry of $\le$\\
                8 & $\vDash\phi$ & 7: Definition\\
            \end{tabular}

            Replacement property $\{\phi\oto\psi,\chi\}\tto\chi[\phi\mapsto\psi]$:

            \begin{tabular}{lll}
                1 & $\vDash\phi\oto\psi$ & Known\\
                2 & $\vDash\chi$ & Known\\
                3 & $x\oto y =1 $ & 1: Definition\\
                4 & $z=1$ & 2: Definition\\
                & $1$ &\\
                5 & \qquad $=1\wedge 1$ & Identity element\\
                6 & \qquad $=(x\oto y)\wedge z$ & 3,4: Replacement\\
                7 & \qquad $\le z[x\mapsto y]$ & Replacement property of $\oto$\\
                8 & \qquad $\le 1$ & Boundedness\\
                9 & $z[x\mapsto y] = 1$ & 5--7--8: Transitivity of $=$ and $\le$, and Antisymmetry of $\le$\\
                10 & $\vDash\chi[\phi\mapsto\psi]$ & 9: Definition\\
            \end{tabular}
        \end{proof}

        \begin{metatheorem}[CP rule preserves tautologousness]
            \label{metathrm:CP rule is tautologousness-preserving}
            Suppose a conditional proof $\phi\vdash\psi$ of IRL satisfies that if a member of its proof sequence is a theorem then it is a tautology, then $\phi\to\psi$ is a tautology. That is, the rule of conditional proof (the CP rule) is tautologousness-preserving.
        \end{metatheorem}

        \begin{proof}
            We prove by induction that, for any integer $n\ge 1$, for a conditional proof $\phi\vdash\psi$ of length $n$ satisfying that \emph{if a member of its proof sequence is a theorem then it is a tautology} (let it be called ``theorem-is-tautology'' condition below), then $\phi\to\psi$ is a tautology.

            Only the three primitive rules need considered in the meta proof: Conjunction introduction $\{\phi,\psi\}\tto\phi\wedge\psi$, Conjunction elimination $\phi\wedge\psi\tto\phi$, and Replacement property $\{\phi\oto\psi,\chi\}\tto\chi[\phi\mapsto\psi]$.

            \vspace{\baselineskip}
            Let $\alpha$, $\beta$ and $\gamma$ be IRL formulas that exist in the cases discussed below, and their IRLA counterparts be $p$, $q$ and $r$ respectively.

            \begin{enumerate}

                \item For $n=1$: The proof sequence of $\phi\vdash\psi$ of length $1$ is just $(\psi)$. There are following two cases for $\psi$.

                \begin{itemize}
                    \item Case 1: $\psi$ is a theorem. In this case, $\psi$ is a tautology by the ``theorem-is-tautology'' condition, so $\phi\to\psi$ must be a tautology as shown below:

                    \begin{tabular}{lll}
                        1 & $\vDash\psi$ & Known\\
                        2 & $y=1$ & 1: Definition\\
                        3 & $x\le 1$ & Boundedness\\
                        4 & $x\le y$ & 2,3: Replacement\\
                        5 & $\vDash\phi\to\psi$ & 4: Definition\\
                    \end{tabular}

                    \item Case 2: $\psi$ is obtained by a one-premise rule. In this case, $\psi$ can only be obtained by the rule $\psi\wedge\alpha\tto\psi$ where $\psi\wedge\alpha=\phi$, so $\phi\to\psi$ must be a tautology as shown below:

                    \begin{tabular}{lll}
                        1 & $\psi\wedge\alpha=\phi$ & Known\\
                        2 & $y\wedge p = x$ & 1: Definition\\
                        3 & $y\wedge p\le y$ & Conjunction elimination\\
                        4 & $x\le y$ & 2,3: Replacement\\
                        5 & $\vDash\phi\to\psi$ & 4: Definition\\
                    \end{tabular}

                \end{itemize}
                Thus, if $\phi\vdash\psi$ of length $1$ satisfies the ``theorem-is-tautology'' condition, then $\phi\to\psi$ is a tautology.

                \item For any $n\ge 1$: Suppose that, for any $k=1,...,n$, if $\phi\vdash\psi$ of length $k$ satisfies the ``theorem-is-tautology'' condition, then $\phi\to\psi$ is a tautology; it need to show that if $\phi\vdash\psi$ of length $n+1$ satisfies the ``theorem-is-tautology'' condition, then $\phi\to\psi$ is a tautology.

                Let $(\psi_1,...,\psi_n,\psi)$ be a proof sequence of the length-$(n+1)$ conditional proof $\phi\vdash\psi$. Then, for any $\chi \in \{\phi,\psi_1,...,\psi_n\}$, $\phi\to\chi$ must be a tautology as shown below:

                \begin{tabular}{lll}
                    1 & $\chi \in \{\phi,\psi_1,...,\psi_n\}$ & Known\\
                    2 & $\chi=\phi$ or $\chi \in \{\psi_1,...,\psi_n\}$ & 1: Equivalent to it\\
                    3 & \qquad $\chi=\phi$ & Assumption\\
                    4 & \qquad $z=x$ & 3: Definition\\
                    5 & \qquad $x\le x$ & Reflexivity of $\le$\\
                    6 & \qquad $x\le z$ & 4,5: Replacement\\
                    7 & \qquad $\vDash\phi\to\chi$ & 6: Definition\\
                    & \qquad  & \\
                    8 & \qquad $\chi\in \{\psi_1,...,\psi_n\}$ & Assumption\\
                    9 & \qquad $\vDash\phi\to\chi$ & 8: Induction hypothesis (as $\vDash\phi\to\psi_k$ for all $k=1,...,n$)\\
                    10 & $\vDash\phi\to\chi$ & 2,3--7,8--9: Proof by exhaustion\\
                \end{tabular}

                This result is used in the following four cases for $\psi$.

                \begin{itemize}
                    \item Case 1: $\psi$ is a theorem. In this case, we proceed as for $n=1$.

                    \item Case 2: $\psi$ is obtained by the rule $\{\alpha,\beta\}\tto \alpha\wedge \beta$ where $\alpha,\beta \in \{\phi,\psi_1,...,\psi_n\}$ and $\alpha\wedge \beta=\psi$. In this case, $\phi\to\alpha$ and $\phi\to\beta$ are tautologies, so $\phi\to\psi$ must be a tautology as shown below:

                    \begin{tabular}{lll}
                        1 & $\vDash\phi\to\alpha$ & Known\\
                        2 & $\vDash\phi\to\beta$ & Known\\
                        3 & $\alpha\wedge\beta=\psi$ & Known\\
                        4 & $x\le p$ & 1: Definition\\
                        5 & $x\le q$ & 2: Definition\\
                        6 & $p\wedge q=y$ & 3: Definition\\
                        7 & $x\le p\wedge q$ & 4,5: One-to-two\\
                        8 & $x\le y$ & 6,7: Replacement\\
                        9 & $\vDash\phi\to\psi$ & 8: Definition\\
                    \end{tabular}

                    \item Case 3: $\psi$ is obtained by the rule $\psi\wedge \alpha\tto\psi$ where $\psi\wedge \alpha \in \{\phi,\psi_1,...,\psi_n\}$. In this case, $\phi\to\psi\wedge \alpha$ is a tautology, so $\phi\to\psi$ must be a tautology as shown below:

                    \begin{tabular}{lll}
                        1 & $\vDash\phi\to\psi\wedge \alpha$ & Known\\
                        2 & $x\le y\wedge p$ & 1: Definition\\
                        3 & $y\wedge p\le y$ & Conjunction elimination\\
                        4 & $x\le y$ & 2,3: Transitivity of $\le$\\
                        5 & $\vDash\phi\to\psi$ & 4: Definition\\
                    \end{tabular}

                    \item Case 4: $\psi$ is obtained by the rule $\{\alpha\oto\beta,\gamma\}\tto\psi$ where $\alpha\oto\beta,\gamma \in \{\phi,\psi_1,...,\psi_n\}$ and $\gamma[\alpha\mapsto\beta]=\psi$. In this case, $\phi\to (\alpha\oto\beta)$ and $\phi\to\gamma$ are tautologies, so $\phi\to\psi$ must be a tautology as shown below:

                    \begin{tabular}{lll}
                        1 & $\vDash\phi\to(\alpha\oto\beta)$ & Known\\
                        2 & $\vDash\phi\to\gamma$ & Known\\
                        3 & $\gamma[\alpha\mapsto\beta]=\psi$ & Known\\
                        4 & $x\le p\oto q$ & 1: Definition\\
                        5 & $x\le r$ & 2: Definition\\
                        6 & $r[p\mapsto q]=y$ & 3: Definition\\
                        7 & $x\le (p\oto q)\wedge r$ & 4,5: One-to-two\\
                        8 & $(p\oto q)\wedge r\le r[p\mapsto q]$ & Replacement property of $\oto$\\
                        9 & $x\le r[p\mapsto q]$ & 7,8: Transitivity of $\le$\\
                        10 & $x\le y$ & 6,9: Replacement\\
                        11 & $\vDash\phi\to\psi$ & 10: Definition\\
                    \end{tabular}

                \end{itemize}
                Thus, for all $n\ge 1$, suppose that for all $k=1,...,n$, any conditional proof $\phi\vdash\psi$ of length $k$ satisfying the ``theorem-is-tautology'' condition results in tautological $\phi\to\psi$, then any length-$(n+1)$ conditional proof $\phi\vdash\psi$ satisfying the ``theorem-is-tautology'' condition leads to tautological $\phi\to\psi$ as well.
            \end{enumerate}
            Therefore, from 1 and 2 by induction, we have that for any $n\ge 1$, if a conditional proof $\phi\vdash\psi$ of length $n$ satisfies the ``theorem-is-tautology'' condition, then $\phi\to\psi$ is a tautology.
        \end{proof}

        From Meta Theorem \ref{metathrm:Axioms are tautologies}, Meta Theorem \ref{metathrm:Primitive rules preserve tautologousness} and Meta Theorem \ref{metathrm:CP rule is tautologousness-preserving}, it can be concluded that IRL is sound with respect to its semantics specified by IRLA.

        \subsubsection*{Completeness of IRL}
        \label{subsubsec:Completeness of IRL}

        In the following two meta theorems, IRLA axiom $x=y\tto x\oto y=1$ (one of the two defining axioms for $\oto$) is considered to be a primitive rule of IRLA, and the universal rule $s=s$ (one of the two primitive universal rules in IRLA) is considered to be an axiom of IRLA.

        \begin{metatheorem}
            \label{metathrm:IRLA axioms have corresponding IRL theorems}
            Any axiom in IRLA has a corresponding theorem in IRL.
        \end{metatheorem}

        \begin{proof}[Sketch of proof]\;
            \begin{enumerate}
                \item All axioms of IRLA that are inherited from Boolean lattice have their counterparts in IRL. For example, IRLA axiom $x\wedge y=y\wedge x$ corresponds to IRL axiom $\phi\wedge\phi\oto\psi\wedge\phi$.
                \item The IRLA axiom $(x\oto y)\wedge z \le z[x\mapsto y]$ (one of the two conditions for defining $\oto$ in IRLA) corresponds to IRL theorem $(\phi\oto\psi)\wedge\chi\to\chi[\phi\mapsto\psi]$ (from the primitive rule $\{\phi\oto\psi,\chi\}\tto\chi[\phi\mapsto\psi]$ of IRL).
                \item The IRLA axiom $x\to y=x\wedge y\oto x$ (the defining axiom for $\to$) corresponds to IRL axiom $(\phi\to\psi)\oto(\phi\wedge\psi\oto\phi)$.
                \item The universal one $s=s$ corresponds to IRL axiom $\phi\oto\phi$.
            \end{enumerate}
        \end{proof}

        \begin{metatheorem}
            \label{metathrm:IRLA primitive rules preserve have-IRL-theorem property}
            For all the two primitive rules of IRLA, namely $x=y\tto x\oto y=1$ and $\{s=t,u=v\}\tto (u=v)[s\mapsto t]$, if their premises have corresponding theorems in IRL, so are their conclusions.
        \end{metatheorem}

        \begin{proof}\;
            \begin{enumerate}
                \item $x=y\tto x\oto y=1$: In IRL, $\phi\oto\psi$ is a theorem if and only if $(\phi\oto\psi)\oto \top$ is a theorem by the rule of theorem representation. Thus, if the premise $x=y$ has a corresponding theorem $\phi\oto\psi$ in IRL, then the conclusion $x\oto y=1$ has the counterpart theorem $(\phi\oto\psi)\oto \top$ in IRL.

                \item $\{s=t,u=v\}\tto (u=v)[s\mapsto t]$: For the premises $s=t$ and $u=v$, if their counterparts $ E \oto F $ and $G\oto H $ are theorems of IRL, then $(G\oto H )[ E \mapsto F ]$, the counterpart in IRL for the conclusion $(u=v)[s\mapsto t]$, must be a theorem in IRL as shown below:

                \begin{tabular}{lll}
                    1&$\vdash E \oto F $&Known\\
                    2&$\vdash G\oto H $&Known\\
                    3&\qquad$ E \oto F $&1 (assumed theorem)\\
                    4&\qquad$G\oto H $&2 (assumed theorem)\\
                    5&\qquad$(G\oto H )[ E \mapsto F ]$&3,4: Replacement (IRL rule)\\
                    6&$\vdash(G\oto H )[ E \mapsto F ]$&3--5: Unconditional proof
                \end{tabular}
            \end{enumerate}
        \end{proof}

        From Meta Theorem \ref{metathrm:IRLA axioms have corresponding IRL theorems} and Meta Theorem \ref{metathrm:IRLA primitive rules preserve have-IRL-theorem property}, it can be concluded that IRL is complete with respect to its semantics specified by IRLA.

        \subsection{Isomorphism between IRL and IRLA}
        \label{subsec:Isomorphism between IRL and IRLA}

        \begin{definition}[Isomorphism between two systems]
            \label{def:Isomorphism between two systems}
            Two proof systems $A$ and $B$ (possibly for different formal languages) are isomorphic (to each other), denoted as $A\cong B$, if the following conditions are met:
            \begin{enumerate}
                \item Expressions of the two systems can be translated to each other. Specifically, there is a bijection between expressions of $A$ and expressions of $B$, where ``expressions'' include not only formulas of object languages but also expressions of meta languages such as that for rules and meta rules.
                \item Axioms and primitive rules of $A$ (after translated to expressions of $B$) can be derived in $B$, and vice versa.
            \end{enumerate}
        \end{definition}

        Isomorphism between two systems can be equivalently expressed as that they are interpretable in each other in terms of the concept of interpretability \citep{tarski}.

        \vspace{\baselineskip}
        Two systems $A$ and $B$ have the following properties related to isomorphism:
        \begin{enumerate}
            \item If $A$ and $B$ are isomorphic, then they are equiconsistent and equidecidable.
            \item Suppose $B$ is the semantic system of $A$. Then $A$ and $B$ are isomorphic if and only if $A$ is sound and complete with regard to $B$.
        \end{enumerate}

        Since IRL is sound and complete with respect to IRLA, there is a natural isomorphism between them as shown in Table \ref{tab:Isomorphism between IRL and IRLA}.

        \begin{table}[H] 
            \centering
            \small
            \setlength{\tabcolsep}{2pt} 
            \renewcommand{\arraystretch}{1.28} 
            \begin{threeparttable}
                \caption{Isomorphism between IRL and IRLA}
                \label{tab:Isomorphism between IRL and IRLA}
                \begin{tabular}{rcl}
                    \hline
                    \emph{IRL} &$\cong$& \emph{IRLA}\\
                    \hline
                    Primitive propositions\quad $A,B,...$ && $a,b,...$\quad Constants\\
                    Formulas\quad $\phi,\psi,...$\quad && $x,y,...$\quad Variables\\
                    Formula schema\quad $F(\phi,\psi,...)$ && $t(x,y,...)$\quad Term \\
                    Operations\quad $\top,\bot,\neg,\wedge,\vee,\oto,\to$ && $1,0,\neg,\wedge,\vee,\oto,\to$\\
                    \hline
                    Theorem\quad $\vdash\phi$ && $x=1$\quad Tautology\\
                    Theorem of negation\quad $\vdash\neg\phi$ && $x=0$\quad Contradiction\\
                    Independent formula\quad $\nvdash\phi$ and $\nvdash\neg\phi$ && $x\ne 1$ and $x\ne 0$\quad Contingency\\
                    \hline
                    Bi-implication theorem\quad $\vdash\phi\oto\psi$ or $\phi\otto\psi$ \tnote{a} && $x\oto y=1$ or $x=y$ \tnote{a}\quad Tautological bi-implication\\
                    Implication theorem\quad $\vdash\phi\to\psi$ or $\phi\vdash\psi$ or $\phi\tto\psi$ \tnote{b} && $x\to y=1$ or $x\le y$  \tnote{b}\quad Tautological implication\\
                    \hline
                    Meta proof (e.g.) \quad $(\phi\vdash\psi)\vdash(\vdash\phi\to\psi)$ \tnote{c} && $x\le y\vdash x\to y=1$  \tnote{c} \quad Conditional proof in IRLA (e.g.)\\
                    Meta Rule (e.g.) \quad $(\phi\vdash\psi)\tto(\vdash\phi\to\psi)$ \tnote{c} && $x\le y\tto x\to y=1$ \tnote{c} \quad Rule in IRLA (e.g.)\\
                    \hline
                \end{tabular}
                \begin{tablenotes}
                    \item[a] In IRL, $\vdash\phi\oto\psi$ and $\phi\otto\psi$ are equivalent; in IRLA, $x\oto y=1$ and $x=y$ are equivalent; so any of the two in IRL can be translated to any of the two in IRLA, and vice versa.
                    \item[b] In IRL, $\vdash\phi\to\psi$, $\phi\vdash\psi$ and  $\phi\tto\psi$ are equivalent; in IRLA, $x\to y=1$ and $x\le y$ are equivalent; so any of the three in IRL can be translated to any of the two in IRLA, and vice versa.
                    \item[c] For simplicity, $\vdash$ are used for both proof and meta proof, and $\tto$ are used for both rule and meta rule, and they can be used in different systems, if they can be differentiated with context.
                \end{tablenotes}
            \end{threeparttable}
        \end{table}

        As IRL is sound and complete, for a formula $\phi$ of IRL, ``$\phi$ is valid or is a tautology'' is equivalent to ``$\phi$ is a theorem'', ``$\phi$ is unsatisfiable or is a contradiction'' is equivalent to ``$\neg\phi$ is a theorem''; so ``$\phi$ is a contingency'' is equivalent to ``both $\phi$ and $\neg\phi$ are not theorems'' i.e. ``$\phi$ is independent of IRL''.

        \subsection{Truth values of implication}
        \label{subsec:Truth values of implication}
        The truth values of an implication $\phi\to\psi$ is represented by its T-set $V_{\phi\to\psi}$ i.e. $x\to y$ in IRLA.

        Define $x<y$ as $x\le y$ and $x\ne y$. Thus $0<x<1$ means that $x$ is (represents) a contingency.

        By Corollary \ref{corollary:Upper bound of implication}, in any cases we have the bounds of implication as in
        $$0\le x\to y\le \neg x\vee y, $$
        and both the lower bound $0$ and the upper bound $\neg x \vee y$ of $x\to y$ are \emph{sharp}, i.e. \emph{achievable}, since we have $1\to 0=0$ and $x\to x=\neg x\vee x$.

        Combined with Theorem \ref{thrm:Tautological implication and disjunction}, there are three specific cases for the truth values of an implication:

        \begin{center}
            \begin{tabular}{ccc}
                $x\to y=1$ & if and only if & $\neg x\vee y=1$\\
                $x\to y=0$ & if & $\neg x\vee y=0$ (i.e. $x=1$ and $y=0$)\\
                $0\le x\to y\le \neg x\vee y<1$ & if & $0<\neg x\vee y<1$\\
            \end{tabular}
        \end{center}

        Thus, generally $V_{\phi\to\psi}$ is a subset of $V_{\neg\phi\vee\psi}$ (including of course $\emptyset$ and $V_{\neg\phi\vee\psi}$). Especially in the case where $\neg\phi\vee\psi$ is contingent, unlike classical logic that always adopts the equality $V_{\phi\to\psi}=V_{\neg\phi\vee\psi}$, $V_{\phi\to\psi}$ can be any subset of $V_{\neg\phi\vee\psi}$, i.e. $V_{\phi\to\psi}\in 2^{V_{\neg\phi\vee\psi}}$, so there are $2^{|V_{\neg\phi\vee\psi}|}$ subsets, of which only one is ``correct'' for $\phi\to\psi$. The correct $V_{\phi\to\psi}$ is determined by the original meaning of implication -- a ``certain mechanism'' (cf. Subsection \ref{subsec:OriginalMeaningOfImplication}) corresponding to the implication $\phi\to\psi$ considered.

        \begin{example}
            \label{exam:PQ}
            Let $A$ and $B$ are primitive propositions having 4 possible truth-value combinations, i.e. the set of interpretations $U=\{\mathrm{TT,TF,FT,FF}\}$. Then $\emptyset\subset V_{\neg A\vee B}=\{\mathrm{TT,FT,FF}\}\subset U$, so $\neg A\vee B$ is contingent, therefore $V_{A\to B}\in 2^{ V_{\neg A\vee B}}=2^{\{\mathrm{TT,FT,FF}\}}$, thus there are $2^{|\{\mathrm{TT,FT,FF}\}|}=2^3=8$ possibilities and only one of them is correct for $V_{A\to B}$.
        \end{example}

        The following is an example that shows which should be the ``correct one'' for the T-set of implication $\phi\to\psi$ as a subset of the T-set of $\neg\phi\vee\psi$.

        \begin{example}
            \label{exam:OPQ}
            Let $O=$``in the open air'', $P=$``it rains'', and $Q=$``the ground gets wet''. Consider the truth values of $P\to Q=$``if it rains, then the ground gets wet''. Suppose: (a) In the open air, whenever it rains the ground must get wet; (b) ``it rains'' is not the only cause of ``the ground gets wet'', no matter it is in the open air or not, e.g., if someone sprinkles water coincidentally on the ground (or the floor), the ground (or the floor) gets wet as well. Then, we obtain the truth values of $P\to Q$ in all 7 possible cases (there is one impossible case that must be excluded: ``in the open air, it rains and the ground does not get wet''), as shown in Table \ref{tab:Truth table for Example OPQ}.

            \begin{table}[H] 
                \centering
                \small
                \setlength{\tabcolsep}{17pt} 
                \renewcommand{\arraystretch}{1.28} 
                \begin{threeparttable}
                    \caption{Truth table for Example \ref{exam:OPQ}}
                    \label{tab:Truth table for Example OPQ}
                    \begin{tabular}{c|c|cc|c|c}
                        \hline
                        & $O$ & $P$ & $Q$ & $\neg P\vee Q$ \tnote{a} & $P\to Q$ \tnote{b} \\
                        \hline
                        1 & T & T & T & \textbf{T} & \textbf{T} \\
                        2 & T & F & T & \textbf{T} & \textbf{T} \\
                        3 & T & F & F & \textbf{T} & \textbf{T} \\
                        4 & F & T & T & \textbf{T} & F \\
                        5 & F & T & F & F & F \\
                        6 & F & F & T & \textbf{T} & F \\
                        7 & F & F & F & \textbf{T} & F \\
                        \hline
                    \end{tabular}
                    \begin{tablenotes}
                        \item[a] T-set of $\neg P\vee Q$ is the set $\{1,2,3,4,6,7\}$ of rows, indicating that it is contingent.
                        \item[b] T-set of $P\to Q$, determined by T-set of $O$ rather than that of $P$ and $Q$, is the set $\{1,2,3\}$ of rows, a proper subset of that of $\neg P\vee Q$.
                    \end{tablenotes}
                \end{threeparttable}
            \end{table}

        \end{example}

        In this example, it is the \emph{certain mechanism} of ``in the open air'' that corresponds to the inevitable implication relation of ``if it rains, then the ground gets wet'', and the \emph{scope} of the mechanism is the T-set of the implication. This is a typical example showing that it is not uncommon that $\neg\phi\vee\psi\tto \phi\to\psi$ is invalid. The example can explain why the classical logic, with the common ``implication'' unreasonably defined by $\phi\to\psi \otto\neg\phi\vee\psi$, is doomed to fail frequently in common situations in real world.


        \section{Relationship between IRL and CL}
        \label{sec:Relationship between IRL and CL}

        \subsection{Definitions and notations}
        \label{subsec:Definitions and notations 2}

        In the following, view a proof system as the set of its axioms \emph{and} primitive rules, and use ``statement'' to mean formula or rule expression of a system; treat an object and the singleton containing the object the same thing.

        \begin{definition}
            \label{def:notation for both theorem and derived rule}
            Let $A$ be a system. Let $X$ be a set of statements of the system. Denote with $\vdash_A X$ that each member $x$ of $X$ is derivable in $A$ (i.e. $x$ is a theorem if it is a formula or $x$ is a derivable rule if it is a rule expression).
        \end{definition}

        \begin{definition}[Supper-system and sub-system]
            \label{def:Supper-system and sub-system}
            Let $A$ and $B$ be systems for a same language.

            $A$ is a supper-system of $B$, or $B$  is a sub-system of $A$, if $\vdash_A B$.

            $A$ and $B$ are (deductively) equivalent if $\vdash_A B$ and $\vdash_B A$, denoted as $A=B$.
        \end{definition}

        Note: For systems $A$ and $B$ for a \emph{same} language, $A\cong B$ and $A=B$ are equivalent.

        \begin{definition}[System extension]
            \label{def:System extension}
            Let $A$ be a system. Let $X$ be a set of statements.

            An extension of $A$ with $X$, denoted as $A+X$, is the system obtained from $A$ by adding each member of $X$ as an axiom or a primitive rule.

        \end{definition}

        \begin{definition}[Provable consequence and equivalence]
            \label{def:Provable consequence and equivalence}
            Let $A$ be a system. Let $X$ and $Y$ be sets of statements.

            Then $Y$ is called a \emph{provable consequence} of $X$ in $A$, if $\vdash_{A+X} Y$.

            $X$ and $Y$ are \emph{provably equivalent} in $A$, if $\vdash_{A+X} Y$ and $\vdash_{A+Y} X$.

            Note that $X\vdash_A Y$ is different from $\vdash_{A+X} Y$: the former is just the conditional proof of $Y$ from $X$ in $A$, while the latter is not a conditional proof in $A$ -- at most we may define the meta conditional proof $(\vdash_A X)\vdash_A (\vdash_A Y)$ as $\vdash_{A+X} Y$.
        \end{definition}

        \begin{definition}[Provable-equivalence class]
            \label{def:Provable-equivalence class}
            Let $A$ be a system and $x$ be a statement.

            The \emph{provable-equivalence class} of $x$ in $A$ is defined as $[x]_A=\{y \mid y \textup{ is provably equivalent to } x \textup{ in } A\}$.
        \end{definition}

        \begin{remark}
            \label{rem:Extension and restriction notations} In terms of the above definitions: an extension of a system is a super-system of that system; the (deductively) equivalent relation ``$=$'' between two systems can be used in a usual way (i.e. it has the properties of reflexivity, symmetry, transitivity and replacement) ; the extension operator ``$+$'' is commutative and associative. Let $A$, $B$ and $C$ be sets of statements, any of which can be a system, then:
            \begin{itemize}
                \item $A=B+C$ implies $\vdash_A B$ and $\vdash_A C$;
                \item $A=A+B$ is equivalent to $\vdash_A B$.
            \end{itemize}
        \end{remark}

        \subsection{Two super-systems of IRLA}
        \label{subsec:Two super-systems of IRLA}

        \begin{definition}[BL0 and BL1]
            \label{def:BL0 and BL1}
            Let $BL=(X,1,0,\neg,\wedge,\vee,\le)$ be the Boolean lattice described in Subsection \ref{subsec:Algebra IRLA}.  $BL0=(X,1,0,\neg,\wedge,\vee,\to,\oto,\le)$ and $BL1=(X,1,0,\neg,\wedge,\vee,\to,\oto,\le)$ are extensions of $BL$ defined as
            \begin{equation*}
                \begin{aligned}
                    BL0 &= BL+\{x\to y=\mu(x,y)\}\text{ where }\mu(x,y)=\begin{cases} 1,\;x\le y\\ 0,\;x\nleq y\end{cases}\text{ and }\\
                    BL1 &= BL+\{x\to y=\neg x\vee y\}
                \end{aligned}
            \end{equation*}
            respectively, with $\oto$ defined as usual by $x\oto y=(x\to y)\wedge(y\to x)$ for both of them.
        \end{definition}

        \begin{remark}[Logics BL0 and BL1 isomorphic to]\;
            \label{rem:Logics BL0 and BL1 isomorphic to}
            \begin{itemize}
                \item BL1 is isomorphic to classical (propositional) logic (CL), i.e. $BL1\cong CL$ in symbols.
                \item BL0 is isomorphic to a ``contra-classical'' logic, denoted as CL0, obtained by translating BL0 to a logic in terms of Table \ref{tab:Isomorphism between IRL and IRLA}, so $BL0\cong CL0$. The logic CL0 is contra-classical because the definitions of implication $x\to y=\mu(x,y)$ in BL0 and $x\to y=\neg x\vee y$ in BL1 are contradictory in \emph{the general case} (cf. Subsection \ref{subsec:D2I is underivable in IRL}), hence $BL0+BL1$ ($CL0+CL$) is inconsistent.
            \end{itemize}
        \end{remark}

        \begin{lemma}[Consistency of BL0 and BL1]
            \label{lemma:Consistency of BL0 and BL1}
            BL0 and BL1 are consistent.
        \end{lemma}

        \begin{proof}\;
            \begin{itemize}
                \item Consistency of BL0. For any given $(x,y)\in X^2$, $\mu(x,y)\in X$ uniquely exists (either $1$ or $0$), that is, $(x,y)\mapsto\mu(x,y)$ is (indeed) a binary operation on $X$, thus the new operations $\to$ and $\oto$ are well-defined. Therefore, BL0 is an extension of BL by definition, hence a conservative extension of BL, thus BL0 is consistent relative to BL (so they are equiconsistent as BL0 is a super-system of BL).
                \item Consistency of BL1. BL1 and CL are equiconsistent since $BL1\cong CL$.
            \end{itemize}
        \end{proof}

        Let
        $$\Delta=\{(x\oto y)\wedge z\le z[x\mapsto y],\;x=y\tto x\oto y=1,\;x\to y=x\wedge y\oto x\},$$
        then
        $$IRLA=BL+\Delta$$
        by definition (cf. Subsection \ref{subsec:Algebra IRLA}).

        \begin{lemma}[BL0 is a super-system of IRLA]
            \label{lemma:BL0 is a super-system of IRLA}
            $$\vdash_{BL0} IRLA.$$
        \end{lemma}

        \begin{proof}
            Since $\vdash_{BL0} BL$ and $IRLA=BL+\Delta$, we need only to prove $\vdash_{BL0} \Delta$, which is shown below.

            \begin{enumerate}
                \item Proof of $x=y\tto x\oto y=1$ (and $x\ne y\tto x\oto y=0$) in BL0. By definition of $\to$ and $\oto$ in BL0:

                If $x=y$, then $x\le y$ and $y\le x$, so $x\to y=1$ and $y\to x=1$, thus $x\oto y=(x\to y)\wedge(y\to x)=1\wedge 1=1$.

                If $x\ne y$, then $x\nleq y$ or $y\nleq x$, so $x\to y=0$ or $y\to x=0$, thus, let's say $x\to y=0$, $x\oto y=(x\to y)\wedge(y\to x)= 0\wedge (x\to y)=0$.

                \item Proof of $x\to y=x\wedge y\oto x$  in BL0.

                If $x\le y$ i.e. $x\wedge y=x$ by definition of $\le$, then $x\to y=1$ and $x\wedge y\oto x=1$, so $x\to y=x\wedge y\oto x$.

                If $x\nleq y$ i.e.  $x\wedge y\ne x$, then $x\to y=0$ and $x\wedge y\oto x=0$, so $x\to y=x\wedge y\oto x$.

                \item Proof of $(x\oto y)\wedge z\le z[x\mapsto y]$  in BL0.

                If $x=y$, then $(x\oto y)\wedge z\le z[x\mapsto y]$ becomes $1\wedge z\le z[x\mapsto x]$, i.e. $z\le z$, so it holds.

                If $x\ne y$, then $(x\oto y)\wedge z\le z[x\mapsto y]$ becomes $0\wedge z\le z[x\mapsto y]$, i.e. $0\le z[x\mapsto y]$, so it also holds.
            \end{enumerate}
        \end{proof}

        \begin{lemma}[BL1 is a super-system of IRLA]
            \label{lemma:BL1 is a super-system of IRLA}
            $$\vdash_{BL1} IRLA.$$
        \end{lemma}

        \begin{proof}
            Since $\vdash_{BL1} BL$ and $IRLA=BL+\Delta$, we need only to prove $\vdash_{BL1} \Delta$, which is shown below.

            \begin{enumerate}
                \item Proof of $x=y\tto x\oto y=1$ in BL1. If $x=y$, then $x\oto y = x\oto x = (x\to x)\wedge(x\to x) = x\to x = \neg x\vee x=1$.

                \item Proof of $x\to y=x\wedge y\oto x$ in BL1. We have $x\wedge y\oto x = (x\wedge y\to x)\wedge(x\to x\wedge y) = (\neg(x\wedge y)\vee x))\wedge(\neg x\vee (x\wedge y)) = (\neg x\vee\neg y\vee x)\wedge(\neg x\vee (x\wedge y)) = 1\wedge(\neg x\vee y) = \neg x\vee y = x\to y$.

                \item Proof sketch of $(x\oto y)\wedge z\le z[x\mapsto y]$ in BL1.

                1 Prove $x\oto y \le \neg x\oto \neg y$, $x\oto y \le z\circ x \oto z\circ y$ and $x\oto y \le x\circ z \oto y\circ z$ where $\circ\in\{\wedge,\vee,\to,\oto\}$;\\
                2 Prove $x\oto y\le z\oto z[x\mapsto y]$ from 1 by induction;\\
                3 Prove $(x\oto y)\wedge x\le y$;\\
                4 Prove $(x\oto y)\wedge z\le z[x\mapsto y]$ from 2 and 3.
            \end{enumerate}
        \end{proof}

        \begin{remark}[CL is a super-system of IRL]
            \label{rem:CL is a super-system of IRL}
            Actually, the idea in the proof of the above Lemma \ref{lemma:BL1 is a super-system of IRLA} can be used in proving directly that classical propositional logic CL is a supper-system of IRL. This means that the consistency of IRL can be proved before any derivations.
        \end{remark}

        \subsubsection*{IRLA-BL0 and IRLA-BL1 relations}

        \begin{metatheorem}[IRLA-BL0 and IRLA-BL1 relations]
            \label{metathrm:IRLA-BL0 and IRLA-BL1 relations}
            \begin{equation*}
                \begin{aligned}
                    IRLA+\{x\nleq y\tto x\to y=0\}=BL0;\\
                    IRLA+\{\neg x\vee y\le x\to y\}=BL1.
                \end{aligned}
            \end{equation*}
        \end{metatheorem}

        \begin{proof}
            By definition, $IRLA = BL+\Delta$, $BL0=BL+\{x\to y = \mu(x,y)\}$ where $\{x\to y = \mu(x,y)\}=\{x\le y\tto x\to y=1,\:x\nleq y\tto x\to y=0\}$, and $BL1=BL+\{x\to y = \neg x\vee y\}$ where $\{x\to y = \neg x\vee y\}=\{x\to y \le \neg x\vee y,\: \neg x\vee y\le x\to y\}$.
            \begin{enumerate}
                \item IRLA-BL0 relation. $\vdash_{IRLA} \{x\le y\tto x\to y=1\}$, i.e. $IRLA = IRLA + \{x\le y\tto x\to y=1\}$; $\vdash_{BL0} \Delta$, i.e. $BL0 + \Delta = BL0$. Thus
                \begin{equation*}
                    \begin{aligned}
                        IRLA & + \{x\nleq y\tto x\to y=0\} \\
                        & = IRLA + \{x\le y\tto x\to y=1\} + \{x\nleq y\tto x\to y=0\}\\
                        & = IRLA + \{x\le y\tto x\to y=1,\:x\nleq y\tto x\to y=0\}\\
                        & = IRLA + \{x\to y = \mu(x,y)\}\\
                        & = BL + \Delta + \{x\to y = \mu(x,y)\}\\
                        & = BL0 + \Delta\\
                        & = BL0.
                    \end{aligned}
                \end{equation*}

                \item IRLA-BL1 relation. $\vdash_{IRLA}\{x\to y\le \neg x\vee y\}$, i.e. $IRLA = IRLA +\{x\to y\le \neg x\vee y\}$; $\vdash_{BL1} \Delta$, i.e. $BL1 + \Delta = BL1$. Thus
                \begin{equation*}
                    \begin{aligned}
                        IRLA & +  \{\neg x\vee y \le x\to y\} \\
                        & = IRLA + \{x\to y\le \neg x\vee y\} + \{\neg x\vee y \le x\to y\}\\
                        & = IRLA + \{x\to y \le \neg x\vee y,\: \neg x\vee y\le x\to y\}\\
                        & = IRLA + \{x\to y = \neg x\vee y\}\\
                        & = BL + \Delta + \{x\to y = \neg x\vee y\}\\
                        & = BL1 + \Delta\\
                        & = BL1.
                    \end{aligned}
                \end{equation*}
            \end{enumerate}
        \end{proof}

        In contrast, Heyting lattice ($HL$) contains $\neg x\vee y\le x\to y$ but not $x\to y\le\neg x\vee y$, as a matter of fact $HL+\{x\to y\le\neg x\vee y\}=BL1$.

        \subsection{IRL-CL relationship}
        \label{subsec:IRL-CL relationship}

        \emph{Disjunction-to-Implication} (D2I) is the formula schema $\neg\phi\vee\psi\to(\phi\to\psi)$ in IRL, whose translation is the term $\neg x\vee y \to (x\to y)$ in IRLA.

        From the facts $IRL\cong IRLA$, $CL\cong BL1$ (cf. Subsection \ref{subsec:Isomorphism between IRL and IRLA}) and $IRLA+\{\neg x\vee y\le x\to y\}=BL1$ (cf. Subsection \ref{subsec:Two super-systems of IRLA}), we have $IRL+\{\neg\phi\vee\psi\to (\phi\to \psi)\}=CL$.

        It is shown that $D2I$ is generally independent of IRL in this section. Since $D2I$ is a theorem of the consistent CL, $\neg D2I$ is not derivable in CL, so $\neg D2I$ is not derivable in IRL either, as $\vdash_{CL} IRL$. Thus, it needs only to show that $D2I$ is underivable in IRL.

        \subsubsection*{Disjunction-to-Implication is underivable in IRL}
        \label{subsec:D2I is underivable in IRL}

        \begin{lemma}[D2I is not provable in BL0]
            \label{lemma:D2I is not provable in BL0}
            In $BL0=(X,1,0,\wedge,\vee,\to,\oto,\le)$, if there exists $c\in X$ such that $c\ne 1$ and $c\ne 0$, then $\neg x\vee y \le x\to y$ is not a theorem of BL0.
        \end{lemma}

        \begin{proof}
            From $c\ne 1$, $1\nleq c$, so $1\to c =0$ by definition; similarly, from $c\ne 0$, $c\to 0 = 0$. Let $t(x,y)=\neg x\vee y \to (x\to y)$. Then $t(1,c) = \neg 1\vee c \to (1\to c) = 0\vee c \to 0 = c\to 0 = 0 \ne 1$, so $t(x,y)=1$ is not a theorem of BL0 since its instance $t(1,c)\ne 1$.
        \end{proof}

        \begin{lemma}[Condition for D2I being a theorem of IRLA]
            \label{lemma:Condition for D2I being a theorem of IRLA}
            In $IRLA=(X,1,0,\wedge,\vee,\oto,\to,\le)$, $\neg x\vee y \to (x\to y)=1$ is a theorem if and only if there exists no $c\in X$ such that $c\ne 1$ and $c\ne 0$.
        \end{lemma}

        \begin{proof}
            Let $t(x,y)=\neg x\vee y \to (x\to y)$.

            \item Necessity. If there exists $c\in X$ such that $c\ne 1$ and $c\ne 0$, then $t(x,y)=1$ is not a theorem of BL0 in terms of Lemma \ref{lemma:D2I is not provable in BL0}, neither is it a theorem of IRLA, since $\vdash_{BL0} IRLA$.

            \item Sufficiency. If there exists no $c\in X$ such that $c\ne 1$ and $c\ne 0$, then for all $x,y\in X$, either $x=1$ or $x=0$ (i.e. $x,y\in\{1,0\}$). In IRLA, $t(1,1)=\neg 1\vee 1 \to (1\to 1) = 1 \to 1 =1$, $t(1,0)=\neg 1\vee 0 \to (1\to 0) = 0 \to 0 =1$, $t(0,1)=\neg 0\vee 1 \to (0\to 1) = 1 \to 1 =1$, and $t(0,0)=\neg 0\vee 0 \to (0\to 0) = 1 \to 1 =1$. Thus, for all $x,y\in X$, $\neg x\vee y \to (x\to y) = 1$, namely it is a theorem of IRLA in this case.
        \end{proof}

        \subsubsection*{The general case and the no-contingency case}

        Propositions are in the three types: tautology, contradiction, and contingency.

        There are two cases according to what types of propositional formulas the system contains, as described below.

        \begin{description}
            \item[The general case] The system contains all the three types of propositional formulas.
            \item[The no-contingency case] The system contains no contingencies.
        \end{description}

        Both IRL and CL always contain tautologies (e.g. $\phi\vee\neg\phi$) and contradictions (e.g. $\phi\wedge\neg\phi$). So the only difference between \emph{the general case} and \emph{the no-contingency case} is whether or not the system contains a contingency.

        For either IRL or CL: if it contains a primitive proposition that is contingent, then, of course, it contains a contingency; if it contains no primitive contingency (i.e. every primitive proposition is either a tautology or a contradiction), then it does not contain any contingencies because the set $\{1,0\}$ (representing all tautologies and contradictions) is closed under all the operations $\{\neg,\wedge,\vee,\to,\oto\}$ (e.g., all $1\to 1=1$, $1\to 0=0$, $0\to 1=1$ and $0\to 0=1$ hold in IRLA).

        Thus, for IRL and CL, it is \emph{the general case} if and only if \emph{there exists a primitive contingency} in the system, and it is \emph{the no-contingency case} if and only if \emph{there exists no primitive contingency} i.e. any formula is either tautological or contradictory.

        For example, suppose the set of primitive propositional symbols (the signature of propositional language) is simply the singleton $\{A\}$, and let $S_A$ be an instance of IRL or CL with this signature. Then: (a) $S_A$ contains all the three types of propositions if neither $A$ nor $\neg A$ is an axiom; (b) $S_A$ contains no contingencies (each formula is either a tautology or a contradiction) if either $A$ or $\neg A$ is an axiom.

        When not specified, we mean it is in the general case.

        \subsubsection*{Summary}

        Based on the previous results, the logic IRL has the relationship to classical logic (CL) as summarized below.

        \vspace{\baselineskip}
        In \emph{any cases} (no matter the general case or the no-contingency case):
        $$IRL+\{\neg\phi\vee\psi\to(\phi\to\psi)\}=CL, $$
        thus IRL is a sub-system of CL, so IRL is consistent (relative to CL), and the negation of Disjunction-to-Implication $\neg(\neg\phi\vee\psi\to(\phi\to\psi))$ is not a theorem of IRL.

        \begin{itemize}
            \item In \emph{the general case} (the system contains all types of propositions including contingencies, which is equivalent to that there is a primitive contingency):
            \begin{equation*}
                \begin{aligned}
                    IRL+\{&\neg\phi\vee\psi\to(\phi\to\psi)\} = CL \textup{, and }\\
                    &\neg\phi\vee\psi\to(\phi\to\psi) \textup{ is independent of } IRL.
                \end{aligned}
            \end{equation*}
            So IRL is a proper sub-system of CL in the general case.

            \item In \emph{the no-contingency case} (the system contains no contingencies, which is equivalent to that there is no primitive contingency):
            $$IRL = CL \textup{, and }\neg\phi\vee\psi\to(\phi\to\psi)\textup{ is a theorem of the system.}$$
        \end{itemize}

        \begin{remark}
            In the no-contingency case, any formula is either tautological or contradictory. However, it does not mean that the system is syntactically complete in this case, because we only know that every primitive proposition is either a tautology or a contradiction, but we do not know which one it is, and this is underivable.
        \end{remark}

        \subsection{Provable-equivalence class of Disjunction-to-Implication}
        \label{subsec:Provable-equivalence class of D2I}

        \begin{metatheorem}[Some members of the provable-equivalence class of Disjunction-to-Implication]
            \label{metathrm:Some members of the provable-equivalence class of Disjunction-to-Implication}
            The following formulas are provably equivalent in IRL.
            \begin{enumerate}
                \item $\phi\to (\top\to \phi)$, $\neg \phi\to (\phi\to \bot)$. (Lower bound of univariate implication)
                \item $\phi\to (\psi\to \phi)$, $\neg \phi\to( \phi\to \psi)$. (Generalized ECQ)
                \item $\neg \phi\vee \psi\to (\phi\to \psi)$. (Disjunction-to-Implication):
                \item $(\phi\wedge \psi\to \chi)\to (\phi\to(\psi\to \chi))$. (Exportation)
                \item $(\phi\to(\psi\to \chi))\oto(\psi\to(\phi\to \chi))$. (Commutativity of antecedents)
                \item $\phi\to ((\phi\to \psi)\to \psi)$. (Assertion)
                \item $(\phi\to \psi)\vee(\psi\to \phi)$. (Totality)
                \item $(\phi\to \psi)\vee(\phi\to\neg \psi)$. (Conditional excluded middle)
                \item $(\top\to \phi)\vee(\phi\to \bot)$. (Top-or-bottom)
            \end{enumerate}
        \end{metatheorem}
        Note that each of 1 and 2 is a pair of dual formulas that are provably equivalent in IRL.
        \begin{proof} It is done in the structure \textbf{1}$\to$\textbf{2}$\to$\textbf{3}$\to$\textbf{4}$\to$\textbf{5}$\to$\textbf{6}$\to$1; 2$\to$\textbf{7}$\to$\textbf{9}$\to$1; 2$\to$\textbf{8}$\to$9 as follows.\\
            \\
            First group \textbf{1}$\to$\textbf{2}$\to$\textbf{3}$\to$\textbf{4}$\to$\textbf{5}$\to$\textbf{6}$\to$1:\\
            (1) Lower bound of univariate implication $\to$ (2) Generalized ECQ:\\
            \begin{tabular}{lll}
                1&$\phi$&Premise\\
                2&$\top\to\phi$& 1: Lower bound of univariate implication\\
                3&$\psi\to\phi$& 2: Safe generalized ECQ
            \end{tabular}\\
            (2) Generalized ECQ $\to$ (3) Disjunction-to-Implication:\\
            \begin{tabular}{lll}
                1&$\neg\phi\to(\phi\to\psi)$& Generalized ECQ\\
                2&$\psi\to(\phi\to\psi)$& Generalized ECQ\\
                3&$\neg\phi\vee\psi\to(\phi\to \psi)$&1,2: Two-to-one\\
            \end{tabular}\\
            (3) Disjunction-to-Implication $\to$ (4) Exportation (with D2I, material implication $\phi\to\psi\otto\neg\phi\vee\psi$ holds):\\
            \begin{tabular}{lll}
                1&$\phi\wedge\psi\to\chi$&Premise\\
                2&$\neg(\phi\wedge\psi)\vee\chi$&1: Material implication\\
                3&$(\neg\phi\vee\neg\psi)\vee\chi$&2: De Morgan's law\\
                4&$\neg\phi\vee(\neg\psi\vee\chi)$&3: Associativity\\
                5&$\neg \phi\vee(\psi\to \chi)$&4: Material implication\\
                6&$\phi\to(\psi\to \chi)$&5: Material implication\\
            \end{tabular}\\
            (4) Exportation $\to$ (5) Commutativity of antecedents:\\
            \begin{tabular}{lll}
                1&$\phi\to(\psi\to\chi)$&Premise\\
                2&$\phi\wedge\psi\to\chi$&1: Importation\\
                3&$\psi\wedge\phi\to\chi$&2: Commutativity\\
                4&$\psi\to(\phi\to\chi)$&3: Exportation\\
            \end{tabular}\\
            (5) Commutativity of antecedents $\to$ (6) Assertion:\\
            \begin{tabular}{lll}
                1&$(\phi\to\psi)\to(\phi\to\psi)$&Reflexivity\\
                2&$\phi\to((\phi\to \psi)\to\psi)$&1: Commutativity of antecedents\\
            \end{tabular}\\
            (6) Assertion $\to$ (1) Lower bound of univariate implication:\\
            \begin{tabular}{lll}
                1&$(\phi\to\phi)\oto\top$&Reflexivity\\
                2&$\phi\to((\phi\to\phi)\to\phi)$&Assertion\\
                3&$\phi\to(\top\to\phi)$&1,2: Replacement
            \end{tabular}\\
            \\
            Second group 2$\to$\textbf{7}$\to$\textbf{9}$\to$1:\\
            (2) Generalized ECQ $\to$ (7) Totality:\\
            \begin{tabular}{lll}
                1&$\psi\to(\phi\to\psi)$& Generalized ECQ\\
                2&$\neg\psi\to(\psi\to\phi)$& Generalized ECQ\\
                3&$\psi\vee\neg\psi\to(\phi\to\psi)\vee(\psi\to\phi)$&1,2: Monotonicity of $\vee$ w.r.t $\to$\\
                4&$\psi\vee\neg \psi$&LEM\\
                5&$(\phi\to\psi)\vee(\psi\to\phi)$&3,4: Modus ponens\\
            \end{tabular}\\
            (7) Totality $\to$ (9) Top-or-bottom:\\
            \begin{tabular}{lll}
                1&$(\neg\phi\to\phi)\oto(\top\to\phi)$& Theorem \ref{thrm:Univariate implication expressions}\\
                2&$(\phi\to\neg\phi)\oto(\phi\to\bot)$& Theorem \ref{thrm:Univariate implication expressions}\\
                3&$(\neg\phi\to\phi)\vee(\phi\to\neg\phi)$&Totality\\
                4&$(\top\to \phi)\vee(\phi\to \bot)$& (1,2),3: Replacement
            \end{tabular}\\
            (9) Top-or-bottom $\to$ (1) Lower bound of univariate implication:\\
            \begin{tabular}{lll}
                1&$\phi$&Premise\\
                2&$(\phi\to\bot)\to\neg\phi$&Upper bound of univariate implication\\
                3&$\neg(\phi\to\bot)$&1,2: Modus tollens\\
                4&$(\top\to\phi)\vee(\phi\to\bot)$&Top-or-bottom\\
                5&$\top\to\phi$& 3,4: Disjunctive syllogism
            \end{tabular}\\
            \\
            Third group 2$\to$\textbf{8}$\to$9:\\
            (2) Generalized ECQ $\to$ (8) Conditional excluded middle:\\
            \begin{tabular}{lll}
                1&$\psi\to(\phi\to \psi)$& Generalized ECQ\\
                2&$\neg \psi\to(\phi\to\neg \psi)$& Generalized ECQ\\
                3&$\psi\vee\neg\psi\to(\phi\to\psi)\vee(\phi\to\neg\psi)$&1,2: Monotonicity of $\vee$ w.r.t. $\to$\\
                4&$\psi\vee\neg\psi$&LEM\\
                5&$(\phi\to\psi)\vee(\phi\to\neg \psi)$&3,4: Modus ponens
            \end{tabular}\\
            (8) Conditional excluded middle $\to$ (9) Top-or-bottom:\\
            \begin{tabular}{lll}
                1&$(\top\to\neg\phi)\oto(\phi\to\bot)$& Contraposition, Double negation\\
                2&$(\top\to\phi)\vee(\top\to\neg\phi)$&Conditional excluded middle\\
                3&$(\top\to\phi)\vee(\phi\to\bot)$&1,2: Replacement
            \end{tabular}\\
        \end{proof}

        From known members of the provable-equivalence class $[\neg\phi\vee\psi\to(\phi\to \psi)]_{IRL}$ of Disjunction-to-Implication such as formulas in Meta Theorem \ref{metathrm:Some members of the provable-equivalence class of Disjunction-to-Implication}, new member(s) of $[\neg\phi\vee \psi\to(\phi\to\psi)]_{IRL}$ can be decided by a chain of proofs like known1$\to$\textbf{new1}$\to$\textbf{new2}$\to$known2. The following meta theorem is an example.

        \begin{metatheorem}[Conjunction-to-equivalence]
            \label{metathrm:Conjunction-to-equivalence}
            The formula $\phi\wedge\psi\to(\phi\oto\psi)$ (Conjunction-to-equivalence) is a member of $[\neg\phi\vee\psi\to (\phi\to\psi)]_{IRL}$.
        \end{metatheorem}

        \begin{proof} It is done in the chain Generalized ECQ (known) $\to$ \textbf{Conjunction-to-equivalence} (new) $\to$ Lower bound of univariate implication (known):

            Generalized ECQ $\to$ \textbf{Conjunction-to-equivalence}:\\
            \begin{tabular}{lll}
                1&$\phi\wedge\psi$& Premise\\
                2&$\phi$& 1: Conjunction elimination\\
                3&$\psi\to\phi$& 2: Generalized ECQ\\
                4&$\phi\to\psi$& 1--3: Similarly\\
                5&$\phi\oto\psi$& 3,4: Bi-implication introduction\\
            \end{tabular}\\

            \textbf{Conjunction-to-equivalence} $\to$ Lower bound of univariate implication:\\
            \begin{tabular}{lll}
                1&$\phi$& Premise\\
                2&$\top\wedge\phi$&1: Identity element\\
                3&$\top\oto\phi$&2: Conjunction-to-equivalence\\
                4&$\top\to\phi$&3: Bi-implication elimination\\
            \end{tabular}\\
        \end{proof}

        All members of $[\neg\phi\vee\psi\to(\phi\to\psi)]_{IRL}$ are not provable in IRL, unless the extra assumption of no-contingency is made.


        \section{The decision problem}

        Since $IRL\cong IRLA$ (cf. Subsection \ref{subsec:Isomorphism between IRL and IRLA}), a formula schema $F(\phi_1,...,\phi_n)$ in IRL corresponds to the term $t(x_1,...,x_n)$ in IRLA (cf. Table \ref{tab:Isomorphism between IRL and IRLA}).

        In a sentence $t(x_1,...,x_n)=1$ of an algebraic system with an underlying set $X$, all variables are implicitly universally quantified over $X$. So the decision problem considered in $IRLA=(X,1,0,\neg,\wedge,\vee,\oto,\to,\le)$ is, given any term $t(x_1,...,x_n)$, decide whether or not $t(x_1,...,x_n)=1$ for all $(x_1,...,x_n)\in X^n$ -- we may simply say whether or not ``$t(x_1,...,x_n)$ is valid'' for convenience.

        A complete decision procedure is not found for IRL/IRLA in this work. A conjecture: IRL/IRLA is undecidable. If it is indeed undecidable, this is a reflection of the nature of implication relation, which comes from the original meaning of implication.

        \vspace{\baselineskip}
        If $t(x_1,...,x_n)$ is valid in IRLA, it is valid in both its super-systems BL0 and BL1 (cf. Subsection \ref{subsec:Two super-systems of IRLA}). So if $t(x_1,...,x_n)$ is invalid in either BL0 or BL1, it is invalid in IRLA. Thus, BL0 and BL1 can be used as two basic validity ``filters'' for IRLA. For example, $\neg (\neg x\vee y\to (x\to y))$ is invalid in BL1, so it is invalid in IRLA.

        Since BL1 represents classical propositional logic (cf. Subsection \ref{subsec:Two super-systems of IRLA}), many ready decision procedures can be utilized, such as the basic methods of truth table, truth tree, and SAT algorithms.

        The following is a simple validity-checking method for BL0.

        \subsubsection*{Validity checking in BL0 on a minimal set}

        Consider the validity-checking problem in $BL0=(X,1,0,\neg,\wedge,\vee,\to,\le)$: given any expression $t(x_1,...,x_n)$, whether or not $t(x_1,...,x_n)$ is valid, i.e. if $t(x_1,...,x_n)=1$ for all $(x_1,...,x_n)\in X^n$?

        Suppose there exists $c\in X$ such that $0<c<1$ (i.e. there exists a contingency in the system). Consider validity checking on the set $\{0,c,1\}$, a minimal subset of $X$ that can reflect \emph{the general case} (there are all the three types of propositions: tautology, contradiction and contingency).

        In BL0, the closure of the set $\{0,c,1\}\subset X$ under all operations $\{\neg,\wedge,\vee,\to\}$ is the set $\{0,c,\neg c,1\}\subseteq X$, as shown in Table \ref{tab:Three T-set values closed at four}.

        \begin{table}[H] 
            \centering
            \small
            \setlength{\tabcolsep}{18pt} 
            \renewcommand{\arraystretch}{1.0} 
            \begin{threeparttable}
                \caption{Closure of the set $\{0,c,1\}$ in BL0 \tnote{*}}
                \label{tab:Three T-set values closed at four}
                \begin{tabular}{cc|c|c|c|c}
                    \hline
                    $x$ & $y$ & $\neg x$ & $x\wedge y$ & $x\vee y$ & $x\to y$\\
                    \hline
                    $0$ & $0$ & $1$ & $0$ & $0$ & $1$ \\
                    $0$ & $c$ & $1$ & $0$ & $c$ & $1$ \\
                    $0$ & $1$ & $1$ & $0$ & $1$ & $1$ \\
                    \hline
                    $c$ & $0$ & $\neg c$ & $0$ & $c$ & $0$ \\
                    $c$ & $c$ & $\neg c$ & $c$ & $c$ & $1$ \\
                    $c$ & $1$ & $\neg c$ & $c$ & $1$ & $1$ \\
                    \hline
                    $1$ & $0$ & $0$ & $0$ & $1$ & $0$ \\
                    $1$ & $c$ & $0$ & $c$ & $1$ & $0$ \\
                    $1$ & $1$ & $0$ & $1$ & $1$ & $1$ \\
                    \hline
                \end{tabular}
                \begin{tablenotes}
                    \item[*] In BL0, the closure of $\{0,c,1\}$ under $\{\neg,\wedge,\vee,\to\}$ is $\{0,c,\neg c,1\}$.
                \end{tablenotes}
            \end{threeparttable}
        \end{table}

        Thus, in BL0, any expression $t(x_1,...,x_n)\in \{0,c,\neg c,1\}$ for all $(x_1,...,x_n)\in \{0,c,1\}^n$. Note that $\neg c$ is also a contingency, i.e. $0<\neg c<1$, since, only in a Boolean algebra, $0<x<1$ if and only if $0< \neg x <1$ for all $x$ in its underlying set.

        \vspace{\baselineskip}
        If $t(x_1,...,x_n)\ne 1$ for some $(x_1,...,x_n)\in \{0,c,1\}^n$, then $t(x_1,...,x_n)\ne 1$ for some $(x_1,...,x_n)\in X^n$, so the answer to the validity question of $t(x_1,...,x_n)$ is ``No''.

        Otherwise $t(x_1,...,x_n)=1$ for all $(x_1,...,x_n)\in \{0,c,1\}^n$, the answer to the validity question of $t(x_1,...,x_n)$ is ``Unknown'', i.e., it is possible for $t$ is valid and it is also possible for $t$ is invalid, since, although $t=1$ for all $(x_1,...,x_n)\in \{0,c,1\}^n$, it is possible that $t\ne 1$ for some $(x_1,...,x_n)\in X^n\setminus \{0,c,1\}^n$.

        In terms of these facts, we have a simple algorithm for the (in)validity-checking problem as described below.

        \vspace{\baselineskip}
        \textbf{Validity checking by exhaustion over $\{0,c,1\}\subset X$}\\

        In BL0:
        \begin{enumerate}
            \item Input an expression $t(x_1,...,x_n)$;
            \item For each $(x_1,...,x_n)\in \{0,c,1\}^n$, compute $t(x_1,...,x_n)$;
            \item At any step:
            \begin{itemize}
                \item If $t(u_1,...,u_n)\ne 1$ (i.e. $t=0,c,\neg c$) for some $(u_1,...,u_n)\in\{0,c,1\}^n$, $t$ is \textbf{Invalid} with a countermodel $(x_1,...,x_n)=(u_1,...,u_n)$;
                \item If $t(x_1,...,x_n)\ne 0$ (i.e. $\neg t\ne 1$) for some $(u_1,...,u_n)\in\{0,c,1\}^n$, $\neg t$ is invalid with a countermodel $(x_1,...,x_n)=(u_1,...,u_n)$;
                \item If both $t$ and $\neg t$ are invalid, $t$ is \textbf{independent} of the system;
            \end{itemize}
            \item If $t(x_1,...,x_n)=1$ for all $(x_1,...,x_n)\in\{0,c,1\}^n$, the validity of $t$ is \textbf{Unkown} (but $\neg t$ is invalid);
            \item If $t(x_1,...,x_n)=0$ for all $(x_1,...,x_n)\in\{0,c,1\}^n$, $t$ is \textbf{Invalid} (but the validity of $\neg t$ is unknown).
        \end{enumerate}

        For IRLA as a sub-system of BL0:
        \begin{itemize}
            \item $t$ is invalid in IRLA if it is invalid in BL0.
            \item $t$ is independent of IRLA if it is independent of BL0.
            \item The validity of $t$ is unknown in IRLA when its validity is unknown in BL0, or $t$ is valid in BL0.
        \end{itemize}

        This simple exhaustive search algorithm, with the worst-case complexity $3^n$, is of course \emph{incomplete}, since it is not for all inputs that the algorithm can give ``Yes'' or ``No'' answer (for some inputs the answer is ``Unknown''). But it is useful as a basic validity ``filter'' for IRLA/IRL.

        \begin{example}
            Is $t(x)=(1\to x)\to x$ valid? For $x\in \{0,c,1\}$:
            \begin{itemize}
                \item[1.] $x=0$: $t(0)=(1\to 0)\to 0=0\to 0=1$.
                \item[2.] $x=c$: $t(c)=(1\to c)\to c=0\to c=1$.
                \item[3.] $x=1$: $t(1)=(1\to 1)\to 1=1\to 1=1$.
            \end{itemize}
            The result is $t(x)=(1\to x)\to x=1$ for all $x$ in $\{0,c,1\}$, so the answer to the validity question of $(1\to x)\to x$ is ``Unknown'' (i.e. $(1\to x)\to x$ is possibly valid and possibly invalid as well).
        \end{example}

        \begin{example}
            Is $t(x)=x\to (1\to x)$ valid? For $x\in \{0,c,1\}$:
            \begin{itemize}
                \item[1.] $x=0$: $t(0)=0\to (1\to 0)=0\to 0=1$.
                \item[2.] $x=c$: $t(c)=c\to(1\to c)=c\to 0=0$. (a countermodel)
            \end{itemize}
            So $t(x)=x\to (1\to x)$ is invalid with a countermodel $x=c$. And since $\neg t(0)=\neg 1=0$, $\neg t(x)$ is also invalid, thus $t(x)$ is independent of IRLA.
        \end{example}

        \begin{example}
            Is $t(x,y)=(x\to y)\to \neg x\vee y$ valid? For $(x,y)\in \{0,c,1\}^2$:
            \begin{itemize}
                \item[1.] $(x,y)=(0,0)$: $t(0,0)=(0\to 0)\to \neg 0 \vee 0=1\to 1=1$.
                \item[2.] $(x,y)=(0,c)$: $t(0,c)=(0\to c)\to \neg 0\vee c=1\to 1=1$.
                \item[3.] $(x,y)=(0,1)$: $t(0,1)=(0\to 1)\to \neg 0\vee 1=1\to 1=1$.
                \item[4.] $(x,y)=(c,0)$: $t(c,0)=(c\to 0)\to \neg c\vee 0=0\to \neg c=1$.
                \item[5.] $(x,y)=(c,c)$: $t(c,c)=(c\to c)\to \neg c\vee c=1\to 1=1$.
                \item[6.] $(x,y)=(c,1)$: $t(c,1)=(c\to 1)\to \neg c\vee 1=1\to 1=1$.
                \item[7.] $(x,y)=(1,0)$: $t(1,0)=(1\to 0)\to \neg 1\vee 0=0\to 0=1$.
                \item[8.] $(x,y)=(1,c)$: $t(1,c)=(1\to c)\to \neg 1\vee c=0\to c=1$.
                \item[9.] $(x,y)=(1,1)$: $t(1,1)=(1\to 1)\to \neg 1\vee 1=1\to 1=1$.
            \end{itemize}
            The result is $t(x,y)=(x\to y)\to \neg x\vee y=1$ for all $(x,y)$ in $\{0,c,1\}^2$, so the answer to the validity question of $(x\to y)\to \neg x\vee y$ is ``Unknown'' (i.e. $(x\to y)\to \neg x\vee y$ is possibly valid and possibly invalid as well).
        \end{example}

        \begin{example}
            Is $t(x,y)=\neg x \vee y\to (x\to y)$ valid? For $(x,y)\in \{0,c,1\}^2$:
            \begin{itemize}
                \item[1.] $(x,y)=(0,0)$: $t(0,0)=\neg 0\vee 0\to (0\to 0)=1\to 1=1$.
                \item[2.] $(x,y)=(0,c)$: $t(0,c)=\neg 0\vee c\to (0\to c)=1\to 1=1$.
                \item[3.] $(x,y)=(0,1)$: $t(0,1)=\neg 0\vee 1\to (0\to 1)=1\to 1=1$.
                \item[4.] $(x,y)=(c,0)$: $t(c,0)=\neg c\vee 0\to (c\to 0)=\neg c\to 0=1\to c=0$. (a countermodel)
            \end{itemize}
            So $t(x,y)=\neg x\vee y\to (x\to y)$ is invalid with a countermodel $(x,y)=(c,0)$. And as we also have e.g.  $t(0,0)=1$, $t(x,y)$ is independent of IRLA.
        \end{example}

        \begin{example}
            \label{exam:Peirce's law}
            Are $t(x,y)=((x\to y)\to x)\to x$ and its negation $\neg t$ valid? For $(x,y)\in \{0,c,1\}^2$:
            \begin{itemize}
                \item[1.] $(x,y)=(0,0)$: $t(0,0)=((0\to 0)\to 0)\to 0=(1\to 0)\to 0=0\to 0=1$. (a countermodel for $\neg t$)
                \item[2.] $(x,y)=(0,c)$: $t(0,c)=((0\to c)\to 0)\to 0=(1\to 0)\to 0=0\to 0=1$.
                \item[3.] $(x,y)=(0,1)$: $t(0,1)=((0\to 1)\to 0)\to 0=(1\to 0)\to 0=0\to 0=1$.
                \item[4.] $(x,y)=(c,0)$: $t(c,0)=((c\to 0)\to c)\to c=(0\to c)\to c=1\to c=0$. (a countermodel for $t$)
            \end{itemize}
            So, $t(x,y)=((x\to y)\to x)\to x$ is invalid with a countermodel $(x,y)=(c,0)$, and $\neg t(x,y)$ is also invalid with a countermodel $(x,y)=(0,0)$. Thus $((x\to y)\to x)\to x$ is independent of IRLA.
        \end{example}

        Example \ref{exam:Peirce's law} is actually about the validity of Peirce’s law, and the result indicates that Peirce’s law, like Disjunction-to-Implication, is independent of IRL. It can be shown that, in IRL, Peirce’s law is a provable consequence of Disjunction-to-Implication, but (very likely) not vice versa.


        \section{Natural deduction}

        The natural deduction for IRL is given by the set of rules below:

        \begin{tabular}{rlll}
            1. & $(\phi\vdash\psi)\tto(\vdash\phi\to\psi)$ & $\to$I & CP (\emph{Restricted}, cf. Remark \ref{rem:Restriction on CP} below) \\
            2. & $\{\phi\to\psi,\phi\}\tto\psi$ & $\to$E & MP \\
            3. & $(\phi\to\psi)\tto (\psi\to\chi)\to(\phi\to\chi)$ &  & Suffixing \\
            4. & $\{\phi,\psi\}\tto\phi\wedge\psi$ & $\wedge$I & Conjunction \\
            5. & $\phi\wedge\psi\tto\phi$, $\phi\wedge\psi\tto\psi$ & $\wedge$E & Simplification \\
            6. & $\{\chi\to\phi,\chi\to\psi\}\tto \chi\to\phi\wedge\psi$ &  & ``One-to-two'' \\
            7. & $\phi\tto\phi\vee\psi$, $\psi\tto\phi\vee\psi$ & $\vee$I & Addition \\
            8. & $\{\phi\to\chi,\psi\to\chi\}\tto \phi\vee\psi\to\chi$ &  & ``Two-to-one'' ($\vee$E follows from this)\\
            9. & $\phi\wedge(\psi\vee\chi)\tto(\phi\wedge\psi)\vee(\phi\wedge\chi)$ &  & Distribution \\
            10. & $\top\tto\phi\vee\neg\phi$ & $\neg$I/$\top$E& LEM \\
            11. & $\phi\wedge\neg\phi\tto\bot$ & $\neg$E/$\bot$I & LNC \\
            12. & $\phi\tto\top$ & $\top$I & Dual of ECQ \\
            13. & $\bot\tto\phi$ & $\bot$E & ECQ \\
            14. & $\{\phi\to\psi,\psi\to\phi\}\tto \phi\oto\psi$ & $\oto$I & Bi-implication introduction\\
            15. & $\phi\oto\psi \tto \phi\to\psi$, $\phi\oto\psi\tto \psi\to\phi$ & $\oto$E & Bi-implication elimination
        \end{tabular}

        \begin{remark}[Restriction on CP]
            \label{rem:Restriction on CP}
            By the definition of proof in Subsection \ref{subsec:Definitions and notations 1}, any member of a proof sequence must be either a theorem (or a theorem indicated in an outer layer), or a result from a subset of the premises and previous members in its own layer by a rule (or by a rule indicated in an outer layer), not allowing any other cases. Attention should be paid to this restriction for nested proofs, especially for nested conditional proofs.
        \end{remark}

        \subsubsection*{Equivalence to Hilbert system}

        The natural deduction for IRL is equivalent to its Hilbert-style system defined in Section \ref{sec:Propositional logic IRL}.

        The primitive rules listed above of the natural deduction are already proved in its Hilbert system (cf. Section \ref{sec:Propositional logic IRL}).

        On the other hand, all axioms and primitive rules of the Hilbert-style IRL can be derived in this natural deduction. Some derivations are shown as examples in the following.

        \begin{derivedrule}[Commutativity of conjunction]
            $$\phi\wedge\psi\tto\psi\wedge\phi.$$
        \end{derivedrule}

        \begin{proof}\;\\
            \begin{tabular}{lll}
                1&\qquad$\phi\wedge\psi$&Premise\\
                2&\qquad$\phi$&1: $\wedge$E\\
                3&\qquad$\psi$&1: $\wedge$E\\
                4&\qquad$\psi\wedge\phi$&3,2: $\wedge$I\\
                5&$\phi\wedge\psi\vdash\psi\wedge\phi$&1--4: Definition of conditional proof\\
                6&$\phi\wedge\psi\tto\psi\wedge\phi$&5: Definition of derivable rule
            \end{tabular}\\
        \end{proof}

        The last two lines will be omitted in the subsequent proofs.

        \begin{derivedrule}[Associativity of conjunction]
            $$(\phi\wedge\psi)\wedge\chi\tto\phi\wedge(\psi\wedge\chi).$$
        \end{derivedrule}

        \begin{proof}\;\\
            \begin{tabular}{lll}
                1&$(\phi\wedge\psi)\wedge\chi$&Premise\\
                2&$\phi\wedge\psi$&1: $\wedge$E\\
                3&$\phi$&2: $\wedge$E\\
                4&$\psi$&2: $\wedge$E\\
                5&$\chi$&1: $\wedge$E\\
                6&$\psi\wedge\chi$&4,5: $\wedge$I\\
                7&$\phi\wedge(\psi\wedge\chi)$&3,6: $\wedge$I
            \end{tabular}\\
        \end{proof}

        With the rules of commutativity and associativity of conjunction, conjunction introduction and elimination can be generalized to $\{\phi_1,...,\phi_n\}\tto\phi_1\wedge\cdots\wedge\phi_n$ and $\phi_1\wedge\cdots\wedge\phi_n\tto\phi_i$ for any $n\ge 2$, respectively.

        \begin{metarule}[Equivalence of the three forms]
            $$(\vdash\phi_1\wedge\cdots\wedge\phi_n\to\psi)\otto(\{\phi_1,...,\phi_n\}\vdash\psi)\otto(\{\phi_1,...,\phi_n\}\tto\psi).$$
        \end{metarule}

        \begin{proof}
            This result is dependent of (cf. Subsection \ref{subsec:Elementary derived results}): (a) the definition of proof (conditional proof and unconditional proof, cf. Subsection \ref{subsec:Definitions and notations 1}), (b) the definition of derivable rule (cf. Subsection \ref{subsec:Definitions and notations 1}), (c) the rule of conditional proof (it is primitive), (d) modus ponens (it is primitive), and (e) rules of commutativity and associativity of conjunction (they are derived above).
        \end{proof}

        With this result, we need not distinguish between a rule $\{\phi_1,...,\phi_n\}\tto\psi$ and a theorem $\phi_1\wedge\cdots\wedge\phi_n\to\psi$ in proving and using them. Specifically, from a conditional proof $\{\phi_1,...,\phi_n\}\vdash\psi$ or an unconditional proof $\vdash\phi_1\wedge\cdots\wedge\phi_n\to\psi$, we obtain both the rule $\{\phi_1,...,\phi_n\}\tto\psi$ and the theorem $\phi_1\wedge\cdots\wedge\phi_n\to\psi$; we can use the rule $\{\phi_1,...,\phi_n\}\tto\psi$ as the theorem $\phi_1\wedge\cdots\wedge\phi_n\to\psi$, and vice versa.

        \begin{theorem}[Reflexivity]
            $$\phi\to\phi.$$
        \end{theorem}

        \begin{proof}\;\\
            \begin{tabular}{lll}
                1&$\phi$&Premise\\
                2&$\phi\vee\phi$&1: $\vee$I\\
                3&$\phi\wedge(\phi\vee\phi)$&1,2: $\wedge$I\\
                4&$\phi$&3: $\wedge$E
            \end{tabular}\\
        \end{proof}

        \begin{derivedrule}[Transitivity]
            $$\{\phi\to\psi,\psi\to\chi\}\tto \phi\to\chi.$$
        \end{derivedrule}

        \begin{proof}\;\\
            \begin{tabular}{lll}
                1&$\phi\to\psi$&Premise\\
                2&$\psi\to\chi$&Premise\\
                3&$(\psi\to\chi)\to(\phi\to\chi)$&1: Suffixing\\
                4&$\phi\to\chi$&2,3: MP
            \end{tabular}\\
        \end{proof}

        \begin{derivedrule}[Proof by exhaustion, a form of $\vee$E]
            $$\{\phi\vee\psi,\phi\to\chi,\psi\to\chi\}\tto\chi.$$

            This is the usual rule of $\vee$E in natural deduction for classical logic.
        \end{derivedrule}

        \begin{proof}\;\\
            \begin{tabular}{lll}
                1&$\phi\vee\psi$&Premise\\
                2&$\phi\to\chi$&Premise\\
                3&$\psi\to\chi$&Premise\\
                4&$\phi\vee\psi\to\chi$&2,3: ``Two-to-one''\\
                5&$\chi$&1,4: MP
            \end{tabular}\\
        \end{proof}

        \begin{derivedrule}[Distribution of $\wedge$ over $\vee$, converse]
            $$(\phi\wedge\psi)\vee(\phi\wedge\chi)\tto\phi\wedge (\psi\vee\chi).$$
        \end{derivedrule}

        \begin{proof}\;\\
            \begin{tabular}{lll}
                1&$(\phi\wedge\psi)\vee(\phi\wedge\chi)$&Premise\\
                2&\qquad$\phi\wedge\psi$&Assumption\\
                3&\qquad$\phi$&2: $\wedge$E\\
                4&\qquad$\psi$&2: $\wedge$E\\
                5&\qquad$\psi\vee\chi$&4: $\vee$I\\
                6&\qquad$\phi\wedge(\psi\vee\chi)$&3,5: $\wedge$I\\
                7&$\phi\wedge\psi\to\phi\wedge(\psi\vee\chi)$&2--6: CP\\
                8&$\phi\wedge\chi\to\phi\wedge(\psi\vee\chi)$&2--7: Similarly\\
                9&$\phi\wedge(\psi\vee\chi)$&1,7,8: Proof by exhaustion
            \end{tabular}\\
        \end{proof}

        \begin{derivedrule}[Distribution of $\vee$ over $\wedge$]
            $$\phi\vee(\psi\wedge\chi)\tto(\phi\vee\psi)\wedge(\phi\vee\chi).$$
        \end{derivedrule}

        \begin{proof}\;\\
            \begin{tabular}{lll}
                1&$\phi\vee(\psi\wedge\chi)$&Premise\\
                2&\qquad$\phi$&Assumption\\
                3&\qquad$\phi\vee\psi$&2: $\vee$I\\
                4&\qquad$\phi\vee\chi$&2: $\vee$I\\
                5&\qquad$(\phi\vee\psi)\wedge(\phi\vee\chi)$&3,4: $\wedge$I\\
                6&$\phi\to(\phi\vee\psi)\wedge(\phi\vee\chi)$&2--5: CP\\
                7&\qquad$\psi\wedge\chi$&Assumption\\
                8&\qquad$\psi$&7: $\wedge$E\\
                9&\qquad$\chi$&7: $\wedge$E\\
                10&\qquad$\phi\vee\psi$&8: $\vee$I\\
                11&\qquad$\phi\vee\chi$&9: $\vee$I\\
                12&\qquad$(\phi\vee\psi)\wedge(\phi\vee\chi)$&10,11: $\wedge$I\\
                13&$\psi\wedge\chi\to(\phi\vee\psi)\wedge(\phi\vee\chi)$&7--12: CP\\
                14&$(\phi\vee\psi)\wedge(\phi\vee\chi)$&1,6,13: Proof by exhaustion
            \end{tabular}\\
        \end{proof}

        \begin{derivedrule}[Distribution of $\vee$ over $\wedge$, converse]
            $$(\phi\vee\psi)\wedge(\phi\vee\chi)\tto\phi\vee(\psi\wedge\chi).$$
        \end{derivedrule}

        \begin{proof}\;\\
            \begin{tabular}{lll}
                1&$(\phi\vee\psi)\wedge(\phi\vee\chi)$&Premise\\
                2&$((\phi\vee\psi)\wedge\phi)\vee((\phi\vee\psi)\wedge\chi)$&1: Distribution (of $\wedge$ over $\vee$, primitive)\\
                3&\qquad$(\phi\vee\psi)\wedge\phi$&Assumption\\
                4&\qquad$\phi$&3: $\wedge$E\\
                5&\qquad$\phi\vee(\psi\wedge\chi)$&4: $\vee$I\\
                6&$(\phi\vee\psi)\wedge\phi\to\phi\vee(\psi\wedge\chi)$&3--5: CP\\
                7&\qquad$(\phi\vee\psi)\wedge\chi$&Assumption\\
                8&\qquad$(\phi\wedge\chi)\vee(\psi\wedge\chi)$&7: Distribution (of $\wedge$ over $\vee$, primitive)\\
                9&\qquad\qquad$\phi\wedge\chi$&Assumption\\
                10&\qquad\qquad$\phi$&9: $\wedge$E\\
                11&\qquad\qquad$\phi\vee(\psi\wedge\chi)$&10: $\vee$I\\
                12&\qquad$\phi\wedge\chi\to\phi\vee(\psi\wedge\chi)$&9--11: CP\\
                13&\qquad$\psi\wedge\chi\to\phi\vee(\psi\wedge\chi)$&$\vee$I\\
                14&\qquad$\phi\vee(\psi\wedge\chi)$&8,12,13: Proof by exhaustion\\
                15&$(\phi\vee\psi)\wedge\chi\to\phi\vee(\psi\wedge\chi)$&7--14: CP\\
                16&$\phi\vee(\psi\wedge\chi)$&2,6,15: Proof by exhaustion
            \end{tabular}\\
        \end{proof}

        \begin{derivedrule}[Compatibility of $\to$ with $\wedge$]
            $$\phi\to\psi \tto \phi\wedge\chi\to\psi\wedge\chi.$$
        \end{derivedrule}

        \begin{proof}\;\\
            \begin{tabular}{lll}
                1&$\phi\to\psi$&Premise\\
                2&$\phi\wedge\chi\to\phi$&$\wedge$E\\
                3&$\phi\wedge\chi\to\psi$&2,1: Transitivity\\
                4&$\phi\wedge\chi\to\chi$&$\wedge$E\\
                5&$\phi\wedge\chi\to\psi\wedge\chi$&3,4: ``One-to-two''
            \end{tabular}\\
        \end{proof}

        \begin{derivedrule}[Compatibility of $\to$ with $\vee$]
            $$\phi\to\psi \tto \phi\vee\chi\to\psi\vee\chi.$$
        \end{derivedrule}

        \begin{proof}\;\\
            \begin{tabular}{lll}
                1&$\phi\to\psi$&Premise\\
                2&$\psi\to\psi\vee\chi$&$\vee$I\\
                3&$\phi\to\psi\vee\chi$&1,2: Transitivity\\
                4&$\chi\to\psi\vee\chi$&$\vee$I\\
                5&$\phi\vee\chi\to\psi\vee\chi$&3,4: ``Two-to-one''
            \end{tabular}\\
        \end{proof}

        \begin{theorem}[$\top$ is a theorem]
            $$\top.$$
        \end{theorem}

        \begin{proof}\;\\
            \begin{tabular}{lll}
                1&$(\phi\to\phi)\to\top$&$\top$I\\
                2&$\phi\to\phi$&Reflexivity\\
                3&$\top$&1,2: MP
            \end{tabular}\\
        \end{proof}

        \begin{theorem}[LEM, theorem form]
            $$\phi\vee\neg\phi.$$
        \end{theorem}

        \begin{proof}\;\\
            \begin{tabular}{lll}
                1&$\top\to\phi\vee\neg\phi$&$\top$E\\
                2&$\top$&$\top$ is a theorem\\
                3&$\phi\vee\neg\phi$&1,2: MP
            \end{tabular}\\
        \end{proof}

        \begin{theorem}[ECQ, a form without $\bot$]
            $$\phi\wedge\neg\phi\to\psi.$$
        \end{theorem}

        \begin{proof}\;\\
            \begin{tabular}{lll}
                1&$\phi\wedge\neg\phi$&Premise\\
                2&$\bot$&1: $\bot$I\\
                3&$\psi$&2: $\bot$E
            \end{tabular}\\
        \end{proof}

        \begin{theorem}[Dual of ECQ, a form without $\top$]
            $$\phi\to\psi\vee\neg\psi.$$
        \end{theorem}

        \begin{proof}\;\\
            \begin{tabular}{lll}
                1&$\phi$&Premise\\
                2&$\top$&1: $\top$I\\
                3&$\psi\vee\neg\psi$&2: $\top$E
            \end{tabular}\\
        \end{proof}

        \begin{derivedrule}[DNI]
            $$\phi\tto\neg\neg\phi.$$
        \end{derivedrule}

        \begin{proof}\;\\
            \begin{tabular}{lll}
                1&$\phi$&Premise\\
                2&$\neg\phi\vee\neg\neg\phi$&LEM\\
                3&$\phi\wedge(\neg\phi\vee\neg\neg\phi)$&1,2: $\wedge$I\\
                4&$(\phi\wedge\neg\phi)\vee(\phi\wedge\neg\neg\phi)$&3: Distribution\\
                5&$\phi\wedge\neg\phi\to\neg\neg\phi$&ECQ\\
                6&$\phi\wedge\neg\neg\phi\to\neg\neg\phi$&$\wedge$E\\
                7&$\neg\neg\phi$&4,5,6: Proof by exhaustion
            \end{tabular}\\
        \end{proof}

        \begin{derivedrule}[DNE]
            $$\neg\neg\phi\tto\phi.$$
        \end{derivedrule}

        \begin{proof}\;\\
            \begin{tabular}{lll}
                1&$\neg\neg\phi$&Premise\\
                2&$\phi\vee\neg\phi$&LEM\\
                3&$\neg\neg\phi\wedge(\phi\vee\neg\phi)$&1,2: $\wedge$I\\
                4&$(\neg\neg\phi\wedge\phi)\vee(\neg\neg\phi\wedge\neg\phi)$&3: Distribution\\
                5&$\neg\neg\phi\wedge\phi\to\phi$&$\wedge$E\\
                6&$\neg\neg\phi\wedge\neg\phi\to\phi$&ECQ\\
                7&$\phi$&4,5,6: Proof by exhaustion
            \end{tabular}\\
        \end{proof}

        \begin{derivedrule}[Contraposition, negation introduced]
            $$\phi\to\psi\tto \neg\psi\to\neg\phi.$$
        \end{derivedrule}

        \begin{proof}\;\\
            \begin{tabular}{lll}
                1&$\phi\to\psi$&Premise\\
                2&$\neg\psi\to\neg\psi$&Reflexivity\\
                3&$\neg\psi\to\phi\vee\neg\phi$&Dual of ECQ\\
                4&$\neg\psi\to\neg\psi\wedge(\phi\vee\neg\phi)$&2,3: ``One-to-two''\\
                5&$\phi\vee\neg\phi\to\psi\vee\neg\phi$&1: Compatibility of $\to$ with $\vee$\\
                6&$\neg\psi\wedge(\phi\vee\neg\phi)\to\neg\psi\wedge(\psi\vee\neg\phi)$&5: Compatibility of $\to$ with $\wedge$\\
                7&$\neg\psi\wedge(\psi\vee\neg\phi)\to(\neg\psi\wedge\psi)\vee(\neg\psi\wedge\neg\phi)$&Distribution\\
                8&$\neg\psi\wedge\psi\to\neg\psi\wedge\neg\phi$&ECQ\\
                9&$\neg\psi\wedge\neg\phi\to\neg\psi\wedge\neg\phi$&Reflexivity\\
                10&$(\neg\psi\wedge\psi)\vee(\neg\psi\wedge\neg\phi)\to\neg\psi\wedge\neg\phi$&8,9: ``Two-to-one''\\
                11&$\neg\psi\wedge\neg\phi\to\neg\phi$&$\wedge$E\\
                12&$\neg\psi\to\neg\phi$&4,6,7,10,11: Transitivity
            \end{tabular}\\
        \end{proof}

        \begin{derivedrule}[Contraposition, negation eliminated]
            $$\neg\phi\to\neg\psi\tto \psi\to\phi.$$
        \end{derivedrule}

        \begin{proof}\;\\
            \begin{tabular}{lll}
                1&$\neg\phi\to\neg\psi$&Premise\\
                2&$\psi\to\neg\neg\psi$&DNI\\
                3&$\neg\neg\psi\to\neg\neg\phi$&1: Contraposition, negation introduced\\
                4&$\neg\neg\phi\to\phi$&DNE\\
                5&$\psi\to\phi$&2,3,4: Transitivity
            \end{tabular}\\
        \end{proof}

        \begin{derivedrule}[Prefixing]
            $$\phi\to\psi \tto (\chi\to\phi)\to(\chi\to\psi).$$
        \end{derivedrule}

        \begin{proof}\;\\
            \begin{tabular}{lll}
                1&$\phi\to\psi$&Premise\\
                2&$(\chi\to\phi)\to(\neg\phi\to\neg\chi)$&Contraposition\\
                3&$\neg\psi\to\neg\phi$&1: Contraposition\\
                4&$(\neg\phi\to\neg\chi)\to(\neg\psi\to\neg\chi)$&3: Suffixing\\
                5&$(\neg\psi\to\neg\chi)\to(\chi\to\psi)$&Contraposition\\
                6&$(\chi\to\phi)\to(\chi\to\psi)$&2,4,5: Transitivity
            \end{tabular}\\
        \end{proof}

        \begin{derivedrule}[RAA, another $\neg$I]
            $$\{\phi\to\psi,\phi\to\neg\psi\}\tto\neg\phi.$$
        \end{derivedrule}

        \begin{proof}\;\\
            \begin{tabular}{lll}
                1&$\phi\to\psi$&Premise\\
                2&$\phi\to\neg\psi$&Premise\\
                3&$\neg\psi\to\neg\phi$&1: Contraposition\\
                4&$\neg\neg\psi\to\neg\phi$&2: Contraposition\\
                5&$\neg\psi\vee\neg\neg\psi\to\neg\phi$&3,4: ``Two-to-one''\\
                6&$\neg\psi\vee\neg\neg\psi$&LEM\\
                7&$\neg\phi$&5,6: MP
            \end{tabular}\\
        \end{proof}

        \begin{derivedrule}[Definition of $\to$ by $\oto$]
            $$\phi\to\psi\otto \phi\wedge\psi\oto\phi.$$
        \end{derivedrule}

        \begin{proof}\;\\
            \begin{tabular}{lll}
                1&$\phi\to\psi$&Premise\\
                2&$\phi\to\phi$&Reflexivity\\
                3&$\phi\to\phi\wedge\psi$&1,2: ``One-to-two''\\
                4&$\phi\wedge\psi\to\phi$&$\wedge$E\\
                5&$\phi\wedge\psi\oto\phi$&3,4: $\oto$I\\
                & &\\
                1&$\phi\wedge\psi\oto\phi$&Premise\\
                2&$\phi\to\phi\wedge\psi$&1: $\oto$E\\
                3&$\phi\wedge\psi\to\psi$&$\wedge$E\\
                4&$\phi\to\psi$&2,3: Transitivity
            \end{tabular}\\
        \end{proof}

        Based on some of the above results, it is not difficult to prove the primitive rule $\{\phi\oto\psi,\chi\} \tto \chi[\phi\mapsto\psi]$ of replacement of the Hilbert-style IRL (cf. Subsection \ref{subsec:Logic IRL}) by induction.


        \section{Extension to first-order logic}
        \label{sec:Extension to first-order logic}

        Extending IRL to first-order logic (FOL) is similar to what is done classically.

        Let $\{x,y,...\}$, $\{s,t,...\}$, $\{P,Q,...\}$ and $\{\phi,\psi,...\}$ be the FOL variables, terms, predicates and formulas, respectively.

        Suppose that in a substitution $\phi\{t/x\}$ or a replacement $\phi[s\mapsto t]$, $s$ and $t$ are not bound variables in $\phi$ (if they are, rename them before substitution or replacement).

        The extension of the propositional logic IRL to first-order logic is defined as a FOL that contains:

        \begin{enumerate}
            \item All theorems and rules of IRL where propositional formulas are extended to be any FOL formulas.
            \item Axioms for equality:
            \begin{itemize}
                \item Reflexivity: $s=s$. Symmetry and transitivity can be derived.
                \item Replacement: $(s=t)\wedge \phi\to \phi[s\mapsto t]$.
            \end{itemize}
            \item Axioms for quantifiers:
            \begin{itemize}
                \item Universal generalization: $\phi(y)\to \forall x(\phi(x))$ if $y$ is arbitrary.
                \item Universal instantiation: $\forall x(\phi(x))\to \phi\{t/x\}$ where $t$ can be any term.
                \item ``One-to-many'': $\forall x(\phi\to\psi(x))\to (\phi\to \forall x(\psi(x)))$ where $x$ does not occur free in $\phi$. Its converse can be derived.
                \item  Duality: $\neg\forall x(\phi)\to \exists x(\neg\phi)$. Another one $\neg\exists x(\phi)\to \forall x(\neg\phi)$ can be derived.
                \item Existential instantiation: $\exists x(\phi(x))\to \phi\{t_0/x\}$ where $t_0$ is unknown, so it must not occur earlier in the proof and it also must not occur in the conclusion of the proof.
                \item Existential generalization (redundant): $\phi\{t_0/x\}\to \exists x(\phi(x))$ where $t_0$ is some (known) term.
                \item ``Many-to-one'' (redundant): $\forall x(\phi(x)\to\psi)\to (\exists x(\phi(x))\to \psi))$ where $x$ does not occur free in $\psi$. Its converse can be derived so is redundant as well.
            \end{itemize}
        \end{enumerate}

        The following are some derived results.

        \begin{theorem}[``One-to-many'', converse]
            $$(\phi\to \forall x(\psi(x))) \to \forall x(\phi\to\psi(x))$$
            where $x$ does not occur free in $\phi$.
        \end{theorem}

        \begin{proof}\;\\
            \begin{tabular}{lll}
                1&$\phi\to \forall x(\psi(x))$& Premise\\
                2&$\forall x(\psi(x))\to \psi(y)$& Universal instantiation ($y$ is arbitrary)\\
                3&$\phi\to\psi(y)$& 1,2: Transitivity\\
                4&$\forall x(\phi\to\psi(x))$& 3: Universal generalization (as $y$ is arbitrary from line 2)
            \end{tabular}\\
        \end{proof}

        \begin{theorem}[Duality, existential to universal quantifier]
            $$\neg\exists x(\phi) \to \forall x(\neg\phi).$$
        \end{theorem}

        \begin{proof}\;\\
            \begin{tabular}{lll}
                1&$\neg\exists x(\phi)$& Premise\\
                2&$\neg\forall x(\neg\phi)\to \exists x(\neg\neg\phi)$& Duality, universal to existential quantifier\\
                3&$\neg\forall x(\neg\phi)\to \exists x(\phi)$& 2: Double negation\\
                4&$\neg\neg\forall x(\neg\phi)$& 1,3: Modus tollens\\
                5&$\forall x(\neg\phi)$& 4: Double negation
            \end{tabular}\\
        \end{proof}

        \begin{theorem}[Existential generalization]
            $$\phi\{t_0/x\}\to \exists x(\phi(x))$$
            where $t_0$ is some (known) term.
        \end{theorem}

        \begin{proof}\;\\
            \begin{tabular}{lll}
                1&$\phi\{t_0/x\}$& Premise\\
                2&$\forall x(\neg\phi(x))\to \neg\phi\{t_0/x\}$& Universal instantiation\\
                3&$\neg\forall x(\neg\phi(x))$& 1,2: Modus tollens\\
                4&$\exists x(\neg\neg\phi(x))$& 3: Duality\\
                5&$\exists x(\phi(x))$& 4: Double negation
            \end{tabular}\\
        \end{proof}

        \begin{theorem}[``Many-to-one'']
            $$\forall x(\phi(x)\to\psi)\to (\exists x(\phi(x))\to \psi))$$
            where $x$ does not occur free in $\psi$.
        \end{theorem}

        \begin{proof}\;\\
            \begin{tabular}{lll}
                1&$\forall x(\phi(x)\to\psi)$& Premise\\
                2&$\phi(y)\to\psi$& 1: Universal instantiation ($y$ is arbitrary)\\
                3&$\neg\psi\to\neg\phi(y)$& 2: Contraposition\\
                4&$\forall x(\neg\psi\to\neg\phi(x))$& 3: Universal generalization (as $y$ is arbitrary from line 2)\\
                5&$\neg\psi\to\forall x(\neg\phi(x))$& 4: ``One-to-many''\\
                6&$\neg\forall x(\neg\phi(x))\to\psi$& 5: Contraposition, Double negation\\
                7&$\exists x(\phi(x))\to\psi$& 6: Duality
            \end{tabular}\\
        \end{proof}

        \begin{theorem}[``Many-to-one'', converse]
            $$(\exists x(\phi(x))\to \psi)\to \forall x(\phi(x)\to\psi)$$
            where $x$ does not occur free in $\psi$.
        \end{theorem}

        \begin{proof}\;\\
            \begin{tabular}{lll}
                1&$\exists x(\phi(x))\to\psi$& Premise\\
                2&$\neg\psi\to\neg\exists x(\phi(x))$& 1: Contraposition\\
                3&$\neg\psi\to\forall x(\neg\phi(x))$& 2: Duality\\
                4&$\forall x(\neg\psi\to\neg\phi(x))$& 3: ``One-to-many'', converse\\
                5&$\neg\psi\to\neg\phi(y)$& 4: Universal instantiation ($y$ is arbitrary)\\
                6&$\phi(y)\to\psi$& 5: Contraposition\\
                7&$\forall x(\phi(x)\to\psi)$& 6: Universal generalization (as $y$ is arbitrary)
            \end{tabular}\\
        \end{proof}

        \begin{theorem}[Universal conjunction]
            $$\forall x(\phi(x)\wedge\psi(x))\oto \forall x(\phi(x)\wedge \forall x(\psi(x).$$
        \end{theorem}

        \begin{proof}\;\\
            \begin{tabular}{lll}
                1&$\forall x(\phi(x)\wedge\psi(x))$& Premise\\
                2&$\phi(y)\wedge\psi(y)$& 1: Universal instantiation ($y$ is arbitrary)\\
                3&$\phi(y)$& 2: Conjunction elimination\\
                4&$\forall x(\phi(x))$& 3: Universal generalization (as $y$ is arbitrary)\\
                5&$\forall x(\psi(x))$& 2--4: Similarly\\
                6&$\forall x(\phi(x)\wedge \forall x(\psi(x)$& 4,5: Conjunction introduction
            \end{tabular}\\
            The converse can be proved in the reverse order.
        \end{proof}

        \begin{theorem}[Existential disjunction]
            $$\exists x(\phi(x)\vee\psi(x))\oto \exists x(\phi(x)\vee\exists x(\psi(x)).$$
        \end{theorem}

        \begin{proof}\;\\
            \begin{tabular}{lll}
                1&$\exists x(\phi(x)\vee\psi(x))$& Premise\\
                2&$\phi(a)\vee\psi(a)$& 1: Existential instantiation\\
                3&$\phi(a)\to\exists x(\phi(x))$& 2: Existential generalization\\
                4&$\phi(a)\to\exists x(\phi(x))\vee\exists x(\psi(x))$& 3: Weakening\\
                5&$\psi(a)\to\exists x(\phi(x))\vee\exists x(\psi(x))$& 2--4: Similarly\\
                6&$\exists x(\phi(x))\vee\exists x(\psi(x))$& 2,4,5: Proof by exhaustion\\
                & &\\
                1 &$\exists x(\phi(x))\vee\exists x(\psi(x))$& Premise\\
                2 &$\exists x(\phi(x))\to \phi(a)$& Existential instantiation\\
                3 &$\exists x(\phi(x))\to \phi(a)\vee\psi(a)$& 2: Weakening\\
                4 &$\phi(a)\vee\psi(a)\to \exists x(\phi(x)\vee\psi(x))$& Existential generalization\\
                5 &$\exists x(\phi(x))\to \exists x(\phi(x)\vee\psi(x))$& 3--4: Transitivity\\
                6 &$\exists x(\psi(x))\to \exists x(\phi(x)\vee\psi(x))$& 2--5: Similarly\\
                7 &$\exists x(\phi(x)\vee\psi(x))$& 1,5,6: Proof by exhaustion
            \end{tabular}\\
        \end{proof}

        \begin{theorem}[Quantifier switching]
            $$\exists x\forall y(\phi(x,y))\to \forall y\exists x(\phi(x,y)).$$
        \end{theorem}

        \begin{proof}\;\\
            \begin{tabular}{lll}
                1&\qquad$\exists x\forall y(\phi(x,y))$& Assumption\\
                2&\qquad$\neg\forall y\exists x(\phi(x,y))$& Assumption\\
                3&\qquad$\forall y(\phi(a,y))$& 1: Existential instantiation\\
                4&\qquad$\exists y\forall x(\neg\phi(x,y))$& 2: Duality\\
                5&\qquad$\forall x(\neg\phi(x,b))$& 4: Existential instantiation\\
                6&\qquad$\phi(a,b)$& 3: Universal instantiation\\
                7&\qquad$\neg\phi(a,b)$& 5: Universal instantiation\\
                8&\qquad$\bot$& 6,7: LNC\\
                9&$\exists x\forall y(\phi(x,y))\wedge\neg\forall y\exists x(\phi(x,y))\to \bot$& 1--8: Conditional proof\\
                10&$\exists x\forall y(\phi(x,y))\to \forall y\exists x(\phi(x,y))$& 9: Implication and inconsistency
            \end{tabular}\\
        \end{proof}

        \begin{metarule}[Equivalence of implication theorem and disjunction theorem in FOL]
            \label{metarule:Equivalence of implication theorem and disjunction theorem in FOL}
            $$(\vdash\forall x\,(\phi(x)\to\psi(x))) \otto (\vdash\forall x\,(\neg\phi(x)\vee\psi(x))).$$
        \end{metarule}

        \begin{proof}\;\\
            \begin{tabular}{lll}
                1&$\vdash\forall x\,(\phi(x)\to\psi(x))$&Premise\\
                2&\qquad$\forall x\,(\phi(x)\to\psi(x))$&1 (assumed theorem)\\
                3&\qquad$\phi(y)\to\psi(y)$&2: Universal instantiation ($y$ is arbitrary)\\
                4&$\vdash\phi(y)\to\psi(y)$&2--3: Unconditional proof\\
                5&$\vdash\neg\phi(y)\vee\psi(y)$&4: Meta Rule \ref{metarule:Equivalence of implication theorem and disjunction theorem} of IRL\\
                6&\qquad$\neg\phi(y)\vee\psi(y)$&5 (theorem indicated in an outer proof)\\
                7&\qquad$\forall x\,(\neg\phi(x)\vee\psi(x))$&6: Universal generalization (as $y$ is arbitrary)\\
                8&$\vdash\forall x\,(\neg\phi(x)\vee\psi(x))$&6--7:  Unconditional proof
            \end{tabular}\\
            The converse can be proved in the reverse order.
        \end{proof}

        \subsubsection*{Bounded quantifiers}

        It is convenient if including \emph{bounded quantifiers} or \emph{restricted quantifiers} in a first-order language.

        \begin{definition}[Bounded quantifiers]
            \label{def:Bounded quantifiers}
            Let $\chi$ be a unary formula, $\phi$ be any formula.
            Let $X=\{x\mid \chi(x)\}$.
            Then $x\in X$ if and only if $\chi(x)$.
            The \emph{bounded quantifiers} $\exists x\in X$ and $\forall x\in X$ are defined by the rules:
            \begin{itemize}
                \item $\exists x\in X \,(\phi) \otto \exists x\, (\chi(x)\wedge\phi)$.
                \item $\forall x\in X\,(\phi) = \neg\exists x\in X\, (\neg\phi) \otto \neg\exists x\,(\chi(x)\wedge\neg\phi) \otto \forall x\,\neg(\chi(x)\wedge\neg\phi)$, i.e.,  $\forall x\in X\,(\phi) \otto \forall x\,(\neg\chi(x)\vee\phi)$.
            \end{itemize}
        \end{definition}

        \begin{remark}
            The universal bounded quantifier $\forall x\in X$ can not be defined by $\forall x\in X\,(\phi)\otto \forall x\, (\chi(x)\to\phi)$ since $\chi(x)\to\phi$ and $\neg\chi(x)\vee\phi$ are not equivalent in IRL -- however, if $\forall x\in X\,(\phi)$ represents a theorem, the two expressions are equivalent according to the above Meta Rule \ref{metarule:Equivalence of implication theorem and disjunction theorem in FOL}.
        \end{remark}


        \section{Discussion}
        \label{sec:Discussion}

        \subsubsection*{Normal situation}

        An arbitrary proposition can be tautological, contradictory, or contingent. Most generally, there are possibly infinitely many T-set values for a contingency ($0<x<1$) whereas there is only one value for a tautology ($x=1$) and one value for a contradiction ($x=0$). This is just similar to the unit interval $[0,1]$ (no matter it is different as totally ordered) where the number of interior points is uncountable while the number of boundary points is only two. In this sense, it is the \emph{normal situation} for a proposition being contingent rather than tautological or contradictory. Subsequently in this section, when saying ``in the normal situation'' we mean propositions mentioned are all contingent.

        \subsubsection*{Interpretation of Generalized ECQ of implication}

        In IRL, the Generalized ECQ $\psi\to(\phi\to\psi)$ does not hold while the Generalized ECQ \emph{of implication} $(\phi\to\psi)\to(\chi\to(\phi\to\psi))$ does.

        Generalized ECQ $\psi\to(\phi\to\psi)$ means that, for any propositions $\phi$ and $\psi$, there is always an implication from $\phi$ to $\psi$ in the normal situation (where $\psi$ can be true so $\phi\to\psi$ can be true). This is not consistent with the original meaning of implication as it is not ensured for a ``certain mechanism'' of implication in this situation.

        Whereas Generalized ECQ of implication $(\phi\to\psi)\to(\chi\to(\phi\to\psi))$ means that if there exists the ``certain mechanism'' of implication from $\phi$ to $\psi$, then this mechanism is irrelevant to any proposition $\chi$. This does not conflict with the original meaning of implication at all.

        The ``worst'' case for Generalized ECQ of implication is when $\chi=\neg(\phi\to\psi)$. From IRL theorems $(\neg\phi\to\phi)\oto(\top\to\phi)$ and $(\top\to\phi)\to\phi$, we have $(\top\to(\phi\to\psi))\oto(\phi\to\psi)$ in the worst case. This can be understood as a way to interpret $\top\to(\phi\to\psi)$ as $\phi\to\psi$. In contrast, since the general form $(\top\to\phi)\oto\phi$ does not hold (as $\phi\to(\top\to\phi)$ does not hold), $\top\to\phi$ can not be interpreted as $\phi$ generally.

        \subsubsection*{Interpretation of ECQ}

        A simple conditional proof of ECQ $\phi\wedge\neg\phi\to\psi$: \\
        \begin{tabular}{lll}
            1 & $\phi$ & Premise\\
            2 & $\neg\phi$ & Premise\\
            3 & $\phi\vee\psi$ & 1: \emph{Disjunction introduction}\\
            4 & $\psi$ & 2,3: \emph{Disjunctive syllogism}\\
        \end{tabular}

        Hence ECQ is necessary for disjunction introduction and disjunctive syllogism which do not conflict to the original meaning of implication. Therefore, ECQ contributes to the goal of strength for the system (cf. Subsection \ref{subsec:Construction of the right logic}).

        Furthermore, ECQ $\bot\to\phi$ and its dual $\phi\to\top$ can be interpreted to intuition as ``if an always-false thing is (was) true, then anything is (would be) true'' and ``in any cases, an always-true thing is true'', respectively, as similarly mentioned, e.g., in \citep{ceniza}.

        ECQ is an abnormal situation, for which a ``standard (set of rules for) interpretation'' may be explored for practical purpose. For example, given a conditional $P\to Q$: if $P$ is a contradiction (e.g. it is well known that $P$ is false), then $P\to Q$ can be interpreted as ``if the false thing like $P$ were true, anything including $Q$ should be true''; if $Q$ is a tautology (e.g. it is a well known fact), then $P\to Q$ can be interpreted as ``in any cases including $P$, an always true thing like $Q$ is true''.


        \section{Concluding remarks}
        \label{sec:ConcludingRemarks}

        By the original meaning of implication, the formula $\neg\phi\vee\psi\to(\phi\to\psi)$ is generally invalid, so the material implication $(\phi\to\psi)\oto\neg\phi\vee\psi$ is not valid in general. Defining implication according to its original meaning is the key to avoid invalid results. The logic IRL of this work has the expected property that the invalid class of $\neg\phi\vee\psi\to(\phi\to\psi)$ is removed while fundamental laws like LEM and double negation, LNC and ECQ, conjunction elimination and disjunction introduction, and transitivity of implication and disjunctive syllogism, are all retained. Based on information provided with IRL, more points are summarized as follows.

        \begin{enumerate}

            \item Minimal functionally complete operator sets. (a) $\phi\to\psi\otto\neg\phi\vee\psi$ is invalid, so $\to$ and $\oto$ cannot be defined by $\{\neg,\wedge,\vee\}$, meaning that one of $\to$ and $\oto$ must be primitive; (b) $\neg\phi\otto \phi\to\bot$ is invalid, so $\neg$ cannot be defined by $\{\bot,\to\}$, meaning that $\neg$ must be primitive; (c) for the same reason as in (a), one of $\wedge$ and $\vee$ must be primitive. Thus, a functionally complete operator set contains at least three members, e.g. $\{\neg,\wedge,\to\}$, $\{\neg,\vee,\to\}$, $\{\neg,\wedge,\oto\}$, $\{\neg,\vee,\oto\}$.

            \item Truth value of implication. Since implication is not truth functional, valuation of $\phi\to\psi$ can not be generally done by valuation of $\phi$ and $\psi$, but it can be ``estimated'' via the procedure: (a) $\phi\to\psi$ is a tautology, if and only if $\neg\phi\vee\psi$ is a tautology; (b) $\phi\to\psi$ is a contradiction, if $\neg\phi\vee\psi$ is a contradiction (i.e. $\phi$ is a tautology and $\psi$ is a contradiction); (c) the T-set of $\phi\to\psi$ is a subset (possibly empty) of the T-set of $\neg\phi\vee\psi$, if $\neg\phi\vee\psi$ is a contingency. T-sets of primitive implications, like any primitive propositions, are given (e.g. by physical mechanisms) rather than derived.

            \item Classical logic is sound and complete with regard to its own semantics. In the sense that it contains theorems that are not generally valid with respect to semantics based on the original meaning of implication, classical logic is less ``sound'' than ``complete'' (actually it is post-complete), while some weaker systems are in the opposite.

            \item The logic IRL of this work, in which the generally invalid classical theorems, including especially the class of $\neg\phi\vee\psi\to(\phi\to\psi)$, are removed, should be used instead of classical logic for the general case (the system contains any types of propositions including especially contingencies). IRL and classical logic are the same if and only if there is no (primitive) contingent propositions in a system.

            \item An interpreted formal system has no contingent statements. Precisely, for a system under a given interpretation, every sentence is simply either true or false, no longer relevant to ``tautological, contradictory or contingent''. It is common that a theory has an \emph{intended interpretation} or a \emph{standard model} or \emph{usual model}. For example, Tarski's axioms for Euclidean geometry has its standard model of Euclidean plane geometry, and the standard model for Peano axioms and Robinson arithmetic is the natural numbers. In such contexts when a system contains no contingencies, classical logic can be used safely.
        \end{enumerate}


    \end{document}